
%
%
%
%
%
%
\magnification=\magstephalf      
%
%
\vsize=7.5truein                 
\hsize=5.2truein                 
\newskip\stdskip                 
\stdskip=6pt plus3pt minus3pt    
\medskipamount=\stdskip          
\parindent=0pt                   
\parskip=\stdskip                
\abovedisplayskip=\stdskip       
\belowdisplayskip=\stdskip       
\mathsurround=0.75pt             
\overfullrule=0pt                
%
%
\def\ppar{\par\goodbreak\vskip 8pt plus 4pt minus 4pt}     
%
%
\def\stdspace{\hskip 0.75em plus 0.15em\ignorespaces}
\let\qua\stdspace 
%
%
%
%
%
%
%
\def\hexnumber#1{\ifcase#1 0\or 1\or 2\or 3\or 4\or 5\or 6\or 7\or 8\or
 9\or A\or B\or C\or D\or E\or F\fi}
%
%
\font\thirtnmsa=msam10 scaled 1315    
\font\tenmsa=msam10          \font\ninemsa=msam9
\font\sevenmsa=msam7         \font\sixmsa=msam6
\font\fivemsa=msam5
%
%
\newfam\msafam                  \textfont\msafam=\tenmsa
\scriptfont\msafam=\sevenmsa    \scriptscriptfont\msafam=\fivemsa
\edef\hexa{\hexnumber\msafam}        
\def\msa{\fam\msafam\tenmsa}         
%
%
\font\thirtnmsb=msbm10 scaled 1315   
\font\tenmsb=msbm10      \font\ninemsb=msbm9
\font\sevenmsb=msbm7     \font\sixmsb=msbm6
\font\fivemsb=msbm5
%
\newfam\msbfam                   \textfont\msbfam=\tenmsb       
\scriptfont\msbfam=\sevenmsb     \scriptscriptfont\msbfam=\fivemsb
\edef\hexb{\hexnumber\msbfam}    
\def\msb{\fam\msbfam\tenmsb}     
%
%
\font\thirtneufm=eufm10 scaled 1315   
\font\teneufm=eufm10                 \font\nineeufm=eufm9
\font\seveneufm=eufm7                \font\sixeufm=eufm6
\font\fiveeufm=eufm5
%
\newfam\eufmfam                    \textfont\eufmfam=\teneufm
\scriptfont\eufmfam=\seveneufm     \scriptscriptfont\eufmfam=\fiveeufm
\edef\hexf{\hexnumber\eufmfam}      
\def\frak{\fam\eufmfam\teneufm}     
%
%
%
\font\thirtnrm=cmr10 scaled 1315    
\font\ninerm=cmr9                   \font\sixrm=cmr6   
%
\font\thirtni=cmmi10 scaled 1315    
\font\ninei=cmmi9                   \font\sixi=cmmi6  
%
\font\thirtnsy=cmsy10 scaled 1315   
\font\ninesy=cmsy9                  \font\sixsy=cmsy6  
%
\font\thirtnbf=cmbx10 scaled 1315   
\font\ninebf=cmbx9                  \font\sixbf=cmbx6  
%
%
\font\thirtnex=cmex10 scaled 1315   
\font\nineex=cmex9                  
%
%
\font\thirtnit=cmti10 scaled 1315  
\font\nineit=cmti9                  
%
\font\thirtnsl=cmsl10 scaled 1315  
\font\ninesl=cmsl9                  
%
\font\thirtntt=cmtt10 scaled 1315  
\font\ninett=cmtt9                  
%
%
%
%
\def\small{%
%
%
\textfont0=\ninerm \scriptfont0=\sixrm \scriptscriptfont0=\fiverm
\def\rm{\fam0\ninerm}
%
%
\textfont1=\ninei \scriptfont1=\sixi \scriptscriptfont1=\fivei
%
%
\textfont2=\ninesy \scriptfont2=\sixsy \scriptscriptfont2=\fivesy
%
%
\textfont3=\nineex \scriptfont3=\nineex \scriptscriptfont3=\nineex
%
%
\textfont\bffam=\ninebf \scriptfont\bffam=\sixbf
\scriptscriptfont\bffam=\fivebf \def\bf{\fam\bffam\ninebf}%
%
%
\textfont\itfam=\nineit \def\it{\fam\itfam\nineit}%
\textfont\slfam=\ninesl \def\sl{\fam\slfam\ninesl}%
\textfont\ttfam=\ninett \def\tt{\fam\ttfam\ninett}%
%
%
%
\textfont\msafam=\ninemsa \scriptfont\msafam=\sixmsa
\scriptscriptfont\msafam=\fivemsa \def\msa{\fam\msafam\ninemsa}%
%
%
\textfont\msbfam=\ninemsb \scriptfont\msbfam=\sixmsb
\scriptscriptfont\msbfam=\fivemsb \def\msb{\fam\msbfam\ninemsb}%
%
%
\textfont\eufmfam=\nineeufm  \scriptfont\eufmfam=\sixeufm
\scriptscriptfont\eufmfam=\fiveeufm \def\frak{\fam\eufmfam\nineeufm}%
%
%
%
\normalbaselineskip=11pt%
\setbox\strutbox=\hbox{\vrule height8pt depth3pt width0pt}%
%
%
\normalbaselines\rm
%
%
\stdskip=4pt plus2pt minus2pt    
\medskipamount=\stdskip          
\parskip=\stdskip                
\abovedisplayskip=\stdskip       
\belowdisplayskip=\stdskip       
\def\ppar{\par\goodbreak\vskip 6pt plus 3pt minus 3pt}%
%
%
\def\section##1{\global\advance\sectionnumber by 1
\vskip-\lastskip\penalty-800\vskip 20pt plus10pt minus5pt 
\egroup{\bf\number\sectionnumber\quad##1}\bgroup\small         
\vskip 6pt plus3pt minus3pt
\nobreak\resultnumber=1}
}    
%
\def\beginsmall{\bgroup\small}
\let\endsmall\egroup
%
%
%
%
\def\large{%
\textfont0=\thirtnrm \scriptfont0=\ninerm \scriptscriptfont0=\sevenrm
\def\rm{\fam0\thirtnrm}%
\textfont1=\thirtni \scriptfont1=\ninei \scriptscriptfont1=\seveni
\textfont2=\thirtnsy \scriptfont2=\ninesy \scriptscriptfont2=\sevensy
\textfont3=\thirtnex \scriptfont3=\thirtnex \scriptscriptfont3=\thirtnex
\textfont\bffam=\thirtnbf \scriptfont\bffam=\ninebf
\scriptscriptfont\bffam=\sevenbf \def\bf{\fam\bffam\thirtnbf}%
\textfont\itfam=\thirtnit \def\it{\fam\itfam\thirtnit}%
\textfont\slfam=\thirtnsl \def\sl{\fam\slfam\thirtnsl}%
\textfont\ttfam=\thirtntt \def\tt{\fam\ttfam\thirtntt}%
\textfont\msafam=\thirtnmsa \scriptfont\msafam=\ninemsa
\scriptscriptfont\msafam=\sevenmsa \def\msa{\fam\msafam\thirtnmsa}%
\textfont\msbfam=\thirtnmsb \scriptfont\msbfam=\ninemsb
\scriptscriptfont\msbfam=\sevenmsb \def\msb{\fam\msbfam\thirtnmsb}%
\textfont\eufmfam=\thirtneufm  \scriptfont\eufmfam=\nineeufm
\scriptscriptfont\eufmfam=\seveneufm \def\frak{\fam\eufmfam\teneufm}%
\normalbaselineskip=16pt%
\setbox\strutbox=\hbox{\vrule height11.5pt depth4.5pt width0pt}%
\normalbaselines\rm}%
\let\Large\large   
%
\def\Bbb#1{{\msb#1}}

%

\def\re{\Bbb R}
%
\mathchardef\plussquare="0\hexa01
\mathchardef\nge="3\hexb0B
\mathchardef\maltesecross="0\hexa7A
\mathchardef\del="0\hexf01
%
%
%
%
\font\sc=cmcsc10
%
%
%
%
\def\sqr#1#2{{\vcenter{\vbox{\hrule  height.#2truept
	\hbox{\vrule width.#2truept height#1truept 
	\kern#1truept \vrule width.#2truept}
	\hrule height.#2truept}}}}
\def\sq{\sqr55}    
%
%
%
%
\newcount\sectionnumber            
\newcount\resultnumber             
\sectionnumber=0\resultnumber=1    
%
%
%
\def\section#1{\global\advance\sectionnumber by 1
\xdef\nextkey{\number\sectionnumber}
\vskip-\lastskip\penalty-800\vskip 20pt plus10pt minus5pt 
{\large\bf\number\sectionnumber\quad#1}         
\vskip 8pt plus4pt minus4pt
\nobreak\resultnumber=1}      
%
%
%
%
%
\def\sh#1{\vskip-\lastskip\ppar{\bf #1}\par\nobreak\medskip}         
%
%
%
%

%
\def\proc#1{\xdef\nextkey{\number\sectionnumber.\number\resultnumber}%
\vskip-\lastskip\ppar\bf%
\noindent#1\ \number\sectionnumber.\number\resultnumber
\stdspace\sl\global\advance\resultnumber by 1\ignorespaces}
\def\endproc{\rm\ppar} 
%
%
\def\prf{\vskip-\lastskip\ppar\noindent{\bf Proof}%
\stdspace\rm}                            
\def\qed{\hfill$\sq$\par\goodbreak\rm}   
\def\endprf{\unskip\stdspace\hbox{}
\hfill$\sq$\par\medskip}                 
%
%
%
%
%
%
%
%
\def\proclaim#1{\vskip-\lastskip\ppar\bf%
\noindent#1\stdspace\sl\ignorespaces} 

%
%
%
%
\def\rk#1{\vskip-\lastskip\ppar{\bf #1}\stdspace\ignorespaces}                

%
%
%
%
%
%
\def\label{\xdef\nextkey{\number\sectionnumber.\number\resultnumber}%
\number\sectionnumber.\number\resultnumber
\global\advance\resultnumber by 1}
%
%
%
%
%
%
%
%
%
%
%
%
%
%
%
%
\newcount\refnumber              
\refnumber=1                     
\long\def\reflist#1\endreflist{%
\long\def\thereflist{#1}{\def\refkey##1##2\par{\xdef##1{\number\refnumber}%
\global\advance\refnumber by 1}%
\def\key##1##2\par{\expandafter\xdef%
\csname##1\endcsname{\number\refnumber}%
\global\advance\refnumber by 1}#1\par}}
\long\def\references{%
\penalty-800\vskip-\lastskip\vskip 15pt plus10pt minus5pt 
{\large\bf References}\ppar 
{\leftskip=25pt\frenchspacing    
\small\parskip=3pt plus2pt       
\def\refkey##1##2\par{\noindent  
\llap{[##1]\stdspace}\ignorespaces##2\par}         
\def\key##1##2\par{\noindent  
\llap{[\ref{##1}]\stdspace}\ignorespaces##2\par}  
\def\,{\thinspace}\thereflist\par}}
%
%
%
\newcount\footnotenumber         
\footnotenumber=1                
\def\fnote#1{\xdef\nextkey{\number\footnotenumber}%
{\small\ifnum\footnotenumber>9\parindent=14pt%
\else\parindent=10pt\fi\footnote{$^{\number\footnotenumber}$}%
{\hglue-5pt#1}\global\advance\footnotenumber by 1}}
%
%
%
%
%
%
%
\newcount\figurenumber          
\figurenumber=1                 
\def\caption#1{\xdef\nextkey{\number\figurenumber}%
\cl{\small Figure \number\figurenumber: #1}%
\global\advance\figurenumber by 1}
\def\figurelabel{\xdef\nextkey{\number\figurenumber}%
\cl{\small Figure \number\figurenumber}%
\global\advance\figurenumber by 1}
\long\def\figure#1\endfigure{{\xdef\nextkey{\number\figurenumber}%
\let\captiontext\relax\def\caption##1{\xdef\captiontext{##1}}%
\midinsert\cl{\ignorespaces#1\unskip\unskip\unskip\unskip}\vglue6pt\cl{\small 
Figure \number\figurenumber\ifx\captiontext\relax\else: \captiontext
\fi}\endinsert\global\advance\figurenumber by 1}}
%
%
%
%
%
%
%
\def\nextkey{??}   
%
\def\key#1{\expandafter\xdef\csname #1\endcsname{\nextkey}}
\def\ref#1{\expandafter\ifx\csname #1\endcsname\relax
\immediate\write16{Reference {#1} undefined}??\else
\csname #1\endcsname\fi}
%
%
%
%
%
%
%
\newread\gtinfile
\newwrite\gtreffile
\def\useforwardrefs{
\openin\gtinfile\jobname.ref
\ifeof\gtinfile
\closein\gtinfile
\immediate\write16{No file \jobname.ref}
\else
\closein\gtinfile
\input \jobname.ref
\fi
\immediate\openout\gtreffile \jobname.ref
%
%
\def\key##1{{\def\\{\noexpand}%
\expandafter\xdef\csname ##1\endcsname{\nextkey}%
\immediate\write\gtreffile{\\\expandafter\\\def\\\csname ##1\\\endcsname%
{\nextkey}}}}
%
%
\long\def\reflist##1\endreflist{%
\long\def\thereflist{##1}{\def\refkey####1####2\par{\xdef####1{%
\number\refnumber}{\def\\{\noexpand}\immediate\write\gtreffile
{\\\def\\####1{\number\refnumber}}}\global\advance\refnumber by 1}%
\def\key####1####2\par{\expandafter\xdef%
\csname####1\endcsname{\number\refnumber}%
{\def\\{\noexpand}\immediate\write\gtreffile
{\\\expandafter\\\def\\\csname ####1\\\endcsname{\number\refnumber}}}
\global\advance\refnumber by 1}##1\par}}
\long\def\biblio##1\endbiblio{\reflist##1\endreflist\references}%
%
%
\def\numkey##1{{\def\\{\noexpand}%
\xdef##1{\number\sectionnumber.\number\resultnumber}
\immediate\write\gtreffile{\\\def\\##1%
{\number\sectionnumber.\number\resultnumber}}}}
\def\seckey##1{{\def\\{\noexpand}\xdef##1{\number\sectionnumber}
\immediate\write\gtreffile{\\\def\\##1{\number\sectionnumber}}}}
\def\figkey##1{\xdef##1{\number\figurenumber}%
{\def\\{\noexpand}\immediate\write\gtreffile%
{\\\def\\##1{\number\figurenumber}}}
\number\figurenumber\global\advance\figurenumber by 1}
}   
%
%
%
%
\def\figkey#1{\xdef#1{\number\figurenumber}%
\number\figurenumber\global\advance\figurenumber by 1}
\def\fig#1#2\endfig{%
\midinsert\cl{#2}\vglue6pt\cl{\small Figure #1}\endinsert}
\def\newfig{\number\figurenumber\global\advance\figurenumber by 1}
\def\numkey#1{\xdef#1{\number\sectionnumber.\number\resultnumber}}
\def\seckey#1{\xdef#1{\number\sectionnumber}}
%
%
%
%
%
%
%
%
%
\def\verb{\catcode`\"=\active}       
\def\brev{\catcode`\"=12}            
\brev                                
\verb                                
{\obeyspaces\gdef {\ }}              
{\catcode`\`=\active\gdef`{\relax\lq}}
\def"{%
\begingroup\baselineskip=12pt\def\par{\leavevmode\endgraf}%
\tt\obeylines\obeyspaces\parskip=0pt\parindent=0pt%
\catcode`\$=12\catcode`\&=12\catcode`\^=12\catcode`\#=12%
\catcode`\_=12\catcode`\~=12%
\catcode`\{=12\catcode`\}=12\catcode`\%=12\catcode`\\=12%
\catcode`\`=\active\let"\endgroup}
\brev      
%
%
%
%
%
%
\def\items{\par\leftskip = 25pt}           
\def\enditems{\par\leftskip = 0pt}         
\def\item#1{\par\leavevmode\llap{#1\stdspace}%
\ignorespaces}                             
%
%

%
%
\def\co{\colon\thinspace}    
\def\np{\vfil\eject}         
\def\nl{\hfil\break}         
\def\cl{\centerline}         
\def\gt{{\mathsurround=0pt\it $\cal G\mskip-2mu$eometry \&\ 
$\cal T\!\!$opology}}        
\def\agt{{\mathsurround=0pt\it$\cal A\mskip-.7mu$lgebraic \&\ 
$\cal G\mskip-2mu$eometric $\cal T\!\!$opology}}  
%
%
%

%
%
%
%
%
\def\title#1{\def\thetitle{#1}}

\def\author#1{\edef\previousauthors{\theauthors}
 \ifx\theauthors\relax\def\theauthors{#1}\else
 \def\theauthors{\previousauthors\par#1}\fi}

%
\def\address#1{\edef\previousaddresses{\theaddress}
 \ifx\theaddress\relax\def\theaddress{#1}\else
 \def\theaddress{\previousaddresses\par\vskip 2pt\par#1}\fi}
\def\secondaddress#1{\edef\previousaddresses{\theaddress}
 \ifx\theaddress\relax\def\theaddress{#1}\else
 \def\theaddress{\previousaddresses\par{\rm and}\par#1}\fi}   

\def\email#1{\edef\previousemails{\theemail}
 \ifx\theemail\relax\def\theemail{#1}\else
 \def\theemail{\previousemails\hskip 0.75em\relax#1}\fi}
\def\secondemail#1{\edef\previousemails{\theemail}
 \ifx\theemail\relax\def\theemail{#1}\else
 \def\theemail{\previousemails\hskip 0.75em{\rm and}\hskip 0.75em
 \relax#1}\fi}
\def\url#1{\edef\previousurls{\theurl}
 \ifx\theurl\relax\def\theurl{#1}\else
 \def\theurl{\previousurls\hskip 0.75em\relax#1}\fi}
\def\secondurl#1{\edef\previousurls{\theurl}
 \ifx\theurl\relax\def\theurl{#1}\else
 \def\theurl{\previousurls\hskip 0.75em{\rm and}\hskip 0.75em
 \relax#1}\fi}
\long\def\abstract#1\endabstract{\long\def\theabstract{#1}}
\def\primaryclass#1{\def\theprimaryclass{#1}}
\def\secondaryclass#1{\def\thesecondaryclass{#1}}
\def\keywords#1{\def\thekeywords{#1}}
%
%
\let\\\par\let\thetitle\relax\let\theshorttitle\relax
\let\theauthors\relax\let\theshortauthors\relax
\let\theaddress\relax\let\theshortaddress\relax
\let\theemail\relax\let\theurl\relax
\let\theabstract\relax\let\theprimaryclass\relax
\let\thesecondaryclass\relax\let\thekeywords\relax
%
%
%
%
\long\def\maketitlepage{    

\vglue 0.2truein   

%
{\parskip=0pt\leftskip 0pt plus 1fil\def\\{\par\smallskip}{\large
\bf\thetitle}\par\medskip}   

\vglue 0.15truein 

%
{\parskip=0pt\leftskip 0pt plus 1fil\def\\{\par}{\sc\theauthors}
\par\medskip}%
 
\vglue 0.1truein 

%
{\small\parskip=0pt
{\leftskip 0pt plus 1fil\def\\{\par}{\sl\theaddress}\par}
\ifx\theemail\relax\else  
\vglue 5pt \def\\{\stdspace{\rm and}\stdspace} 
\cl{Email:\stdspace\tt\theemail}\fi
\ifx\theurl\relax\else    
\vglue 5pt \def\\{\stdspace{\rm and}\stdspace} 
\cl{URL:\stdspace\tt\theurl}\fi\par}

\vglue 7pt 

{\bf Abstract}

\vglue 5pt

\theabstract

\vglue 7pt 

{\bf AMS Classification numbers}\quad Primary:\quad \theprimaryclass\par

Secondary:\quad \thesecondaryclass

\vglue 5pt 

{\bf Keywords:}\quad \thekeywords

\np  

}    
%
%
\long\def\makeshorttitle{    


%
{\parskip=0pt\leftskip 0pt plus 1fil\def\\{\par\smallskip}{\large
\bf\thetitle}\par\medskip}   

\vglue 0.05truein 

%
{\parskip=0pt\leftskip 0pt plus 1fil\def\\{\par}{\sc\theauthors}
\par\medskip}%
 
\vglue 0.03truein 

%
{\small\parskip=0pt
{\leftskip 0pt plus 1fil\def\\{\par}{\sl\ifx\theshortaddress\relax
\theaddress\else\theshortaddress\fi}\par}
\ifx\theemail\relax\else  
\vglue 5pt \def\\{\stdspace{\rm and}\stdspace} 
\cl{Email:\stdspace\tt\theemail}\fi
\ifx\theurl\relax\else    
\vglue 5pt \def\\{\stdspace{\rm and}\stdspace} 
\cl{URL:\stdspace\tt\theurl}\fi\par}

\vglue 10pt 


{\small\leftskip 25pt\rightskip 25pt{\bf Abstract}\stdspace\theabstract

{\bf AMS Classification}\stdspace\theprimaryclass
\ifx\thesecondaryclass\relax\else; \thesecondaryclass\fi\par
{\bf Keywords}\stdspace \thekeywords\par}
\vglue 7pt
}    
\let\maketitle\makeshorttitle        
%
%

\def\volumenumber#1{\def\thevolumenumber{#1}}
\def\volumeyear#1{\def\thevolumeyear{#1}}
\def\pagenumbers#1#2{\def\startpage{#1}\def\finishpage{#2}}
\def\published#1{\def\publishdate{#1}}
\def\received#1{\def\receiveddate{#1}}
\def\revised#1{\def\reviseddate{#1}}
\let\reviseddate\relax
\volumenumber{X}
\volumeyear{20XX}
\pagenumbers{1}{XXX}
\published{XX Xxxember 20XX}

\long\def\makeagttitle{   
\agt\hfill      
\hbox to 60truept{\vbox to 0pt{\vglue -14truept{\bf [Logo here]}\vss}\hss}
\break
{\small Volume \thevolumenumber\ (\thevolumeyear)
\startpage--\finishpage\nl
Published: \publishdate}

\vglue .2truein

{\parskip=0pt\leftskip 0pt plus 1fil\def\\{\par\smallskip}{\large
\bf\thetitle}\par\medskip}   
\vglue 0.05truein 

%
{\parskip=0pt\leftskip 0pt plus 1fil\def\\{\par}{\sc\theauthors}
\par\medskip}%
 
\vglue 0.03truein 


{\small\leftskip 25truept\rightskip 25truept{\bf Abstract}\stdspace\theabstract

{\bf AMS Classification}\stdspace\theprimaryclass
\ifx\thesecondaryclass\relax\else; \thesecondaryclass\fi\par
{\bf Keywords}\stdspace \thekeywords\par}\vglue 7truept

}   


\def\Addresses{\bigskip
{\small \parskip 0pt \leftskip 0pt \rightskip 0pt plus 1fil \def\\{\par}
\sl\theaddress\par\medskip \rm Email:\stdspace\tt\theemail\par
\ifx\theurl\relax\else\smallskip \rm URL:\stdspace\tt\theurl\par\fi}}

\def\agtart{
\hoffset 14truemm
\voffset 31truemm
\font\phead=cmsl9 scaled 950
\font\pnum=cmbx10 scaled 913
\font\pfoot=cmsl9 scaled 950
\headline{\vbox to 0pt{\vskip -4.5mm\line{\small\phead\ifnum
\count0=\startpage ISSN numbers are printed here
\hfill {\pnum\folio}\else\ifodd\count0\def\\{ }%
\ifx\theshorttitle\relax\thetitle\else\theshorttitle\fi\hfill{\pnum\folio}
\else\def\\{ and }{\pnum\folio}\hfill\ifx\theshortauthors\relax\theauthors
\else\theshortauthors\fi\fi\fi}\vss}}
\footline{\vbox to 0pt{\vglue 0mm\line{\small\pfoot\ifnum\count0=\startpage
Copyright declaration is printed here\hfill\else
\agt, Volume \thevolumenumber\ (\thevolumeyear)\hfill\fi}\vss}}
\let\maketitle\makeagttitle\let\makeshorttitle\makeagttitle}

\input pictex          
\chardef\newinsCatAt\the\catcode `\@
\catcode `\@=11
%
%
%
\newskip\insertskipamount\newskip\inserthardskipamount
\insertskipamount 12pt plus2pt     
\inserthardskipamount 4pt          
\def\insertskip{\vskip\insertskipamount}
%
%
\newskip\LastSkip
\def\SaveLastSkip{\LastSkip\lastskip}
\def\RestoreLastSkip{\nobreak\vskip-\LastSkip\vskip\LastSkip}
%
%
\newcount\SplitTest
\def\SetSplitTest{\SplitTest\insertpenalties
  \insert\topins{\floatingpenalty1}%
  \advance\SplitTest-\insertpenalties}
%
%
\def\midinsert{\par
 \SaveLastSkip\penalty-150\SetSplitTest\RestoreLastSkip
 \ifnum\SplitTest=-1
  \@midfalse\p@gefalse\else\@midtrue\fi\@ins}
\def\@ins{\par\begingroup\setbox\z@\vbox\bgroup%
  \vglue\inserthardskipamount}
\def\endinsert{\egroup 
  \if@mid \dimen@\ht\z@ \advance\dimen@\dp\z@
    \advance\dimen@\insertskipamount
    \advance\dimen@\pagetotal\advance\dimen@-\pageshrink
    \ifdim\dimen@>\pagegoal\@midfalse\p@gefalse\fi\fi
  \if@mid%
    \ifdim\lastskip<\insertskipamount\removelastskip\insertskip\fi
    \nointerlineskip\box\z@\penalty-200\insertskip
  \else%
    \SaveLastSkip
    \insert\topins{\penalty100 
    \splittopskip\z@skip
    \splitmaxdepth\maxdimen \floatingpenalty\z@
    \ifp@ge \dimen@\dp\z@
    \vbox to\vsize{\unvbox\z@\kern-\dimen@}
    \else \box\z@\nobreak\insertskip\fi}
    \RestoreLastSkip
   \fi\endgroup}
%
\catcode `\@=\newinsCatAt
 


\def\ifplaintex{\expandafter\ifx\csname documentclass\endcsname\relax}


\ifplaintex 
\hoffset 14truemm
\voffset 31truemm
\else
\headsep 23pt
\footskip 35pt
\hoffset -4truemm
\voffset 12.5truemm
\fi

\expandafter\ifx\csname beginpicture\endcsname\relax
\expandafter\ifx\csname documentclass\endcsname\relax
\input pictex \else\font\fiverm=cmr5
\input prepictex \input pictex \input postpictex \fi\fi

\def\gt{{\mathsurround=0pt\it $\cal G\mskip-2mu$eometry \&\ 
$\cal T\!\!$opology}}        

\def\gtp{{\mathsurround=0pt\it $\cal G\mskip-2mu$eometry \&\ 
$\cal T\!\!$opology $\cal P\!$ublications}}  


\def\lognumber#1{\def\thelognumber{#1}}
\def\volumenumber#1{\def\thevolumenumber{#1}}
\def\papernumber#1{\def\thepapernumber{#1}}
\def\volumeyear#1{\def\thevolumeyear{#1}}

\def\pagenumbers#1#2{\def\startpage{#1}\def\finishpage{#2}}
\def\published#1{\def\publishdate{#1}}
\def\proposed#1{\def\theproposer{#1}}
\def\seconded#1{\def\theseconders{#1}}
\def\received#1{\def\receiveddate{#1}}
\def\revised#1{\def\reviseddate{#1}}
\def\accepted#1{\def\accepteddate{#1}}

\long\def\asciiabstract#1{\long\def\theasciiabstract{#1}}


\let\\\par\let\thelognumber\relax
\let\thevolumenumber\relax\let\thepapernumber\relax
\let\thevolumeyear\relax\let\thesamplenumber\relax\let\startpage\relax
\let\finishpage\relax\let\publishdate\relax\let\receiveddate\relax
\let\reviseddate\relax\let\accepteddate\relax\let\theasciititle\relax
\let\theasciiauthors\relax
\let\theasciiabstract\relax
\let\theasciiemail\relax\let\theshortauthors\relax\let\theshorttitle\relax

\long\def\maketitlep{   

\count0=\startpage

\gt\hfill      
\beginpicture
\setcoordinatesystem units <0.33truein, 0.33truein> point at 2.2 0.9
\setplotsymbol ({$\cal G$})
\plotsymbolspacing=9truept
\circulararc 315 degrees from 0 1 center at 0 0
\setplotsymbol ({$\cal T$})
\circulararc 315 degrees from 1 -1 center at 1 0
\endpicture
%
\break
{\small\ifx\thesamplenumber\relax 
Volume \else Sample
\fi\thevolumenumber\ (\thevolumeyear)
\startpage--\finishpage\nl
Published: \publishdate}
\vglue 0.5truein plus 0.4fil minus 0.1truein

{\parskip=0pt\leftskip 0pt plus 1fil\def\\{\par\smallskip}{\ifplaintex\large
\else\Large\fi\bf\thetitle}\par\medskip}   

\vglue 0pt plus 0.1fil 

{\parskip=0pt\leftskip 0pt plus 1fil\def\\{\par}{\sc\theauthors}
\par\medskip}

\vglue 0pt plus 0.1fil 

{\small\parskip=0pt\let\newline\\
{\leftskip 0pt plus 1fil\def\\{\par}{\sl\theaddress}\par}
\expandafter\ifx\theemail\relax    
\relax\else\vglue 5pt plus 0.02fil minus 2pt\def\\{\stdspace{\rm 
and}\stdspace} 
\cl{Email:\stdspace\tt\theemail}\fi
\ifx\theurl\relax                  
\relax\else\vglue 5pt plus 0.02fil minus 2pt\def\\{\stdspace{\rm 
and}\stdspace}
\cl{URL:\stdspace\tt\theurl}\fi\par}

\vglue 7pt plus 0.3fil minus 3pt

{\bf Abstract}
\vglue 5pt plus 0.1fil minus 2pt

\theabstract

\vglue 7pt plus 0.3fil minus 3pt

{\bf AMS Classification numbers}\quad Primary:\quad \theprimaryclass

Secondary:\quad \thesecondaryclass

\vglue 5pt plus 0.3fil minus 2pt

{\bf Keywords}\quad \thekeywords

\vglue 10pt plus 0.5fil minus 5pt

{\small  Proposed: \theproposer\hfill Received: \receiveddate\nl
Seconded: \theseconders\hfill 
\ifx\reviseddate\relax                         
Accepted: \accepteddate                        
\else
Revised: \reviseddate                          
\fi}
\eject
}       

\let\maketitlepage\maketitlep
\let\maketitle\maketitlepage


\font\phead=cmsl9 scaled 950
\font\lhead=cmsl9 scaled 1050
\font\pnum=cmbx10 scaled 913
\font\lnum=cmbx10 
\font\pfoot=cmsl9 scaled 950
\font\lfoot=cmsl9 scaled 1050
\ifplaintex
\headline{\vbox to 0pt{\vskip -4.5mm\line{\small\phead\ifnum
\count0=\startpage ISSN 1364-0380 (on line)
1465-3060 (printed) \hfill {\pnum\folio}\else\ifodd\count0\def\\{ }%
\ifx\theshorttitle\relax\thetitle\else\theshorttitle\fi\hfill{\pnum\folio}
\else\def\\{ and }{\pnum\folio}\hfill\ifx\theshortauthors\relax\theauthors
\else\theshortauthors\fi\fi\fi}\vss}}
\footline{\vbox to 0pt{\vglue 0mm\line{\small\pfoot\ifnum\count0=\startpage
\copyright\ \gtp\hfill\else
\gt, Volume \thevolumenumber\ (\thevolumeyear)\hfill\fi}\vss
}}
\else
\makeatletter
\def\@oddhead{{\small\lhead\ifnum\count0=\startpage ISSN 1364-0380 (on line)
1465-3060 (printed) \hfill {\lnum\number\count0}\else\ifodd\count0
\def\\{ }\ifx\theshorttitle\relax \thetitle \else\theshorttitle\fi\hfill
{\lnum\number\count0}\else\def\\{ and }{\lnum\number\count0}
\hfill\ifx\theshortauthors\relax 
\theauthors\else\theshortauthors\fi\fi\fi}}\def\@evenhead{\@oddhead}
\def\@oddfoot{\small\lfoot\ifnum\count0=\startpage\copyright\ \gtp\hfill\else
\gt, Volume \thevolumenumber\ (\thevolumeyear)\hfill\fi}
\def\@evenfoot{\@oddfoot}
\makeatother
\fi


\newwrite\gtoutfile
\long\gdef\makeheadfile{  
{\def\\{, }\def\s{ }
\immediate\openout\gtoutfile head.xxx
\immediate\write\gtoutfile{To: math@arxiv.org}
\immediate\write\gtoutfile{Subject: put or rep NNNNN:pppp}
\immediate\write\gtoutfile{--text follows this line--}
\immediate\write\gtoutfile{Proxy-for: \ifx\theasciiauthors\relax
\theauthors\else\theasciiauthors\fi\s<\ifx\theasciiemail\relax\theemail\else\theasciiemail\fi>}
\immediate\write\gtoutfile{\noexpand\\}
\immediate\write\gtoutfile{Authors: \ifx\theasciiauthors\relax
\theauthors\else\theasciiauthors\fi}
\immediate\write\gtoutfile{Title: \ifx\theasciititle\relax
\thetitle\else\theasciititle\fi}
\immediate\write\gtoutfile{Subj-class: GT or SG or MG etc}
\immediate\write\gtoutfile{MSC-class: \theprimaryclass\ifx\thesecondaryclass\relax\else, \thesecondaryclass\fi}
\immediate\write\gtoutfile{Journal-ref: Geom. Topol. \thevolumenumber
(\thevolumeyear) \startpage-\finishpage}
\immediate\write\gtoutfile{Comments: Published by Geometry and Topology at}
\immediate\write\gtoutfile{\s\s http://www.maths.warwick.ac.uk/gt/GTVol\thevolumenumber/paper\thepapernumber.abs.html}
\immediate\write\gtoutfile{\noexpand\\}
\immediate\write\gtoutfile{}
\ifx\theasciiabstract\relax
\immediate\write\gtoutfile{\theabstract}\else
\immediate\write\gtoutfile{\theasciiabstract}\fi
\immediate\write\gtoutfile{}
\immediate\write\gtoutfile{\noexpand\\}
\immediate\write\gtoutfile{}
\immediate\closeout\gtoutfile}}  

\def\maketitlepage{\maketitlep\makeheadfile}
\let\maketitle\maketitlepage


\def\ifplaintex{\expandafter\ifx\csname documentclass\endcsname\relax}


\ifplaintex 
\hoffset 14truemm
\voffset 31truemm
\else
\headsep 23pt
\footskip 35pt
\hoffset -4truemm
\voffset 12.5truemm
\fi

\expandafter\ifx\csname beginpicture\endcsname\relax
\expandafter\ifx\csname documentclass\endcsname\relax
\input pictex \else\font\fiverm=cmr5
\input prepictex \input pictex \input postpictex \fi\fi

\def\gt{{\mathsurround=0pt\it $\cal G\mskip-2mu$eometry \&\ 
$\cal T\!\!$opology}}        

\def\gtp{{\mathsurround=0pt\it $\cal G\mskip-2mu$eometry \&\ 
$\cal T\!\!$opology $\cal P\!$ublications}}  


\def\lognumber#1{\def\thelognumber{#1}}
\def\volumenumber#1{\def\thevolumenumber{#1}}
\def\papernumber#1{\def\thepapernumber{#1}}
\def\volumeyear#1{\def\thevolumeyear{#1}}

\def\pagenumbers#1#2{\def\startpage{#1}\def\finishpage{#2}}
\def\published#1{\def\publishdate{#1}}
\def\proposed#1{\def\theproposer{#1}}
\def\seconded#1{\def\theseconders{#1}}
\def\received#1{\def\receiveddate{#1}}
\def\revised#1{\def\reviseddate{#1}}
\def\accepted#1{\def\accepteddate{#1}}

\long\def\asciiabstract#1{\long\def\theasciiabstract{#1}}


\let\\\par\let\thelognumber\relax
\let\thevolumenumber\relax\let\thepapernumber\relax
\let\thevolumeyear\relax\let\thesamplenumber\relax\let\startpage\relax
\let\finishpage\relax\let\publishdate\relax\let\receiveddate\relax
\let\reviseddate\relax\let\accepteddate\relax\let\theasciititle\relax
\let\theasciiauthors\relax
\let\theasciiabstract\relax
\let\theasciiemail\relax\let\theshortauthors\relax\let\theshorttitle\relax

\long\def\maketitlep{   

\count0=\startpage

\gt\hfill      
\beginpicture
\setcoordinatesystem units <0.33truein, 0.33truein> point at 2.2 0.9
\setplotsymbol ({$\cal G$})
\plotsymbolspacing=9truept
\circulararc 315 degrees from 0 1 center at 0 0
\setplotsymbol ({$\cal T$})
\circulararc 315 degrees from 1 -1 center at 1 0
\endpicture
%
\break
{\small\ifx\thesamplenumber\relax 
Volume \else Sample
\fi\thevolumenumber\ (\thevolumeyear)
\startpage--\finishpage\nl
Published: \publishdate}
\vglue 0.5truein plus 0.4fil minus 0.1truein

{\parskip=0pt\leftskip 0pt plus 1fil\def\\{\par\smallskip}{\ifplaintex\large
\else\Large\fi\bf\thetitle}\par\medskip}   

\vglue 0pt plus 0.1fil 

{\parskip=0pt\leftskip 0pt plus 1fil\def\\{\par}{\sc\theauthors}
\par\medskip}

\vglue 0pt plus 0.1fil 

{\small\parskip=0pt\let\newline\\
{\leftskip 0pt plus 1fil\def\\{\par}{\sl\theaddress}\par}
\expandafter\ifx\theemail\relax    
\relax\else\vglue 5pt plus 0.02fil minus 2pt\def\\{\stdspace{\rm 
and}\stdspace} 
\cl{Email:\stdspace\tt\theemail}\fi
\ifx\theurl\relax                  
\relax\else\vglue 5pt plus 0.02fil minus 2pt\def\\{\stdspace{\rm 
and}\stdspace}
\cl{URL:\stdspace\tt\theurl}\fi\par}

\vglue 7pt plus 0.3fil minus 3pt

{\bf Abstract}
\vglue 5pt plus 0.1fil minus 2pt

\theabstract

\vglue 7pt plus 0.3fil minus 3pt

{\bf AMS Classification numbers}\quad Primary:\quad \theprimaryclass

Secondary:\quad \thesecondaryclass

\vglue 5pt plus 0.3fil minus 2pt

{\bf Keywords}\quad \thekeywords

\vglue 10pt plus 0.5fil minus 5pt

{\small  Proposed: \theproposer\hfill Received: \receiveddate\nl
Seconded: \theseconders\hfill 
\ifx\reviseddate\relax                         
Accepted: \accepteddate                        
\else
Revised: \reviseddate                          
\fi}
\eject
}       

\let\maketitlepage\maketitlep
\let\maketitle\maketitlepage


\font\phead=cmsl9 scaled 950
\font\lhead=cmsl9 scaled 1050
\font\pnum=cmbx10 scaled 913
\font\lnum=cmbx10 
\font\pfoot=cmsl9 scaled 950
\font\lfoot=cmsl9 scaled 1050
\ifplaintex
\headline{\vbox to 0pt{\vskip -4.5mm\line{\small\phead\ifnum
\count0=\startpage ISSN 1364-0380 (on line)
1465-3060 (printed) \hfill {\pnum\folio}\else\ifodd\count0\def\\{ }%
\ifx\theshorttitle\relax\thetitle\else\theshorttitle\fi\hfill{\pnum\folio}
\else\def\\{ and }{\pnum\folio}\hfill\ifx\theshortauthors\relax\theauthors
\else\theshortauthors\fi\fi\fi}\vss}}
\footline{\vbox to 0pt{\vglue 0mm\line{\small\pfoot\ifnum\count0=\startpage
\copyright\ \gtp\hfill\else
\gt, Volume \thevolumenumber\ (\thevolumeyear)\hfill\fi}\vss
}}
\else
\makeatletter
\def\@oddhead{{\small\lhead\ifnum\count0=\startpage ISSN 1364-0380 (on line)
1465-3060 (printed) \hfill {\lnum\number\count0}\else\ifodd\count0
\def\\{ }\ifx\theshorttitle\relax \thetitle \else\theshorttitle\fi\hfill
{\lnum\number\count0}\else\def\\{ and }{\lnum\number\count0}
\hfill\ifx\theshortauthors\relax 
\theauthors\else\theshortauthors\fi\fi\fi}}\def\@evenhead{\@oddhead}
\def\@oddfoot{\small\lfoot\ifnum\count0=\startpage\copyright\ \gtp\hfill\else
\gt, Volume \thevolumenumber\ (\thevolumeyear)\hfill\fi}
\def\@evenfoot{\@oddfoot}
\makeatother
\fi


\newwrite\gtoutfile
\long\gdef\makeheadfile{  
{\def\\{, }\def\s{ }
\immediate\openout\gtoutfile head.xxx
\immediate\write\gtoutfile{To: math@arxiv.org}
\immediate\write\gtoutfile{Subject: put or rep NNNNN:pppp}
\immediate\write\gtoutfile{--text follows this line--}
\immediate\write\gtoutfile{Proxy-for: \ifx\theasciiauthors\relax
\theauthors\else\theasciiauthors\fi\s<\ifx\theasciiemail\relax\theemail\else\theasciiemail\fi>}
\immediate\write\gtoutfile{\noexpand\\}
\immediate\write\gtoutfile{Authors: \ifx\theasciiauthors\relax
\theauthors\else\theasciiauthors\fi}
\immediate\write\gtoutfile{Title: \ifx\theasciititle\relax
\thetitle\else\theasciititle\fi}
\immediate\write\gtoutfile{Subj-class: GT or SG or MG etc}
\immediate\write\gtoutfile{MSC-class: \theprimaryclass\ifx\thesecondaryclass\relax\else, \thesecondaryclass\fi}
\immediate\write\gtoutfile{Journal-ref: Geom. Topol. \thevolumenumber
(\thevolumeyear) \startpage-\finishpage}
\immediate\write\gtoutfile{Comments: Published by Geometry and Topology at}
\immediate\write\gtoutfile{\s\s http://www.maths.warwick.ac.uk/gt/GTVol\thevolumenumber/paper\thepapernumber.abs.html}
\immediate\write\gtoutfile{\noexpand\\}
\immediate\write\gtoutfile{}
\ifx\theasciiabstract\relax
\immediate\write\gtoutfile{\theabstract}\else
\immediate\write\gtoutfile{\theasciiabstract}\fi
\immediate\write\gtoutfile{}
\immediate\write\gtoutfile{\noexpand\\}
\immediate\write\gtoutfile{}
\immediate\closeout\gtoutfile}}  

\def\maketitlepage{\maketitlep\makeheadfile}
\let\maketitle\maketitlepage

\lognumber{178}
\volumenumber{5}\papernumber{14}\volumeyear{2001}
\pagenumbers{399}{429}
\accepted{23 April 2001}
\received{25 January 2001}
\revised{2 April 2001}
\published{24 April 2001}
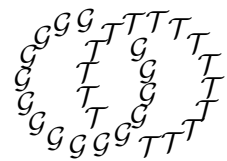
\proposed{Robion Kirby}
\seconded{Yasha Eliashberg, David Gabai}
%
%
%
\def\tw{{\tt \char'176}}  

\def\G{{\cal G}}
\def\ep{\varepsilon}

\font\spec=cmtex10 scaled 1095 
\def\d{\hbox{\spec \char'017\kern 0.05em}} 
\def\inv{^{-1}}
\def\tr{\hbox{tr}}
\def\verts{\thinspace\vert\thinspace}
%
%
\long\def\mbox#1{$$\vbox{#1}$$}
\newdimen\unitlength
\unitlength= 1.000cm
\def\makebox(0,0)[l]#1{\hbox to 0pt{#1\hss}}
\def\symbol#1{\char#1}
\font\tencirc=lcircle10
\def\SetFigFont#1#2#3#4{\small$#4$}
\font\thinlinefont=cmr7
%
%
%
\reflist

\refkey\Boa
{\bf J\,M Boardman}, {\it Singularities of differentiable maps},
Publ. IHES, {33} (1967) 21--57

\refkey\BH 
{\bf W Browder},  {\bf W\,C Hsiang}, {\it Some problems on homotopy theory,
manifolds and transformation groups}, (Proceedings of the Stanford
conference in algebraic and geometric topology, 1976),
AMS Proc. Sympos. Pure Math. XXXII, {2} (1978) 251--267

\refkey\EliMis
{\bf Y\,M Eliashberg}, {\bf N\,M Mishachev}, {\it Wrinkling of smooth
maps II: wrinkling of embeddings and K Igusa's theorem}, Topology, 39
(2000) 711--732

\refkey\FRS
{\bf R Fenn}, {\bf C Rourke}, {\bf B Sanderson}, {\it James bundles
and applications},\nl {\tt http://www.maths.warwick.ac.uk/\tw
cpr/ftp/james.ps}

\refkey\Gr
{\bf M Gromov}, {\it Partial differential relations}, Springer--Verlag (1986)

\refkey\GP
{\bf V Guillemin}, {\bf A Pollack}, {\it Differential Topology}, Prentice--Hall
(1974)

\refkey\Hi
{\bf M Hirsch}, {\it Immersions of manifolds}, Trans. Amer. Math. Soc.
{93} (1959) 242--276

\refkey\Ja 
{\bf I James}, {\it Reduced product spaces}, Annals of Math. {62} (1955) 
170--197     

\refkey\KSa 
{\bf U Koschorke}, {\bf B Sanderson}, {\it Geometric interpretations
of the generalised Hopf invariant}, Math. Scand. {41} (1977) 199--217

\refkey\Ma 
{\bf J\,P May}, {\it The geometry of iterated loop spaces},
Springer Lectures Notes Series, {271} Springer-Verlag (1972)

\refkey\Mat
{\bf J\,N Mather}, {\it Generic projections}, Annals of Math. {98}
(1973) 226--245

\refkey\Mi 
{\bf R\,J Milgram}, {\it Iterated loop spaces}, Annals of Math. 
{84} (1966) 386--403

\refkey\Ph
{\bf A Philips}, {\it Turning a surface inside out}, Scientific American,
{214} (1966) 112--120

\refkey\CompII
{\bf Colin Rourke}, {\bf Brian Sanderson}, {\it The compression theorem II:
directed embeddings}, Geometry and Topology, 5 (2001) 431--440,
{\tt arxiv:math.GT/0003026}

\refkey\CompIII
{\bf Colin Rourke}, {\bf Brian Sanderson}, {\it The compression theorem III:
applications}, to appear, see {\tt http://www.maths.warwick.ac.uk/\tw
cpr/ftp/compIII.ps} for a preliminary version

\refkey\Se 
{\bf G Segal}, {\it Configuration-spaces and iterated loop-spaces},
Invent. Math. {21} (1973) 213--221

\refkey\Sm 
{\bf S Smale}, {\it The classification of immersions of spheres in
Euclidean spaces}, Annals of Math. {69} (1959) 327--344

\refkey\Wi
{\bf B Wiest}, {\it Loop spaces and the compression theorem},
Proc. Amer. Math. Soc. 128 (2000) 3741--3747
   
\refkey\Wie
{\bf B Wiest}, {\it Rack spaces and loop spaces}, J. Knot Theory
Ramifications, 8 (1999) 99--114

\endreflist
%
%
\title{The compression theorem I}               

\author{Colin Rourke\\Brian Sanderson}  

\address{Mathematics Institute, University of Warwick\\
Coventry, CV4 7AL, UK}
\email{cpr@maths.warwick.ac.uk\\bjs@maths.warwick.ac.uk}
\url{http://www.maths.warwick.ac.uk/\tw cpr/\qua {\rm and}\qua \tw bjs/}

\abstract 

This the first of a set of three papers about the {\it Compression
Theorem\/}: if $M^m$ is embedded in $Q^q\times\re$ with a normal
vector field and if $q-m\ge1$, then the given vector field can be {\sl
straightened} (ie, made parallel to the given $\re$ direction) by an
isotopy of $M$ and normal field in $Q\times\re$.

The theorem can be deduced from Gromov's theorem on directed
embeddings [\Gr; 2.4.5 $(\rm C')$] and is implicit in the preceeding
discussion.  Here we give a direct proof that leads to an explicit
description of the finishing embedding.

In the second paper in the series we give a proof in the spirit of
Gromov's proof and in the third part we give applications.

\endabstract

\asciiabstract{
This the first of a set of three papers about the Compression Theorem:
if M^m is embedded in Q^q X R with a normal vector field and if q-m >
0, then the given vector field can be straightened (ie, made parallel
to the given R direction) by an isotopy of M and normal field in Q X R
.  The theorem can be deduced from Gromov's theorem on directed
embeddings [M Gromov, Partial differential relations, Springer--Verlag
(1986); 2.4.5 C'] and is implicit in the preceeding discussion.  Here
we give a direct proof that leads to an explicit description of the
finishing embedding.  In the second paper in the series we give a
proof in the spirit of Gromov's proof and in the third part we give
applications.}

\primaryclass{57R25}

\secondaryclass{57R27, 57R40, 57R42, 57R52}

\keywords{Compression, embedding, isotopy, immersion, straightening,
vector field}

\maketitlepage

\section{Introduction}  
  
We work throughout in the smooth ($C^\infty$) category.  Embeddings,
immersions, regular homotopies etc will be assumed either to take
boundary to boundary or to meet the boundary in a codimension 0
submanifold of the boundary.  Thus for example if $f\co M\to Q$ is an
immersion then we assume that either $f\inv{\d Q}=\d M$ or $f\inv{\d
Q}$ is a codimension 0 submanifold of $\d M$.  In the latter case we
speak of $M$ having {\sl relative boundary}.  The tangent bundle of a
manifold $W$ is denoted $T(W)$ and the tangent space at $x\in W$ is
denoted $T_x(W)$. Throughout the paper, ``normal'' means independent
(as in the usual meaning of ``normal bundle'') and not necessarily
perpendicular.

This paper is about the following result:

\proclaim{Compression Theorem}Suppose that $M^m$ is embedded in 
$Q^q\times\re$ with a normal vector field and suppose that
$q-m\ge1$. Then the vector field can be {\sl straightened} (ie,
made parallel to the given $\re$ direction) by an isotopy of $M$ and
normal field in $Q\times\re$.\endproc

Thus the theorem moves $M$ to a position where it projects
by vertical projection (ie ``compresses'') to an immersion 
in $Q$.

The theorem can be deduced from Gromov's theorem on directed
embeddings [\Gr; 2.4.5 $(\rm C')$] and is implicit in the discussion
which precedes Gromov's theorem.  Here we present a proof that is
different in character from Gromov's proof.

Gromov uses the technique of ``convex integration'' which although
simple in essence leads to very complicated embeddings.  In the second
paper in this series [\CompII] we give a proof of Gromov's theorem
(and deduce the compression theorem).  This proof is in the spirit of
Gromov's and makes clear the complicated (and uncontrolled) nature of
the resulting embeddings.

By contrast the proof that we present in this paper is completely
constructive: given a particular embedding and vector field the
resulting compressed embedding can be described explicitly.  A third,
and again different, proof of the Compression Theorem has been given
by Eliashberg and Mishachev [\EliMis].

The method of proof allows for a number of natural addenda to be
proved and in particular we can straighten a sequence of vector
fields.  More precisely suppose that $M$ is embedded in $Q\times
\re^n$ with $n$ independent normal vector fields, then $M$ is isotopic  
to an embedding in which each vector field is parallel to the
corresponding copy of $\re$.  This result solves an old problem (posed
by Bruce Williams at the 1976 Stanford Conference [\BH; problem 6])
though it should be noted that this solution could also have been
deduced from Gromov's theorem at any time since the publication of
his book.

The third paper in this series [\CompIII] concerns applications of the
Compression Theorem (and its addenda): we give new and constructive
proofs for immersion theory [\Hi, \Sm] and for the loops--suspension
theorem of James, May, Milgram and Segal [\Ja, \Ma, \Mi, \Se].  We
give a new approach to classifying embeddings of manifolds in
codimension one or more, which leads to theoretical solutions, and we
consider the general problem of simplifying (or specifying) the
singularities of a smooth projection up to $C^0$--small isotopy and
give a theoretical solution in the codimension $\ge1$ case.

Two examples of immediate application of the compression theorem are
the following:

\proc{Corollary}Let $\pi$ be a group.  There is a classifying
space $BC(\pi)$ such that the set of homotopy classes $[Q,BC(\pi)]$
is in natural bijection with the set of cobordism classes of 
framed submanifolds $L$ of $Q\times\re$ of codimension 2
equipped with a homomorphism $\pi_1(Q\times\re-L)\to\pi$.

\proc{Corollary}Let $Q$ be a connected manifold with basepoint $*$
and let $M$ be any collection of disjoint submanifolds of $Q-\{*\}$
each of which has codimension $\ge2$ and is equipped with a normal
vector field.  Define the {\rm vertical} loop space of $Q$ denoted
$\Omega^{\rm vert}(Q)$ to comprise loops which meet given tubular
neighbourhoods of manifolds in $M$ in straight line segments parallel
to the given vector field.  Then the natural inclusion $\Omega^{\rm
vert}(Q)\subset\Omega(Q)$ (where $\Omega(Q)$ is the usual loop space)
is a weak homotopy equivalence.\rm

\rk{Proofs} The first corollary is a special case of the 
classification theorem for links in codimension 2 given in
[\FRS; theorem 4.15], the space $BC(\pi)$ being the rack space of the
conjugacy rack of $\pi$.  To prove the second corollary suppose 
given a based map $f\co S^n\to\Omega(Q)$ then the adjoint
of $f$ can be regarded as a map $g\co S^n\times\re\to Q$ which
takes the ends of $S^n\times\re$ and $\{*\}\times\re$  to the 
basepoint.  Make $g$
transverse to $M$ to create a number of manifolds embedded in
$S^n\times\re$ and equipped with normal vector fields.  Apply
the compression theorem to each of these (the local version
proved in section 4 of this paper).  The result is to deform
$g$ into the adjoint of a map  $S^n\to\Omega^{\rm vert}(Q)$.
This shows that $\Omega^{\rm vert}(Q)\subset\Omega(Q)$ induces
a surjection on $\pi_n$.  A similar argument applied to a 
homotopy, using the relative compression theorem, proves injectivity.
(The vertical loop space is introduced in Wiest [\Wi]; for applications
and related results see [\Wi, \Wie].) \qed

The result for one vector field (sufficient for the above
applications) has a particularly simple global proof given in the next
section (section 2).  In section 3 we describe some of the geometry
which results from this proof and in section 4 we localise the proof
and show that the isotopy can be assumed to be arbitrarily small in
the $C^0$ sense.  This small version is needed for straightening
multiple vector fields and leads at once to deformation and bundle
versions, which lead to the connection with Gromov's results.

Finally section A is an appendix which contains proofs of the general
position and transversality results that are needed in sections 2 and
4.

\rk{Acknowledgements}We are grateful to Bert Wiest for observing
corollary 1.2 and to David Mond for help with smooth general position.
We are also extremely grateful to Chris French who has carefully read
the main proof and found some technical errors which have now been
corrected.  We are also grateful to David Spring and Yasha Eliashberg
for comments on the connection of our results with those of Gromov.  

\section{The global proof}    

We think of $\re$ as vertical and the positive $\re$ direction    
as upwards.  We call $M$ {\sl compressible} if the    
vector field always points vertically up.  Note that a    
compressible embedding covers an immersion in $Q$.    
    
Assume that $Q$ is equipped     
with a Riemannian metric and use the product metric on $Q\times\re$.    
Call a normal field {\sl perpendicular} if it is everywhere 
orthogonal to $M$.

A perpendicular vector field $\alpha$ is said to be {\sl grounded} if    
it never points vertically down.  More generally $\alpha$ is said
to be {\sl $\ep$--grounded} if it always makes an angle of at
least $\ep$ with the downward vertical, where $\ep>0$.

\numkey\GComp
\proc{Compression Theorem}Let $M^m$ be a compact manifold embedded    
in $Q^q\times\re$ and equipped with a normal vector    
field.  Assume $q-m\ge 1$ then $M$ is isotopic    
to a compressible embedding.\rm

The method of proof allows a number of extensions, but note that
the addenda to the local proof (in section 4) give stronger
statements for the first two of these addenda:

\rk{Addenda}\sl    
\items    
\item{\rm(i)}{\rm (Relative version)}\qua    
Let $C$ be a compact set in $Q$. If $M$ is already compressible in a    
neighbourhood of $C\times \re$ then the isotopy can be    
assumed fixed on $C\times \re$.    
    
\item{\rm(ii)}{\rm (Parametrised version)}\qua 
Given a parametrised family $M^m_t\subset Q^q_t$    
of embeddings with a normal vector field, $t\in K$,    
where $K$ is a compact manifold of dimension $k$ and    
$q-m-k\ge 1$, then there is a parametrised family of    
isotopies to compressible embeddings (and there is a relative version    
similar to (i)).    
    
\item{\rm(iii)}If in theorem \GComp\ or in addendum
(ii) the fields are all perpendicular and grounded, then the    
dimension condition can be relaxed to $q-m\ge 0$ and there is no dimension    
condition on $K$.\rm    
\enditems
    
We need the following lemma.

\numkey\Ground
\proc{Lemma} Under the hypotheses of theorem \GComp\ the normal field 
may be assumed to be perpendicular and grounded.\rm

The lemma follows from general position.  Call the vector field $\alpha$ 
and without loss assume that $\alpha$ has unit length everywhere.      
Note that the fact that $\alpha$ is normal (ie,     
independent of the tangent plane at each point of $M$) does not imply     
that it is perpendicular; however we can isotope    
$\alpha$ without further moving $M$ to make it perpendicular.  
$\alpha$ now defines a section of $T(Q\times\re)\verts M$ and
vertically down defines another section.  The condition that
$q-m>0$ implies that these two sections are not expected to
meet in general position.  A formal proof can be found in the appendix,
see corollary A.5.

\rk{Proof of theorem \GComp}   
We shall prove the main result first and then the addenda.  The first
move is to apply the lemma which results in the normal vector field
$\alpha$ being perpendicular and grounded.  By compactness of $M$,
$\alpha$ is in fact $\ep$--grounded for some $\ep>0$.

We now define an operation on $\alpha$ given by rotating it towards 
the upward vertical:   Choose a real number $\mu$
with $0<\mu<\ep$.  Consider a point $p\in M$ at which $\alpha(p)$ does not    
point vertically up and consider the plane $P(p)$ in $T_p(Q\times\re)$
defined by the    
vector $\alpha(p)$ and the vertical.  Define the vector $\beta(p)$ to    
be the vector in the plane $P(p)$ obtained by rotating $\alpha(p)$    
through an angle ${\pi\over2}-{\mu}$ in the direction    
towards vertically up, unless this rotation carries $\alpha(p)$ past    
vertically up, when we define $\beta(p)$ to {\it be} vertically up.  If    
$\alpha(p)$ is already vertically up, then we again define $\beta(p)$    
to be vertically up.  The rotation of $\alpha$ to $\beta$ is called 
{\sl upwards rotation}.  If we wish to be precise and refer to the
chosen real number $\mu$ then we say {\sl $\mu$--upwards rotation}.

Figure \figkey\Rotate\ shows the extreme case when $\alpha$ is 
pointing as far down as possible:    

\fig{\Rotate: upwards rotation}
\beginpicture    
\setcoordinatesystem units < 1.000cm, 1.000cm>    
\unitlength= 1.000cm    
\linethickness=1pt    
\setplotsymbol ({\makebox(0,0)[l]{\tencirc\symbol{'160}}})    
\setshadesymbol ({\thinlinefont .})    
\setlinear    
%
%
\linethickness= 0.500pt    
\setplotsymbol ({\thinlinefont .})    
\circulararc 108.157 degrees from  2.477 24.575 center at  2.811 24.381    
%
%
\linethickness= 0.500pt    
\setplotsymbol ({\thinlinefont .})    
\circulararc 67.380 degrees from  2.762 23.844 center at  2.921 23.979    
%
%
\linethickness= 0.500pt    
\setplotsymbol ({\thinlinefont .})    
\putrule from  2.985 26.289 to  2.985 22.606    
%
%
\linethickness= 0.500pt    
\setplotsymbol ({\thinlinefont .})    
\plot  0.953 25.337  5.207 23.590 /    
%
%
\linethickness=1pt    
\setplotsymbol ({\makebox(0,0)[l]{\tencirc\symbol{'160}}})    
\plot  2.985 24.479  2.223 22.701 /    
%
%
\plot  2.264 22.960  2.223 22.701  2.381 22.910 /    
%
%
%
\linethickness=1pt    
\setplotsymbol ({\makebox(0,0)[l]{\tencirc\symbol{'160}}})    
\plot  2.985 24.479  0.953 25.019 /    
%
%
\plot  1.214 25.015  0.953 25.019  1.182 24.892 /    
%
%
%
\linethickness= 0.500pt    
\setplotsymbol ({\thinlinefont .})    
\plot  2.445 24.575  2.350 24.448 /    
%
%
\linethickness= 0.500pt    
\setplotsymbol ({\thinlinefont .})    
\plot  2.477 24.575  2.508 24.416 /    
%
%
\put{\SetFigFont{12}{14.4}{rm}M} [lB] at  4.635 24.003    
%
%
\put{\SetFigFont{12}{14.4}{rm}{\beta}} [lB] at  0.730 24.638    
%
%
\put{\SetFigFont{12}{14.4}{rm}{\alpha}} <-2pt, 0pt> [lB] at  1.937 22.796    
%
%
\put{\SetFigFont{12}{14.4}{rm}{{\pi\over2}-\mu}} <3pt,    
0pt> [lB] at  1.365 23.971    
%
%
\put{\SetFigFont{12}{14.4}{rm}{\varepsilon}} [lB] at  2.730 23.463    
\linethickness=0pt    
\putrectangle corners at  0.730 26.289 and  5.207 22.606    
\endpicture
\endfig    

The operation of upwards rotation, just defined, yields a {\it continuous}
but not in general {\it smooth} vector field.  However the operation may
be altered to yield a smooth vector field by using a bump function 
to phase out the amount of rotation as the rotated vector approaches
vertical.  The properties of the resulting vector field for the proof 
(below) are not altered by this smoothing.

The second move is to perform $\mu$--upwards rotation of $\alpha$.
The resulting vector field $\beta$ has the property that 
it is still normal (though not now perpendicular) to $M$ and has a 
positive vertical component of at least $\sin(\ep-\mu)$.  These facts 
can both be seen in figure \Rotate; $\beta$ makes an angle of at least    
$\ep-\mu$ above the horizontal and at least $\mu$    
with $M$.  The line marked $M$ in the figure is where the tangent plane    
to $M$ at $p$ might meet $P(p)$.

Now extend $\beta$ to a global unit vector field $\gamma$ 
on $Q\times\re$ by taking  $\beta$ on $M$ and the vertically up 
field outside a tubular neighbourhood $N$ of $M$ in $Q\times\re$ 
and interpolating by rotating $\beta$ to vertical    
along radial lines in the tubular neighbourhood. 

Call the flow defined by $\gamma$ the {\sl global flow} on
$Q\times\re$.  Notice that this global flow is determined by
$\alpha$ and the two choices of $\mu<\ep$ and of the tubular
neighbourhood $N$ of $M$.

Since the vertical component of $\gamma$ is positive and bounded 
away from zero (by $\sin(\ep-\mu)$) any point will flow upwards 
in the global flow as far as we like in finite time.    
    
Now let $M$ flow in the global flow.  In finite time    
we reach a region where $\gamma$ is vertically up.  Since 
$\gamma\vert M=\beta$ is normal to $M$ at the start of the
flow, $\gamma\vert M$ remains normal to $M$ throughout the flow
and we have isotoped $M$ together with its normal vector field 
to a compressible embedding. \qed

\rk{Proofs of the addenda}

To prove addendum (i) we modify the global flow to be stationary 
on $C\times \re$ as follows.  

Suppose that $x\in Q\times\re$ and $s\in\re$ then let $x+s$ 
denote the point obtained by moving $x$ vertically by $s$.  Let
the global flow defined by $\gamma$ be given by $x\mapsto f_t(x)$
then the {\sl modified global flow} is given by $x\mapsto f_t(x)-t$.
In words, this flow is obtained by flowing along $\gamma$ for time
$t$ and then flowing back down the unit downward flow again for
time $t$.

We apply this flow to $M$ and the vector field $\gamma$.  To avoid 
confusion use the notation $\gamma_t$, $M_t$ for the effect of this 
modified global flow on $\gamma$, $M$ at time $t$ respectively.  
Note that $\gamma_t(x)=\gamma(x+t)$, where we have identified the
tangent spaces at $x$ and $x+t$ in the canonical way, and that
$\gamma_t\vert M_t$ is normal to $M_t$ for all $t$.

Now the modified flow is stationary 
whenever and wherever $\gamma_t$ is vertically upwards and in 
particular it is stationary on $C\times\re$ for all $t$.
Moreover it has similar properties to the unmodified global flow in
that, after finite time, any compact set reaches a region where the 
vector field $\gamma_t$ is vertically up.

For (ii) the dimension condition implies that the vector fields
can all be assumed to be both perpendicular and grounded.  This
is a general position argument and the details are in the appendix,
see corollary A.5. 
By compactness there is a global $\varepsilon>0$ such that each    
vector field is $\ep$--grounded.  Now apply the main proof for 
each $t\in K$ and observe that the resulting flows vary smoothly 
with $t\in K$.  

For (iii) observe that the dimension conditions were only used to    
prove that the vector fields were all grounded and only compactness    
was used thereafter. \qed    
    
\numkey\TDvf\proc{Remark}\rm

The modified global flow defined in the proof of addendum (i) can
be regarded as given by a time-dependent vector field as follows.
Let $u$ denote the unit vertical vector field and $\gamma$ the
vector field which defines the global flow.  Then the modified
global flow is given by the vector field $\gamma^*$ where
$$\gamma^*(x,t)=\gamma(x+t)-u.$$
Note also that the vector field carried along by $\gamma^*$ (denoted
$\gamma_t$ above) is $\gamma^*+u$.  Here we have again identified 
the tangent spaces at $x$ and $x+t$ in the canonical way.

\section{Pictures}    
    
The global proof given in the last section hides a wealth of
geometry.  In this section we reveal some of this geometry.
We start by drawing a sequence of pictures, where $M$ has codimension
2 in $Q\times\re$, which
are the end result of the isotopy given by the proof in
particular situations.  These pictures contain all the critical
information for constructing an isotopy of a general manifold
of codimension 2 with normal vector field to a compressible 
embedding.  We shall explain how this works in the local setting 
of the next section.  

After the sequence of codimension 2
pictures, we describe some higher codimension situations and 
then we describe how the compression desingularises a map
in a particular case (the removal of a Whitney umbrella).
We finish
with an explicit compression of an (immersed) projective plane 
in $\re^4$ which uses many of the earlier pictures.  The
image of the projection on $\re^3$ changes from a sphere with cross-cap to 
Boy's surface.

\sh{1 in 3}

Consider the vector field on an angled line in $\re^3$ which 
rotates once around the line as illustrated in figure \figkey\vectp.

\fig{\vectp: perpendicular field}
\beginpicture
\setcoordinatesystem units <.70000cm,.70000cm>
\unitlength=1.00000cm
\linethickness=1pt
\setplotsymbol ({\makebox(0,0)[l]{\tencirc\symbol{'160}}})
\setshadesymbol ({\thinlinefont .})
\setlinear
%
%
\linethickness=1pt
\setplotsymbol ({\makebox(0,0)[l]{\tencirc\symbol{'160}}})
\plot 27.616  5.256 31.500  8.134 /
%
%
\linethickness=1pt
\setplotsymbol ({\makebox(0,0)[l]{\tencirc\symbol{'160}}})
\plot 31.966  8.526 33.786  9.976 /
%
%
\linethickness= 0.500pt
\setplotsymbol ({\thinlinefont .})
\plot 33.553  9.775 33.109 10.336 /
%
%
\plot 33.316 10.176 33.109 10.336 33.217 10.097 /
%
%
%
\linethickness= 0.500pt
\setplotsymbol ({\thinlinefont .})
\plot 32.876  9.203 32.432  9.775 /
%
%
\plot 32.638  9.613 32.432  9.775 32.537  9.535 /
%
%
%
\linethickness= 0.500pt
\setplotsymbol ({\thinlinefont .})
\plot 32.167  8.674 31.839  9.161 /
%
%
\plot 32.034  8.986 31.839  9.161 31.928  8.915 /
%
%
%
\linethickness= 0.500pt
\setplotsymbol ({\thinlinefont .})
\plot 31.426  8.018 31.722  7.552 /
%
%
\plot 31.533  7.732 31.722  7.552 31.640  7.801 /
%
%
%
\linethickness= 0.500pt
\setplotsymbol ({\thinlinefont .})
\plot 30.738  7.489 31.162  6.886 /
\put {$\bullet$} <-.5pt, .5pt> at 30.738 7.489
\put {\small $D$} [br] <-1.5pt, 3pt> at 30.738 7.489
%
%
\plot 30.964  7.057 31.162  6.886 31.068  7.130 /
%
%
%
\linethickness= 0.500pt
\setplotsymbol ({\thinlinefont .})
\plot 30.071  7.034 30.336  6.610 /
%
%
\plot 30.148  6.792 30.336  6.610 30.255  6.859 /
%
%
%
\linethickness= 0.500pt
\setplotsymbol ({\thinlinefont .})
\plot 29.119  6.409 29.119  6.409 /
%
%
\linethickness= 0.500pt
\setplotsymbol ({\thinlinefont .})
\plot 29.108  6.325 28.791  6.854 /
%
%
\plot 28.976  6.669 28.791  6.854 28.867  6.603 /
%
%
%
\linethickness= 0.500pt
\setplotsymbol ({\thinlinefont .})
\plot 28.420  5.891 27.997  6.430 /
%
%
\plot 28.204  6.270 27.997  6.430 28.104  6.191 /
%
%
%
\linethickness= 0.500pt
\setplotsymbol ({\thinlinefont .})
\plot 27.754  5.319 27.309  5.891 /
%
%
\plot 27.515  5.729 27.309  5.891 27.415  5.651 /
\linethickness= 0.500pt
\setplotsymbol ({\thinlinefont .})
%
%
%
\plot	33.352 10.569 32.945 10.240
 	32.842 10.158
	32.738 10.074
	32.633  9.990
	32.525  9.904
	32.417  9.818
	32.307  9.730
	32.195  9.642
	32.082  9.553
	31.978  9.457
	31.892  9.349
	31.825  9.230
	31.799  9.166
	31.777  9.099
	31.760  9.029
	31.747  8.956
	31.739  8.880
	31.736  8.801
	31.738  8.720
	31.744  8.635
	31.755  8.548
	31.770  8.457
	31.786  8.366
	31.798  8.278
	31.807  8.191
	31.812  8.107
	31.813  8.025
	31.810  7.946
	31.804  7.868
	31.794  7.793
	31.780  7.720
	31.763  7.649
	31.741  7.581
	31.717  7.515
	31.656  7.389
	31.580  7.272
	31.493  7.164
	31.400  7.066
	31.301  6.977
	31.195  6.897
	31.082  6.827
	30.963  6.767
	30.838  6.716
	30.773  6.694
	30.706  6.674
	30.640  6.656
	30.577  6.640
	30.460  6.611
	30.355  6.587
	30.261  6.569
	30.179  6.557
	30.108  6.550
	30.003  6.552
	 /
\plot 30.003  6.552 29.839  6.579 /
\linethickness= 0.500pt
\setplotsymbol ({\thinlinefont .})
%
%
%
\plot	29.384  6.759 29.241  6.827
 	29.160  6.857
	29.061  6.878
	28.943  6.889
	28.877  6.891
	28.807  6.891
	28.728  6.882
	28.639  6.857
	28.538  6.816
	28.426  6.760
	28.302  6.688
	28.237  6.646
	28.168  6.600
	28.096  6.550
	28.022  6.497
	27.945  6.439
	27.865  6.378
	 /
\plot 27.865  6.378 27.214  5.870 /
\linethickness=0pt
\putrectangle corners at 27.197 10.585 and 33.833  5.209
\endpicture
\endfig

Upwards rotation replaces this vector field by one
which is vertically up outside an interval and rotates just
under the line as illustrated in figure \figkey\vectu.

\fig{\vectu: upwards rotated field}
\beginpicture
\setcoordinatesystem units <.70000cm,.70000cm>
\unitlength=1.00000cm
\linethickness=1pt
\setplotsymbol ({\makebox(0,0)[l]{\tencirc\symbol{'160}}})
\setshadesymbol ({\thinlinefont .})
\setlinear
%
%
\linethickness=1pt
\setplotsymbol ({\makebox(0,0)[l]{\tencirc\symbol{'160}}})
\plot 27.616  5.256 31.500  8.134 /
%
%
\linethickness=1pt
\setplotsymbol ({\makebox(0,0)[l]{\tencirc\symbol{'160}}})
\plot 31.966  8.526 33.786  9.976 /
%
%
\linethickness= 0.500pt
\setplotsymbol ({\thinlinefont .})
\plot 33.553  9.775 33.585 10.674 /
%
%
\plot 33.640 10.418 33.585 10.674 33.513 10.423 /
%
%
%
\linethickness= 0.500pt
\setplotsymbol ({\thinlinefont .})
\plot 32.876  9.203 32.961 10.209 /
%
%
\plot 33.003  9.950 32.961 10.209 32.876  9.961 /
%
%
%
\linethickness= 0.500pt
\setplotsymbol ({\thinlinefont .})
\plot 32.167  8.674 32.241  9.658 /
%
%
\plot 32.285  9.400 32.241  9.658 32.159  9.410 /
%
%
%
\linethickness= 0.500pt
\setplotsymbol ({\thinlinefont .})
\plot 31.426  8.018 31.786  8.653 /
%
%
\plot 31.716  8.401 31.786  8.653 31.606  8.463 /
%
%
%
\linethickness= 0.500pt
\setplotsymbol ({\thinlinefont .})
\plot 30.738  7.489 31.807  7.669 /
%
%
\plot 31.567  7.564 31.807  7.669 31.546  7.689 /
%
%
%
\linethickness= 0.500pt
\setplotsymbol ({\thinlinefont .})
\plot 30.071  7.034 31.394  7.108 /
%
%
\plot 31.144  7.030 31.394  7.108 31.137  7.157 /
%
%
%
\linethickness= 0.500pt
\setplotsymbol ({\thinlinefont .})
\plot 29.119  6.409 29.119  6.409 /
%
%
\linethickness= 0.500pt
\setplotsymbol ({\thinlinefont .})
\plot 29.108  6.325 29.299  6.769 /
%
%
\plot 29.257  6.511 29.299  6.769 29.140  6.561 /
%
%
%
\linethickness= 0.500pt
\setplotsymbol ({\thinlinefont .})
\plot 28.420  5.891 28.495  6.748 /
%
%
\plot 28.536  6.489 28.495  6.748 28.409  6.500 /
%
%
%
\linethickness= 0.500pt
\setplotsymbol ({\thinlinefont .})
\plot 27.754  5.319 27.870  6.335 /
%
%
\plot 27.904  6.076 27.870  6.335 27.778  6.090 /
%
%
%
\linethickness= 0.500pt
\setplotsymbol ({\thinlinefont .})
\plot 31.236  7.933 31.881  8.156 /
%
%
\plot 31.662  8.013 31.881  8.156 31.620  8.133 /
%
%
%
\linethickness= 0.500pt
\setplotsymbol ({\thinlinefont .})
\plot 29.669  6.748 30.484  6.716 /
%
%
\plot 30.228  6.663 30.484  6.716 30.233  6.790 /
\linethickness= 0.500pt
\setplotsymbol ({\thinlinefont .})
%
%
%
\plot	33.680 10.812 33.109 10.362
 	33.038 10.306
	32.968 10.251
	32.899 10.197
	32.830 10.143
	32.763 10.090
	32.697 10.037
	32.631  9.986
	32.567  9.935
	32.440  9.835
	32.317  9.738
	32.198  9.644
	32.082  9.553
	31.979  9.458
	31.896  9.353
	31.835  9.239
	31.794  9.116
	31.781  9.051
	31.774  8.983
	31.772  8.913
	31.775  8.840
	31.783  8.765
	31.796  8.688
	31.815  8.608
	31.839  8.526
	31.863  8.443
	31.882  8.362
	31.897  8.281
	31.906  8.202
	31.911  8.124
	31.911  8.047
	31.907  7.971
	31.897  7.896
	31.883  7.823
	31.864  7.750
	31.840  7.679
	31.811  7.609
	31.778  7.540
	31.739  7.473
	31.696  7.406
	31.648  7.341
	31.546  7.218
	31.438  7.108
	31.327  7.011
	31.211  6.928
	31.090  6.858
	30.965  6.801
	30.901  6.777
	30.835  6.757
	30.769  6.740
	30.701  6.727
	30.635  6.715
	30.571  6.705
	30.452  6.685
	30.345  6.669
	30.250  6.657
	30.167  6.647
	30.095  6.640
	29.987  6.637
	 /
\plot 29.987  6.637 29.817  6.642 /
\linethickness= 0.500pt
\setplotsymbol ({\thinlinefont .})
%
%
%
\plot	29.384  6.759 29.241  6.827
 	29.160  6.857
	29.061  6.878
	28.943  6.889
	28.877  6.891
	28.807  6.891
	28.728  6.882
	28.639  6.857
	28.538  6.816
	28.426  6.760
	28.302  6.688
	28.237  6.646
	28.168  6.600
	28.096  6.550
	28.022  6.497
	27.945  6.439
	27.865  6.378
	 /
\plot 27.865  6.378 27.214  5.870 /
\linethickness=0pt
\putrectangle corners at 27.197 10.829 and 33.833  5.209
\endpicture
\endfig

Seen from on top this vector field has the form illustrated in 
figure \figkey\vectt.

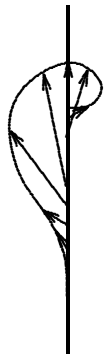
\fig{\vectt: the upwards rotated field seen from on top}
\beginpicture
\setcoordinatesystem units <.70000cm,.70000cm>
\unitlength=1.00000cm
\linethickness=1pt
\setplotsymbol ({\makebox(0,0)[l]{\tencirc\symbol{'160}}})
\setshadesymbol ({\thinlinefont .})
\setlinear
%
%
\linethickness=1pt
\setplotsymbol ({\makebox(0,0)[l]{\tencirc\symbol{'160}}})
\putrule from 30.065 11.235 to 30.065  5.203
%
%
\linethickness= 0.500pt
\setplotsymbol ({\thinlinefont .})
\plot 30.086  9.309 30.414  9.489 /
%
%
\plot 30.222  9.311 30.414  9.489 30.161  9.423 /
%
%
%
\linethickness= 0.500pt
\setplotsymbol ({\thinlinefont .})
\plot 30.044  8.970 30.404 10.124 /
%
%
\plot 30.389  9.863 30.404 10.124 30.268  9.900 /
%
%
%
\linethickness= 0.500pt
\setplotsymbol ({\thinlinefont .})
\putrule from 30.076  8.431 to 30.076 10.262
%
%
\plot 30.139 10.008 30.076 10.262 30.012 10.008 /
%
%
%
\linethickness= 0.500pt
\setplotsymbol ({\thinlinefont .})
\plot 30.033  8.124 29.652 10.050 /
%
%
\plot 29.764  9.813 29.652 10.050 29.639  9.788 /
%
%
%
\linethickness= 0.500pt
\setplotsymbol ({\thinlinefont .})
\plot 30.044  7.785 29.081  9.066 /
%
%
\plot 29.284  8.901 29.081  9.066 29.183  8.825 /
%
%
%
\linethickness= 0.500pt
\setplotsymbol ({\thinlinefont .})
\plot 30.044  7.446 29.652  7.690 /
%
%
\plot 29.902  7.610 29.652  7.690 29.835  7.502 /
%
%
%
\linethickness= 0.500pt
\setplotsymbol ({\thinlinefont .})
\plot 30.055  7.023 29.875  7.319 /
%
%
\plot 30.061  7.135 29.875  7.319 29.952  7.069 /
\linethickness= 0.500pt
\setplotsymbol ({\thinlinefont .})
%
%
%
\plot	30.055  6.208 30.055  6.526
 	30.052  6.605
	30.046  6.686
	30.036  6.767
	30.021  6.848
	30.003  6.931
	29.980  7.014
	29.953  7.098
	29.922  7.182
	29.886  7.268
	29.843  7.357
	29.794  7.449
	29.738  7.544
	29.676  7.642
	29.608  7.743
	29.533  7.847
	29.451  7.954
	29.371  8.066
	29.300  8.184
	29.237  8.309
	29.210  8.374
	29.184  8.440
	29.161  8.508
	29.140  8.577
	29.121  8.649
	29.104  8.721
	29.090  8.796
	29.077  8.872
	29.067  8.949
	29.060  9.029
	29.055  9.108
	29.056  9.184
	29.061  9.259
	29.070  9.332
	29.084  9.402
	29.103  9.470
	29.126  9.536
	29.154  9.600
	29.223  9.722
	29.311  9.834
	29.417  9.938
	29.541 10.034
	29.607 10.077
	29.672 10.116
	29.796 10.179
	29.914 10.224
	30.025 10.250
	30.131 10.257
	30.230 10.245
	30.323 10.215
	30.409 10.166
	30.485 10.108
	30.548 10.051
	30.631  9.936
	30.659  9.823
	30.631  9.711
	30.569  9.614
	30.495  9.543
	30.408  9.500
	30.309  9.484
	 /
\plot 30.309  9.484 30.097  9.478 /
\linethickness=0pt
\putrectangle corners at 29.015 11.282 and 30.704  5.156
\endpicture
\endfig

Now apply the global flow.  The interval where the vector field is not 
vertically up flows upwards more slowly and at the same time
it flows under and to both sides on the original line.  The result is
the twist illustrated in figure \figkey\twist.

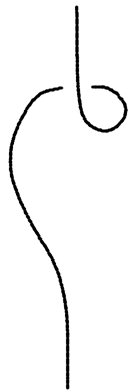
\fig{\twist: the effect of the global flow}
\beginpicture
\setcoordinatesystem units <.70000cm,.70000cm>
\unitlength=1.00000cm
\linethickness=1pt
\setplotsymbol ({\makebox(0,0)[l]{\tencirc\symbol{'160}}})
\setshadesymbol ({\thinlinefont .})
\setlinear
\linethickness=1pt
\setplotsymbol ({\makebox(0,0)[l]{\tencirc\symbol{'160}}})
%
%
%
\plot	33.500  4.917 33.495  5.854
 	33.493  5.969
	33.488  6.082
	33.480  6.193
	33.469  6.300
	33.455  6.405
	33.438  6.508
	33.418  6.607
	33.395  6.704
	33.369  6.799
	33.340  6.890
	33.308  6.979
	33.273  7.065
	33.235  7.149
	33.194  7.229
	33.150  7.308
	33.104  7.383
	33.056  7.457
	33.011  7.530
	32.967  7.602
	32.925  7.673
	32.885  7.744
	32.847  7.814
	32.810  7.883
	32.776  7.952
	32.743  8.020
	32.712  8.087
	32.683  8.153
	32.655  8.218
	32.630  8.283
	32.606  8.347
	32.564  8.473
	32.532  8.598
	32.512  8.724
	32.505  8.851
	32.506  8.914
	32.510  8.978
	32.517  9.043
	32.527  9.107
	32.540  9.172
	32.557  9.237
	32.576  9.302
	32.599  9.368
	32.625  9.434
	32.654  9.500
	32.715  9.624
	32.776  9.734
	32.837  9.829
	32.899  9.908
	33.022 10.022
	33.146 10.076
	 /
\plot 33.146 10.076 33.395 10.124 /
\linethickness=1pt
\setplotsymbol ({\makebox(0,0)[l]{\tencirc\symbol{'160}}})
%
%
%
\plot	33.934 10.145 34.098 10.140
 	34.180 10.126
	34.262 10.088
	34.345 10.028
	34.427  9.944
	34.489  9.848
	34.513  9.750
	34.497  9.649
	34.442  9.547
	34.364  9.460
	34.276  9.404
	34.178  9.381
	34.072  9.388
	33.969  9.423
	33.883  9.478
	33.813  9.555
	33.760  9.653
	33.739  9.717
	33.720  9.800
	33.704  9.902
	33.691 10.023
	33.685 10.091
	33.680 10.164
	33.676 10.242
	33.672 10.324
	33.669 10.411
	33.667 10.503
	33.666 10.599
	33.665 10.701
	 /
\plot 33.665 10.701 33.659 11.532 /
\linethickness=0pt
\putrectangle corners at 32.455 11.563 and 34.548  4.885
\endpicture
\endfig

\sh{2 in 4}

Now observe that the twist constructed above has two possible forms
depending on the slope of the original line.   So now consider the
surface in 4--space (with normal vector field) which is described by
the moving line as it changes slope in 3--space from one side of 
horizontal to
the other.  The vector field so described is not grounded since in 
the middle of the movement the line is horizontal and one vector
points vertically down.  But a small general
position shift moves this normal vector one side (ie, into the past
or the future) and makes the field grounded.  We can then draw
the end result of the isotopy provided by the compression theorem
as the sequence of pictures in figure \figkey\mtwist\ 
which describe an embedded 2--space in 4--space.

\fig{\mtwist: embedding of 2 in 4}
\beginpicture
\setcoordinatesystem units <0.4000cm,0.4000cm>
\unitlength=1.00000cm
\linethickness=1pt
\setplotsymbol ({\makebox(0,0)[l]{\tencirc\symbol{'160}}})
\setshadesymbol ({\thinlinefont .})
\setlinear
%
%
\linethickness=1pt
\setplotsymbol ({\makebox(0,0)[l]{\tencirc\symbol{'160}}})
\circulararc 90.101 degrees from 14.362 20.278 center at 14.366 19.321
%
%
\linethickness=1pt
\setplotsymbol ({\makebox(0,0)[l]{\tencirc\symbol{'160}}})
\circulararc 67.007 degrees from 15.399 19.251 center at 14.319 19.253
%
%
\linethickness=1pt
\setplotsymbol ({\makebox(0,0)[l]{\tencirc\symbol{'160}}})
\circulararc 56.099 degrees from 14.679 18.341 center at 14.412 19.299
%
%
\linethickness=1pt
\setplotsymbol ({\makebox(0,0)[l]{\tencirc\symbol{'160}}})
\circulararc 75.103 degrees from 13.420 18.976 center at 14.293 19.214
%
%
\linethickness=1pt
\setplotsymbol ({\makebox(0,0)[l]{\tencirc\symbol{'160}}})
\circulararc 173.987 degrees from 21.273 17.230 center at 21.235 19.356
%
%
\linethickness=1pt
\setplotsymbol ({\makebox(0,0)[l]{\tencirc\symbol{'160}}})
\circulararc 174.705 degrees from 20.870 21.484 center at 20.990 19.347
%
%
\linethickness=1pt
\setplotsymbol ({\makebox(0,0)[l]{\tencirc\symbol{'160}}})
\ellipticalarc axes ratio  0.383:0.383  360 degrees 
	from  8.490 19.579 center at  8.107 19.579
%
%
\linethickness= 0.500pt
\setplotsymbol ({\thinlinefont .})
\putrule from 18.881 19.558 to 19.082 19.558
%
%
\linethickness= 0.500pt
\setplotsymbol ({\thinlinefont .})
\plot 18.965 19.939 19.135 19.928 /
%
%
\linethickness= 0.500pt
\setplotsymbol ({\thinlinefont .})
\plot 19.188 20.468 19.378 20.405 /
%
%
\linethickness= 0.500pt
\setplotsymbol ({\thinlinefont .})
\plot 19.092 20.246 19.272 20.193 /
%
%
\linethickness= 0.500pt
\setplotsymbol ({\thinlinefont .})
\plot 19.389 20.701 19.558 20.616 /
%
%
\linethickness= 0.500pt
\setplotsymbol ({\thinlinefont .})
\plot 19.632 20.934 19.748 20.817 /
%
%
\linethickness= 0.500pt
\setplotsymbol ({\thinlinefont .})
\plot 19.897 21.114 19.981 20.987 /
%
%
\linethickness= 0.500pt
\setplotsymbol ({\thinlinefont .})
\plot 20.161 21.273 20.235 21.124 /
%
%
\linethickness= 0.500pt
\setplotsymbol ({\thinlinefont .})
\plot 20.436 21.400 20.489 21.209 /
%
%
\linethickness= 0.500pt
\setplotsymbol ({\thinlinefont .})
\plot 20.669 21.452 20.690 21.262 /
%
%
\linethickness= 0.500pt
\setplotsymbol ({\thinlinefont .})
\putrule from 20.860 21.463 to 20.860 21.251
%
%
\linethickness= 0.500pt
\setplotsymbol ({\thinlinefont .})
\putrule from 21.421 21.463 to 21.421 21.220
%
%
\linethickness= 0.500pt
\setplotsymbol ({\thinlinefont .})
\plot 21.717 21.400 21.685 21.188 /
%
%
\linethickness= 0.500pt
\setplotsymbol ({\thinlinefont .})
\plot 22.013 21.294 21.929 21.114 /
%
%
\linethickness= 0.500pt
\setplotsymbol ({\thinlinefont .})
\plot 23.315 19.801 23.315 19.801 /
%
%
\linethickness= 0.500pt
\setplotsymbol ({\thinlinefont .})
\plot 23.305 19.854 23.082 19.780 /
%
%
\linethickness= 0.500pt
\setplotsymbol ({\thinlinefont .})
\plot 22.924 20.627 22.733 20.479 /
%
%
\linethickness= 0.500pt
\setplotsymbol ({\thinlinefont .})
\plot 22.246 21.198 22.151 21.008 /
%
%
\linethickness= 0.500pt
\setplotsymbol ({\thinlinefont .})
\plot 22.511 21.008 22.394 20.828 /
%
%
\linethickness= 0.500pt
\setplotsymbol ({\thinlinefont .})
\plot 22.744 20.817 22.606 20.648 /
%
%
\linethickness= 0.500pt
\setplotsymbol ({\thinlinefont .})
\plot 23.061 20.405 22.849 20.278 /
%
%
\linethickness= 0.500pt
\setplotsymbol ({\thinlinefont .})
\plot 23.167 20.204 22.945 20.087 /
%
%
\linethickness= 0.500pt
\setplotsymbol ({\thinlinefont .})
\plot 23.336 19.558 23.114 19.537 /
%
%
\linethickness= 0.500pt
\setplotsymbol ({\thinlinefont .})
\plot  7.885 19.865  7.779 19.452 /
%
%
\linethickness= 0.500pt
\setplotsymbol ({\thinlinefont .})
\plot  8.096 19.928  7.906 19.283 /
%
%
\linethickness= 0.500pt
\setplotsymbol ({\thinlinefont .})
\plot  8.340 19.854  8.128 19.198 /
%
%
\linethickness= 0.500pt
\setplotsymbol ({\thinlinefont .})
\plot  8.456 19.653  8.371 19.325 /
\linethickness=1pt
\setplotsymbol ({\makebox(0,0)[l]{\tencirc\symbol{'160}}})
%
%
%
\plot	27.318 23.654 27.318 21.511
 	27.318 21.444
	27.318 21.379
	27.318 21.314
	27.318 21.251
	27.318 21.127
	27.319 21.006
	27.320 20.890
	27.320 20.778
	27.321 20.669
	27.323 20.565
	27.324 20.465
	27.325 20.368
	27.327 20.275
	27.329 20.187
	27.331 20.102
	27.333 20.021
	27.335 19.944
	27.338 19.872
	27.340 19.803
	27.343 19.738
	27.349 19.619
	27.355 19.517
	27.362 19.430
	27.378 19.304
	27.397 19.241
	27.496 19.131
	27.561 19.090
	27.635 19.057
	27.735 19.034
	27.801 19.025
	27.877 19.019
	27.964 19.015
	28.061 19.014
	28.168 19.014
	28.286 19.017
	28.406 19.020
	28.518 19.019
	28.623 19.016
	28.721 19.009
	28.811 19.000
	28.894 18.987
	28.971 18.972
	29.040 18.954
	29.158 18.902
	29.252 18.825
	29.320 18.723
	29.365 18.597
	29.374 18.524
	29.371 18.444
	29.355 18.357
	29.327 18.262
	29.286 18.160
	29.232 18.051
	29.166 17.935
	29.087 17.812
	29.001 17.689
	28.914 17.576
	28.825 17.472
	28.734 17.377
	28.642 17.292
	28.549 17.215
	28.454 17.148
	28.357 17.090
	28.263 17.040
	28.174 16.996
	28.092 16.960
	28.015 16.930
	27.945 16.906
	27.881 16.890
	27.770 16.876
	 /
\plot 27.770 16.876 27.572 16.876 /
\linethickness=1pt
\setplotsymbol ({\makebox(0,0)[l]{\tencirc\symbol{'160}}})
%
%
%
\plot	27.032 16.876 26.810 16.908
 	26.688 16.938
	26.618 16.964
	26.544 16.997
	26.464 17.037
	26.378 17.085
	26.287 17.139
	26.191 17.201
	26.094 17.270
	26.001 17.345
	25.913 17.426
	25.829 17.514
	25.750 17.608
	25.676 17.709
	25.605 17.817
	25.540 17.930
	25.481 18.043
	25.431 18.147
	25.390 18.243
	25.359 18.331
	25.337 18.410
	25.324 18.481
	25.325 18.597
	25.357 18.690
	25.416 18.770
	25.500 18.837
	25.611 18.891
	25.675 18.914
	25.746 18.934
	25.823 18.953
	25.907 18.969
	25.996 18.982
	26.092 18.994
	26.194 19.003
	26.302 19.010
	26.409 19.015
	26.507 19.018
	26.597 19.018
	26.679 19.017
	26.752 19.014
	26.817 19.010
	26.921 18.994
	27.002 18.961
	27.072 18.903
	27.129 18.818
	27.175 18.708
	27.193 18.619
	27.202 18.554
	27.210 18.476
	27.217 18.384
	27.223 18.278
	27.229 18.158
	27.232 18.094
	27.234 18.025
	27.237 17.954
	27.239 17.879
	27.241 17.800
	27.243 17.718
	27.245 17.633
	27.246 17.544
	27.248 17.452
	27.249 17.357
	27.250 17.258
	27.251 17.155
	27.252 17.050
	27.253 16.940
	27.254 16.828
	27.254 16.712
	27.254 16.593
	27.254 16.470
	 /
\plot 27.254 16.470 27.254 14.478 /
\linethickness=1pt
\setplotsymbol ({\makebox(0,0)[l]{\tencirc\symbol{'160}}})
%
%
%
\plot	14.516 14.446 14.516 16.589
 	14.516 16.656
	14.516 16.721
	14.516 16.786
	14.516 16.849
	14.515 16.973
	14.515 17.094
	14.514 17.210
	14.513 17.322
	14.512 17.431
	14.511 17.535
	14.510 17.635
	14.508 17.732
	14.507 17.825
	14.505 17.913
	14.503 17.998
	14.501 18.079
	14.499 18.156
	14.496 18.228
	14.494 18.297
	14.491 18.362
	14.485 18.481
	14.479 18.583
	14.471 18.670
	14.455 18.796
	14.437 18.860
	14.338 18.969
	14.273 19.010
	14.199 19.043
	14.098 19.066
	14.033 19.075
	13.957 19.081
	13.870 19.085
	13.773 19.086
	13.666 19.086
	13.548 19.083
	13.428 19.080
	13.316 19.081
	13.211 19.084
	13.113 19.091
	13.023 19.100
	12.939 19.113
	12.863 19.128
	12.794 19.146
	12.676 19.198
	12.582 19.275
	12.513 19.377
	12.469 19.503
	12.460 19.576
	12.463 19.656
	12.479 19.743
	12.507 19.838
	12.548 19.940
	12.601 20.049
	12.668 20.165
	12.747 20.288
	12.833 20.411
	12.920 20.524
	13.009 20.628
	13.100 20.723
	13.192 20.808
	13.285 20.885
	13.380 20.952
	13.477 21.010
	13.571 21.060
	13.660 21.104
	13.742 21.140
	13.818 21.170
	13.889 21.194
	13.953 21.210
	14.064 21.224
	 /
\plot 14.064 21.224 14.262 21.224 /
\linethickness=1pt
\setplotsymbol ({\makebox(0,0)[l]{\tencirc\symbol{'160}}})
%
%
%
\plot	14.802 21.224 15.024 21.192
 	15.146 21.162
	15.215 21.136
	15.290 21.103
	15.370 21.063
	15.456 21.015
	15.547 20.961
	15.643 20.899
	15.740 20.830
	15.833 20.755
	15.921 20.674
	16.004 20.586
	16.084 20.492
	16.158 20.391
	16.228 20.283
	16.294 20.170
	16.353 20.057
	16.403 19.953
	16.443 19.857
	16.475 19.769
	16.497 19.690
	16.510 19.619
	16.509 19.503
	16.476 19.410
	16.418 19.330
	16.334 19.263
	16.223 19.209
	16.158 19.186
	16.088 19.166
	16.010 19.147
	15.927 19.131
	15.838 19.118
	15.742 19.106
	15.640 19.097
	15.532 19.090
	15.425 19.085
	15.327 19.082
	15.237 19.082
	15.155 19.083
	15.082 19.086
	15.017 19.090
	14.913 19.106
	14.832 19.139
	14.762 19.197
	14.705 19.282
	14.659 19.392
	14.640 19.481
	14.632 19.546
	14.624 19.624
	14.617 19.716
	14.611 19.822
	14.605 19.942
	14.602 20.006
	14.599 20.075
	14.597 20.146
	14.595 20.221
	14.593 20.300
	14.591 20.382
	14.589 20.467
	14.587 20.556
	14.586 20.648
	14.585 20.743
	14.583 20.842
	14.582 20.945
	14.582 21.050
	14.581 21.160
	14.580 21.272
	14.580 21.388
	14.580 21.507
	14.580 21.630
	 /
\plot 14.580 21.630 14.580 23.622 /
\linethickness=1pt
\setplotsymbol ({\makebox(0,0)[l]{\tencirc\symbol{'160}}})
%
%
%
\plot	21.107 14.478 21.107 16.621
 	21.107 16.688
	21.107 16.753
	21.107 16.818
	21.107 16.881
	21.107 17.005
	21.106 17.125
	21.105 17.242
	21.105 17.354
	21.104 17.462
	21.102 17.567
	21.101 17.667
	21.100 17.764
	21.098 17.856
	21.096 17.945
	21.094 18.030
	21.092 18.110
	21.090 18.187
	21.088 18.260
	21.085 18.329
	21.082 18.394
	21.076 18.512
	21.070 18.615
	21.063 18.702
	21.047 18.828
	21.028 18.891
	20.929 19.001
	20.864 19.042
	20.790 19.074
	20.690 19.098
	20.624 19.106
	20.548 19.112
	20.461 19.116
	20.364 19.118
	20.257 19.117
	20.139 19.115
	20.020 19.112
	19.907 19.113
	19.802 19.116
	19.705 19.122
	19.614 19.132
	19.531 19.144
	19.454 19.160
	19.385 19.178
	19.267 19.230
	19.174 19.307
	19.105 19.408
	19.061 19.535
	19.051 19.608
	19.054 19.688
	19.070 19.775
	19.098 19.870
	19.139 19.971
	19.193 20.080
	19.259 20.197
	19.338 20.320
	19.424 20.442
	19.511 20.556
	19.600 20.660
	19.691 20.754
	19.783 20.840
	19.876 20.917
	19.972 20.984
	20.068 21.042
	20.162 21.092
	20.251 21.135
	20.333 21.172
	20.410 21.202
	20.480 21.226
	20.545 21.242
	20.655 21.256
	 /
\plot 20.655 21.256 20.853 21.256 /
\linethickness=1pt
\setplotsymbol ({\makebox(0,0)[l]{\tencirc\symbol{'160}}})
%
%
%
\plot	21.393 21.256 21.615 21.224
 	21.737 21.194
	21.807 21.168
	21.881 21.135
	21.961 21.094
	22.047 21.047
	22.138 20.992
	22.235 20.931
	22.332 20.862
	22.424 20.787
	22.512 20.706
	22.596 20.618
	22.675 20.523
	22.749 20.422
	22.820 20.315
	22.885 20.201
	22.944 20.089
	22.994 19.984
	23.035 19.888
	23.066 19.801
	23.088 19.722
	23.101 19.651
	23.100 19.535
	23.068 19.442
	23.009 19.362
	22.925 19.295
	22.814 19.240
	22.750 19.218
	22.679 19.197
	22.602 19.179
	22.518 19.163
	22.429 19.150
	22.333 19.138
	22.231 19.129
	22.123 19.122
	22.016 19.117
	21.918 19.114
	21.828 19.113
	21.746 19.114
	21.673 19.117
	21.608 19.122
	21.504 19.138
	21.423 19.170
	21.353 19.229
	21.296 19.313
	21.250 19.424
	21.232 19.513
	21.223 19.577
	21.216 19.656
	21.208 19.748
	21.202 19.854
	21.196 19.973
	21.193 20.038
	21.191 20.106
	21.188 20.178
	21.186 20.253
	21.184 20.332
	21.182 20.414
	21.180 20.499
	21.179 20.588
	21.177 20.680
	21.176 20.775
	21.175 20.874
	21.174 20.976
	21.173 21.082
	21.172 21.191
	21.172 21.304
	21.171 21.420
	21.171 21.539
	21.171 21.662
	 /
\plot 21.171 21.662 21.171 23.654 /
\linethickness=1pt
\setplotsymbol ({\makebox(0,0)[l]{\tencirc\symbol{'160}}})
%
%
%
\plot	 2.940 14.351  2.940 16.494
 	 2.940 16.561
	 2.940 16.626
	 2.940 16.691
	 2.940 16.754
	 2.939 16.878
	 2.939 16.998
	 2.938 17.115
	 2.937 17.227
	 2.936 17.335
	 2.935 17.440
	 2.934 17.540
	 2.932 17.637
	 2.931 17.729
	 2.929 17.818
	 2.927 17.903
	 2.925 17.983
	 2.923 18.060
	 2.920 18.133
	 2.918 18.202
	 2.915 18.267
	 2.909 18.385
	 2.903 18.488
	 2.895 18.575
	 2.879 18.701
	 2.861 18.764
	 2.761 18.874
	 2.697 18.915
	 2.623 18.947
	 2.522 18.971
	 2.457 18.979
	 2.380 18.985
	 2.294 18.989
	 2.197 18.991
	 2.089 18.990
	 1.972 18.988
	 1.852 18.985
	 1.740 18.986
	 1.635 18.989
	 1.537 18.995
	 1.447 19.005
	 1.363 19.017
	 1.287 19.033
	 1.218 19.051
	 1.100 19.103
	 1.006 19.180
	 0.937 19.281
	 0.893 19.408
	 0.884 19.481
	 0.887 19.561
	 0.902 19.648
	 0.931 19.743
	 0.972 19.844
	 1.025 19.953
	 1.092 20.070
	 1.171 20.193
	 1.256 20.315
	 1.344 20.429
	 1.433 20.533
	 1.523 20.627
	 1.616 20.713
	 1.709 20.790
	 1.804 20.857
	 1.901 20.915
	 1.995 20.965
	 2.084 21.008
	 2.166 21.045
	 2.242 21.075
	 2.313 21.099
	 2.377 21.115
	 2.488 21.129
	 /
\plot  2.488 21.129  2.686 21.129 /
\linethickness=1pt
\setplotsymbol ({\makebox(0,0)[l]{\tencirc\symbol{'160}}})
%
%
%
\plot	 3.226 21.129  3.448 21.097
 	 3.570 21.067
	 3.639 21.041
	 3.714 21.008
	 3.794 20.967
	 3.880 20.920
	 3.971 20.865
	 4.067 20.804
	 4.164 20.735
	 4.257 20.660
	 4.345 20.579
	 4.428 20.491
	 4.507 20.396
	 4.582 20.295
	 4.652 20.188
	 4.718 20.074
	 4.777 19.962
	 4.827 19.857
	 4.867 19.761
	 4.899 19.674
	 4.921 19.595
	 4.934 19.524
	 4.933 19.408
	 4.900 19.315
	 4.842 19.235
	 4.757 19.168
	 4.647 19.114
	 4.582 19.091
	 4.511 19.070
	 4.434 19.052
	 4.351 19.036
	 4.262 19.023
	 4.166 19.011
	 4.064 19.002
	 3.956 18.995
	 3.849 18.990
	 3.751 18.987
	 3.661 18.986
	 3.579 18.987
	 3.506 18.990
	 3.441 18.995
	 3.337 19.011
	 3.256 19.043
	 3.186 19.102
	 3.129 19.186
	 3.083 19.297
	 3.064 19.386
	 3.056 19.450
	 3.048 19.529
	 3.041 19.621
	 3.035 19.727
	 3.029 19.846
	 3.026 19.911
	 3.023 19.979
	 3.021 20.051
	 3.019 20.126
	 3.017 20.205
	 3.015 20.287
	 3.013 20.372
	 3.011 20.461
	 3.010 20.553
	 3.009 20.648
	 3.007 20.747
	 3.006 20.849
	 3.005 20.955
	 3.005 21.064
	 3.004 21.177
	 3.004 21.293
	 3.004 21.412
	 3.004 21.535
	 /
\plot  3.004 21.535  3.004 23.527 /
\linethickness=1pt
\setplotsymbol ({\makebox(0,0)[l]{\tencirc\symbol{'160}}})
%
%
%
\plot	 8.587 14.351  8.587 16.494
 	 8.587 16.561
	 8.587 16.626
	 8.587 16.691
	 8.587 16.754
	 8.587 16.878
	 8.586 16.998
	 8.585 17.115
	 8.585 17.227
	 8.584 17.335
	 8.582 17.440
	 8.581 17.540
	 8.580 17.637
	 8.578 17.729
	 8.576 17.818
	 8.574 17.903
	 8.572 17.983
	 8.570 18.060
	 8.567 18.133
	 8.565 18.202
	 8.562 18.267
	 8.556 18.385
	 8.550 18.488
	 8.543 18.575
	 8.527 18.701
	 8.508 18.764
	 8.409 18.874
	 8.344 18.915
	 8.270 18.947
	 8.170 18.971
	 8.104 18.979
	 8.028 18.985
	 7.941 18.989
	 7.844 18.991
	 7.737 18.990
	 7.619 18.988
	 7.500 18.985
	 7.387 18.986
	 7.282 18.989
	 7.184 18.995
	 7.094 19.005
	 7.011 19.017
	 6.934 19.033
	 6.865 19.051
	 6.747 19.103
	 6.653 19.180
	 6.585 19.281
	 6.541 19.408
	 6.531 19.481
	 6.534 19.561
	 6.550 19.648
	 6.578 19.743
	 6.619 19.844
	 6.673 19.953
	 6.739 20.070
	 6.818 20.193
	 6.904 20.315
	 6.991 20.429
	 7.080 20.533
	 7.171 20.627
	 7.263 20.713
	 7.356 20.790
	 7.451 20.857
	 7.548 20.915
	 7.642 20.965
	 7.731 21.008
	 7.813 21.045
	 7.890 21.075
	 7.960 21.099
	 8.024 21.115
	 8.135 21.129
	 /
\plot  8.135 21.129  8.333 21.129 /
\linethickness=1pt
\setplotsymbol ({\makebox(0,0)[l]{\tencirc\symbol{'160}}})
%
%
%
\plot	 8.873 21.129  9.095 21.097
 	 9.217 21.067
	 9.287 21.041
	 9.361 21.008
	 9.441 20.967
	 9.527 20.920
	 9.618 20.865
	 9.714 20.804
	 9.811 20.735
	 9.904 20.660
	 9.992 20.579
	10.076 20.491
	10.155 20.396
	10.229 20.295
	10.300 20.188
	10.365 20.074
	10.424 19.962
	10.474 19.857
	10.515 19.761
	10.546 19.674
	10.568 19.595
	10.581 19.524
	10.580 19.408
	10.548 19.315
	10.489 19.235
	10.405 19.168
	10.294 19.114
	10.230 19.091
	10.159 19.070
	10.082 19.052
	 9.998 19.036
	 9.909 19.023
	 9.813 19.011
	 9.711 19.002
	 9.603 18.995
	 9.496 18.990
	 9.398 18.987
	 9.308 18.986
	 9.226 18.987
	 9.153 18.990
	 9.088 18.995
	 8.984 19.011
	 8.903 19.043
	 8.833 19.102
	 8.776 19.186
	 8.730 19.297
	 8.712 19.386
	 8.703 19.450
	 8.695 19.529
	 8.688 19.621
	 8.682 19.727
	 8.676 19.846
	 8.673 19.911
	 8.671 19.979
	 8.668 20.051
	 8.666 20.126
	 8.664 20.205
	 8.662 20.287
	 8.660 20.372
	 8.659 20.461
	 8.657 20.553
	 8.656 20.648
	 8.655 20.747
	 8.654 20.849
	 8.653 20.955
	 8.652 21.064
	 8.652 21.177
	 8.651 21.293
	 8.651 21.412
	 8.651 21.535
	 /
\plot  8.651 21.535  8.651 23.527 /
%
%
\put{\SetFigFont{12}{14.4}{rm}{}\mtwist.1} [lB] <-4pt,0pt> at  8.223 13.049
%
%
\put{\SetFigFont{12}{14.4}{rm}{}\mtwist.2} [lB] <-4pt,0pt> at 14.161 13.081
%
%
\put{\SetFigFont{12}{14.4}{rm}{}\mtwist.3} [lB] <-4pt,0pt> at 20.923 13.018
%
%
\put{\SetFigFont{12}{14.4}{rm}{}\mtwist.4} [lB] <-4pt,2pt> at 27.083 12.954
%
%
\put{\SetFigFont{12}{14.4}{rm}{}\mtwist.0} [lB] <-4pt,0pt> at  2.762 13.018
\linethickness=0pt
\putrectangle corners at  0.830 23.685 and 29.428 12.897
\endpicture
\endfig

Notice that the pictures are not accurate
but diagrammatic;  the end result of the actual flow produces a surface
which, described as a moving picture starts like figure \twist\ and 
ends with a rotation of figure \twist.  However the combinatorial
structure of the pictures is the same as that produced by the flow.
Read the figure as a moving picture from left to right.  The picture 
is static
before time 0 (figure \mtwist.0); at time 1 a small disc appears 
(a 0--handle) and the boundary circle of this disc grows until it 
surrounds the main twist (figure \mtwist.3).  At this
point a 1--handle bridges across from the circle to the twist.  The
1--handle has been drawn very wide, because then the effect of the
bridge can be clearly seen as the replacement of the upwards twist
by the downwards twist in figure \mtwist.4.  Note that the 1--handle
in figure \mtwist.3 cancels with the 0--handle in figure \mtwist.1
so the topology of the surface is unaltered and the whole sequence
describes an embedded 2--plane in 4--space.  Indeed we can see
the small disc, whose boundary grows from times 1 to 3, as a little
finger pushed into and under the surface pointing into the past.

Notice that there is a triple point in the projected immersion
between \mtwist.2 and \mtwist.3 as the circle grows past the
double point.

\sh{3 in 5}

The embedding of 2 in 4 just constructed again has an asymmetry 
(corresponding to the choice of past or future for the general
position shift).  The pictures drawn were for the choice of shift
to the past.  The choice of shift to the future produces a similar
picture with the finger pointing to the right.  We can move from 
one set of pictures to the other by 
a similar construction illustrated in figure \figkey\xtwist, 
which should be
thought of as a moving sequence of moving pictures of 1 in 3 and
hence describes an embedding of 3 in 5.

\fig{\xtwist: embedding of 3 in 5}
\beginpicture
\setcoordinatesystem units <.350000cm,.350000cm>
\unitlength=1.00000cm
\linethickness=1pt
\setplotsymbol ({\makebox(0,0)[l]{\tencirc\symbol{'160}}})
\setshadesymbol ({\thinlinefont .})
\setlinear
%
%
\linethickness=1pt
\setplotsymbol ({\makebox(0,0)[l]{\tencirc\symbol{'160}}})
\circulararc 90.101 degrees from 14.302  9.451 center at 14.307  8.494
%
%
\linethickness=1pt
\setplotsymbol ({\makebox(0,0)[l]{\tencirc\symbol{'160}}})
\circulararc 67.007 degrees from 15.339  8.424 center at 14.260  8.426
%
%
\linethickness=1pt
\setplotsymbol ({\makebox(0,0)[l]{\tencirc\symbol{'160}}})
\circulararc 56.099 degrees from 14.620  7.514 center at 14.353  8.473
%
%
\linethickness=1pt
\setplotsymbol ({\makebox(0,0)[l]{\tencirc\symbol{'160}}})
\circulararc 75.103 degrees from 13.360  8.149 center at 14.233  8.387
%
%
\linethickness= 0.500pt
\setplotsymbol ({\thinlinefont .})
\circulararc 90.785 degrees from 14.323  9.669 center at 14.331  8.507
%
%
\linethickness= 0.500pt
\setplotsymbol ({\thinlinefont .})
\circulararc 84.590 degrees from 13.180  8.145 center at 14.191  8.310
%
%
\linethickness= 0.500pt
\setplotsymbol ({\thinlinefont .})
\circulararc 76.765 degrees from 15.551  8.452 center at 14.364  8.490
%
%
\linethickness= 0.500pt
\setplotsymbol ({\thinlinefont .})
\circulararc 70.628 degrees from 14.630  7.309 center at 14.500  8.344
%
%
\linethickness=1pt
\setplotsymbol ({\makebox(0,0)[l]{\tencirc\symbol{'160}}})
\ellipticalarc axes ratio  0.383:0.383  360 degrees 
	from  8.431  8.752 center at  8.048  8.752
%
%
\linethickness= 0.500pt
\setplotsymbol ({\thinlinefont .})
\ellipticalarc axes ratio  0.508:0.508  360 degrees 
	from  8.545  8.769 center at  8.037  8.769
%
%
\linethickness= 0.500pt
\setplotsymbol ({\thinlinefont .})
\plot  7.825  9.038  7.719  8.625 /
%
%
\linethickness= 0.500pt
\setplotsymbol ({\thinlinefont .})
\plot  8.037  9.102  7.846  8.456 /
%
%
\linethickness= 0.500pt
\setplotsymbol ({\thinlinefont .})
\plot  8.280  9.028  8.069  8.371 /
%
%
\linethickness= 0.500pt
\setplotsymbol ({\thinlinefont .})
\plot  8.397  8.826  8.312  8.498 /
\linethickness=1pt
\setplotsymbol ({\makebox(0,0)[l]{\tencirc\symbol{'160}}})
%
%
%
\plot	14.457  3.620 14.457  5.763
 	14.457  5.829
	14.457  5.895
	14.457  5.959
	14.457  6.023
	14.456  6.147
	14.456  6.267
	14.455  6.383
	14.454  6.495
	14.453  6.604
	14.452  6.708
	14.451  6.809
	14.449  6.905
	14.447  6.998
	14.446  7.086
	14.444  7.171
	14.442  7.252
	14.439  7.329
	14.437  7.402
	14.434  7.471
	14.432  7.536
	14.426  7.654
	14.419  7.756
	14.412  7.843
	14.396  7.969
	14.377  8.033
	14.278  8.142
	14.214  8.183
	14.139  8.216
	14.039  8.239
	13.973  8.248
	13.897  8.254
	13.811  8.258
	13.714  8.259
	13.606  8.259
	13.488  8.256
	13.369  8.254
	13.257  8.254
	13.152  8.258
	13.054  8.264
	12.963  8.273
	12.880  8.286
	12.804  8.301
	12.735  8.320
	12.617  8.372
	12.523  8.448
	12.454  8.550
	12.410  8.676
	12.400  8.749
	12.404  8.829
	12.419  8.917
	12.448  9.011
	12.489  9.113
	12.542  9.222
	12.608  9.338
	12.687  9.462
	12.773  9.584
	12.861  9.697
	12.950  9.801
	13.040  9.896
	13.132  9.982
	13.226 10.058
	13.321 10.125
	13.418 10.183
	13.512 10.233
	13.600 10.277
	13.683 10.314
	13.759 10.344
	13.830 10.367
	13.894 10.384
	14.005 10.397
	 /
\plot 14.005 10.397 14.203 10.397 /
\linethickness=1pt
\setplotsymbol ({\makebox(0,0)[l]{\tencirc\symbol{'160}}})
%
%
%
\plot	14.743 10.397 14.965 10.365
 	15.087 10.335
	15.156 10.309
	15.231 10.276
	15.311 10.236
	15.396 10.189
	15.487 10.134
	15.584 10.072
	15.681 10.004
	15.773  9.929
	15.862  9.847
	15.945  9.759
	16.024  9.665
	16.099  9.564
	16.169  9.457
	16.235  9.343
	16.294  9.230
	16.344  9.126
	16.384  9.030
	16.416  8.942
	16.438  8.863
	16.451  8.793
	16.450  8.676
	16.417  8.583
	16.359  8.503
	16.274  8.436
	16.164  8.382
	16.099  8.359
	16.028  8.339
	15.951  8.321
	15.868  8.305
	15.778  8.291
	15.683  8.280
	15.581  8.270
	15.473  8.263
	15.366  8.259
	15.267  8.256
	15.177  8.255
	15.096  8.256
	15.023  8.259
	14.958  8.264
	14.854  8.279
	14.772  8.312
	14.703  8.370
	14.645  8.455
	14.600  8.565
	14.581  8.654
	14.573  8.719
	14.565  8.797
	14.558  8.890
	14.551  8.995
	14.545  9.115
	14.543  9.180
	14.540  9.248
	14.538  9.320
	14.536  9.395
	14.533  9.473
	14.531  9.555
	14.530  9.640
	14.528  9.729
	14.527  9.821
	14.525  9.917
	14.524 10.016
	14.523 10.118
	14.522 10.224
	14.522 10.333
	14.521 10.445
	14.521 10.561
	14.520 10.681
	14.520 10.803
	 /
\plot 14.520 10.803 14.520 12.795 /
\linethickness=1pt
\setplotsymbol ({\makebox(0,0)[l]{\tencirc\symbol{'160}}})
%
%
%
\plot	 8.528  3.524  8.528  5.667
 	 8.528  5.734
	 8.528  5.799
	 8.528  5.864
	 8.528  5.927
	 8.527  6.051
	 8.527  6.172
	 8.526  6.288
	 8.525  6.400
	 8.524  6.509
	 8.523  6.613
	 8.522  6.713
	 8.520  6.810
	 8.519  6.903
	 8.517  6.991
	 8.515  7.076
	 8.513  7.157
	 8.511  7.234
	 8.508  7.306
	 8.506  7.375
	 8.503  7.440
	 8.497  7.559
	 8.491  7.661
	 8.483  7.748
	 8.467  7.874
	 8.449  7.938
	 8.349  8.047
	 8.285  8.088
	 8.211  8.121
	 8.110  8.144
	 8.045  8.153
	 7.968  8.159
	 7.882  8.163
	 7.785  8.164
	 7.677  8.164
	 7.560  8.161
	 7.440  8.158
	 7.328  8.159
	 7.223  8.162
	 7.125  8.169
	 7.035  8.178
	 6.951  8.191
	 6.875  8.206
	 6.806  8.224
	 6.688  8.276
	 6.594  8.353
	 6.525  8.455
	 6.481  8.581
	 6.472  8.654
	 6.475  8.734
	 6.490  8.821
	 6.519  8.916
	 6.560  9.018
	 6.613  9.127
	 6.680  9.243
	 6.759  9.366
	 6.844  9.489
	 6.932  9.602
	 7.021  9.706
	 7.111  9.801
	 7.204  9.886
	 7.297  9.963
	 7.392 10.030
	 7.489 10.088
	 7.583 10.138
	 7.672 10.182
	 7.754 10.218
	 7.830 10.248
	 7.901 10.272
	 7.965 10.288
	 8.076 10.302
	 /
\plot  8.076 10.302  8.274 10.302 /
\linethickness=1pt
\setplotsymbol ({\makebox(0,0)[l]{\tencirc\symbol{'160}}})
%
%
%
\plot	 8.814 10.302  9.036 10.270
 	 9.158 10.240
	 9.227 10.214
	 9.302 10.181
	 9.382 10.141
	 9.468 10.093
	 9.559 10.039
	 9.655  9.977
	 9.752  9.908
	 9.845  9.833
	 9.933  9.752
	10.016  9.664
	10.095  9.570
	10.170  9.469
	10.240  9.361
	10.306  9.248
	10.365  9.135
	10.415  9.031
	10.455  8.935
	10.487  8.847
	10.509  8.768
	10.522  8.697
	10.521  8.581
	10.488  8.488
	10.430  8.408
	10.345  8.341
	10.235  8.287
	10.170  8.264
	10.099  8.244
	10.022  8.225
	 9.939  8.209
	 9.850  8.196
	 9.754  8.184
	 9.652  8.175
	 9.544  8.168
	 9.437  8.163
	 9.339  8.160
	 9.249  8.160
	 9.167  8.161
	 9.094  8.164
	 9.029  8.168
	 8.925  8.184
	 8.844  8.217
	 8.774  8.275
	 8.717  8.360
	 8.671  8.470
	 8.652  8.559
	 8.644  8.624
	 8.636  8.702
	 8.629  8.794
	 8.623  8.900
	 8.617  9.020
	 8.614  9.084
	 8.611  9.153
	 8.609  9.224
	 8.607  9.299
	 8.605  9.378
	 8.603  9.460
	 8.601  9.545
	 8.599  9.634
	 8.598  9.726
	 8.597  9.821
	 8.595  9.920
	 8.594 10.023
	 8.593 10.128
	 8.593 10.238
	 8.592 10.350
	 8.592 10.466
	 8.592 10.585
	 8.592 10.708
	 /
\plot  8.592 10.708  8.592 12.700 /
%
%
\linethickness= 0.500pt
\setplotsymbol ({\thinlinefont .})
\circulararc 61.928 degrees from 20.769 18.855 center at 21.420 18.681
%
%
\linethickness= 0.500pt
\setplotsymbol ({\thinlinefont .})
\circulararc 22.620 degrees from 21.150 19.300 center at 21.388 18.839
%
%
\linethickness= 0.500pt
\setplotsymbol ({\thinlinefont .})
\circulararc 37.032 degrees from 21.954 19.236 center at 21.458 18.541
%
%
\linethickness= 0.500pt
\setplotsymbol ({\thinlinefont .})
\circulararc 126.859 degrees from 21.308 18.040 center at 21.486 18.691
%
%
\linethickness=1pt
\setplotsymbol ({\makebox(0,0)[l]{\tencirc\symbol{'160}}})
\circulararc 90.101 degrees from 14.362 20.278 center at 14.366 19.321
%
%
\linethickness=1pt
\setplotsymbol ({\makebox(0,0)[l]{\tencirc\symbol{'160}}})
\circulararc 67.007 degrees from 15.399 19.251 center at 14.319 19.253
%
%
\linethickness=1pt
\setplotsymbol ({\makebox(0,0)[l]{\tencirc\symbol{'160}}})
\circulararc 56.099 degrees from 14.679 18.341 center at 14.412 19.299
%
%
\linethickness=1pt
\setplotsymbol ({\makebox(0,0)[l]{\tencirc\symbol{'160}}})
\circulararc 75.103 degrees from 13.420 18.976 center at 14.293 19.214
%
%
\linethickness=1pt
\setplotsymbol ({\makebox(0,0)[l]{\tencirc\symbol{'160}}})
\circulararc 173.987 degrees from 21.273 17.230 center at 21.235 19.356
%
%
\linethickness=1pt
\setplotsymbol ({\makebox(0,0)[l]{\tencirc\symbol{'160}}})
\circulararc 174.705 degrees from 20.870 21.484 center at 20.990 19.347
%
%
\linethickness=1pt
\setplotsymbol ({\makebox(0,0)[l]{\tencirc\symbol{'160}}})
\circulararc 173.987 degrees from 21.213  6.403 center at 21.176  8.529
%
%
\linethickness=1pt
\setplotsymbol ({\makebox(0,0)[l]{\tencirc\symbol{'160}}})
\circulararc 174.705 degrees from 20.811 10.657 center at 20.931  8.521
%
%
\linethickness=1pt
\setplotsymbol ({\makebox(0,0)[l]{\tencirc\symbol{'160}}})
\circulararc 90.001 degrees from 28.029 -3.941 center at 28.025 -2.982
%
%
\linethickness=1pt
\setplotsymbol ({\makebox(0,0)[l]{\tencirc\symbol{'160}}})
\circulararc 66.846 degrees from 26.992 -2.915 center at 28.073 -2.915
%
%
\linethickness=1pt
\setplotsymbol ({\makebox(0,0)[l]{\tencirc\symbol{'160}}})
\circulararc 56.102 degrees from 27.711 -2.004 center at 27.980 -2.961
%
%
\linethickness=1pt
\setplotsymbol ({\makebox(0,0)[l]{\tencirc\symbol{'160}}})
\circulararc 75.195 degrees from 28.971 -2.639 center at 28.099 -2.877
%
%
\linethickness=1pt
\setplotsymbol ({\makebox(0,0)[l]{\tencirc\symbol{'160}}})
\circulararc 173.949 degrees from 21.118 -0.893 center at 21.156 -3.018
%
%
\linethickness=1pt
\setplotsymbol ({\makebox(0,0)[l]{\tencirc\symbol{'160}}})
\circulararc 174.785 degrees from 21.520 -5.148 center at 21.402 -3.010
%
%
\linethickness= 0.500pt
\setplotsymbol ({\thinlinefont .})
\circulararc 171.556 degrees from 20.769 10.844 center at 20.960  8.464
%
%
\linethickness= 0.500pt
\setplotsymbol ({\thinlinefont .})
\circulararc 168.985 degrees from 21.277  6.081 center at 21.122  8.518
%
%
\linethickness= 0.500pt
\setplotsymbol ({\thinlinefont .})
\circulararc 90.837 degrees from 28.103  6.890 center at 28.097  8.052
%
%
\linethickness= 0.500pt
\setplotsymbol ({\thinlinefont .})
\circulararc 76.762 degrees from 26.875  8.107 center at 28.063  8.069
%
%
\linethickness= 0.500pt
\setplotsymbol ({\thinlinefont .})
\circulararc 64.702 degrees from 29.246  8.414 center at 28.007  8.040
%
%
\linethickness= 0.500pt
\setplotsymbol ({\thinlinefont .})
\circulararc 63.437 degrees from 27.785  9.288 center at 28.036  8.151
%
%
\linethickness= 0.500pt
\setplotsymbol ({\thinlinefont .})
\ellipticalarc axes ratio  0.428:0.428  360 degrees 
	from 15.958 18.220 center at 15.530 18.220
%
%
\linethickness= 0.500pt
\setplotsymbol ({\thinlinefont .})
\ellipticalarc axes ratio  0.425:0.425  360 degrees 
	from 28.592 17.935 center at 28.166 17.935
%
%
\linethickness=1pt
\setplotsymbol ({\makebox(0,0)[l]{\tencirc\symbol{'160}}})
\ellipticalarc axes ratio  0.383:0.383  360 degrees 
	from  8.490 19.579 center at  8.107 19.579
%
%
\linethickness=1pt
\setplotsymbol ({\makebox(0,0)[l]{\tencirc\symbol{'160}}})
\ellipticalarc axes ratio  0.383:0.383  360 degrees 
	from 34.667 -3.243 center at 34.284 -3.243
%
%
\linethickness= 0.500pt
\setplotsymbol ({\thinlinefont .})
\ellipticalarc axes ratio  0.508:0.508  360 degrees 
	from 34.897  7.789 center at 34.389  7.789
%
%
\linethickness= 0.500pt
\setplotsymbol ({\thinlinefont .})
\putrule from 18.881 19.558 to 19.082 19.558
%
%
\linethickness= 0.500pt
\setplotsymbol ({\thinlinefont .})
\plot 18.965 19.939 19.135 19.928 /
%
%
\linethickness= 0.500pt
\setplotsymbol ({\thinlinefont .})
\plot 19.188 20.468 19.378 20.405 /
%
%
\linethickness= 0.500pt
\setplotsymbol ({\thinlinefont .})
\plot 19.092 20.246 19.272 20.193 /
%
%
\linethickness= 0.500pt
\setplotsymbol ({\thinlinefont .})
\plot 19.389 20.701 19.558 20.616 /
%
%
\linethickness= 0.500pt
\setplotsymbol ({\thinlinefont .})
\plot 19.632 20.934 19.748 20.817 /
%
%
\linethickness= 0.500pt
\setplotsymbol ({\thinlinefont .})
\plot 19.897 21.114 19.981 20.987 /
%
%
\linethickness= 0.500pt
\setplotsymbol ({\thinlinefont .})
\plot 20.161 21.273 20.235 21.124 /
%
%
\linethickness= 0.500pt
\setplotsymbol ({\thinlinefont .})
\plot 20.436 21.400 20.489 21.209 /
%
%
\linethickness= 0.500pt
\setplotsymbol ({\thinlinefont .})
\plot 20.669 21.452 20.690 21.262 /
%
%
\linethickness= 0.500pt
\setplotsymbol ({\thinlinefont .})
\putrule from 20.860 21.463 to 20.860 21.251
%
%
\linethickness= 0.500pt
\setplotsymbol ({\thinlinefont .})
\putrule from 21.421 21.463 to 21.421 21.220
%
%
\linethickness= 0.500pt
\setplotsymbol ({\thinlinefont .})
\plot 21.717 21.400 21.685 21.188 /
%
%
\linethickness= 0.500pt
\setplotsymbol ({\thinlinefont .})
\plot 22.013 21.294 21.929 21.114 /
%
%
\linethickness= 0.500pt
\setplotsymbol ({\thinlinefont .})
\plot 23.315 19.801 23.315 19.801 /
%
%
\linethickness= 0.500pt
\setplotsymbol ({\thinlinefont .})
\plot 23.305 19.854 23.082 19.780 /
%
%
\linethickness= 0.500pt
\setplotsymbol ({\thinlinefont .})
\plot 22.924 20.627 22.733 20.479 /
%
%
\linethickness= 0.500pt
\setplotsymbol ({\thinlinefont .})
\plot 22.246 21.198 22.151 21.008 /
%
%
\linethickness= 0.500pt
\setplotsymbol ({\thinlinefont .})
\plot 22.511 21.008 22.394 20.828 /
%
%
\linethickness= 0.500pt
\setplotsymbol ({\thinlinefont .})
\plot 22.744 20.817 22.606 20.648 /
%
%
\linethickness= 0.500pt
\setplotsymbol ({\thinlinefont .})
\plot 23.061 20.405 22.849 20.278 /
%
%
\linethickness= 0.500pt
\setplotsymbol ({\thinlinefont .})
\plot 23.167 20.204 22.945 20.087 /
%
%
\linethickness= 0.500pt
\setplotsymbol ({\thinlinefont .})
\plot 23.336 19.558 23.114 19.537 /
%
%
\linethickness= 0.500pt
\setplotsymbol ({\thinlinefont .})
\plot  7.885 19.865  7.779 19.452 /
%
%
\linethickness= 0.500pt
\setplotsymbol ({\thinlinefont .})
\plot  8.096 19.928  7.906 19.283 /
%
%
\linethickness= 0.500pt
\setplotsymbol ({\thinlinefont .})
\plot  8.340 19.854  8.128 19.198 /
%
%
\linethickness= 0.500pt
\setplotsymbol ({\thinlinefont .})
\plot  8.456 19.653  8.371 19.325 /
%
%
\linethickness= 0.500pt
\setplotsymbol ({\thinlinefont .})
\putrule from 18.821  8.731 to 19.022  8.731
%
%
\linethickness= 0.500pt
\setplotsymbol ({\thinlinefont .})
\plot 18.906  9.112 19.075  9.102 /
%
%
\linethickness= 0.500pt
\setplotsymbol ({\thinlinefont .})
\plot 19.128  9.641 19.319  9.578 /
%
%
\linethickness= 0.500pt
\setplotsymbol ({\thinlinefont .})
\plot 19.033  9.419 19.213  9.366 /
%
%
\linethickness= 0.500pt
\setplotsymbol ({\thinlinefont .})
\plot 19.329  9.874 19.499  9.790 /
%
%
\linethickness= 0.500pt
\setplotsymbol ({\thinlinefont .})
\plot 19.573 10.107 19.689  9.991 /
%
%
\linethickness= 0.500pt
\setplotsymbol ({\thinlinefont .})
\plot 19.837 10.287 19.922 10.160 /
%
%
\linethickness= 0.500pt
\setplotsymbol ({\thinlinefont .})
\plot 20.102 10.446 20.176 10.298 /
%
%
\linethickness= 0.500pt
\setplotsymbol ({\thinlinefont .})
\plot 20.377 10.573 20.430 10.382 /
%
%
\linethickness= 0.500pt
\setplotsymbol ({\thinlinefont .})
\plot 20.610 10.626 20.631 10.435 /
%
%
\linethickness= 0.500pt
\setplotsymbol ({\thinlinefont .})
\putrule from 20.800 10.636 to 20.800 10.425
%
%
\linethickness= 0.500pt
\setplotsymbol ({\thinlinefont .})
\putrule from 21.361 10.636 to 21.361 10.393
%
%
\linethickness= 0.500pt
\setplotsymbol ({\thinlinefont .})
\plot 21.658 10.573 21.626 10.361 /
%
%
\linethickness= 0.500pt
\setplotsymbol ({\thinlinefont .})
\plot 21.954 10.467 21.869 10.287 /
%
%
\linethickness= 0.500pt
\setplotsymbol ({\thinlinefont .})
\plot 23.256  8.975 23.256  8.975 /
%
%
\linethickness= 0.500pt
\setplotsymbol ({\thinlinefont .})
\plot 23.245  9.028 23.023  8.954 /
%
%
\linethickness= 0.500pt
\setplotsymbol ({\thinlinefont .})
\plot 22.864  9.800 22.674  9.652 /
%
%
\linethickness= 0.500pt
\setplotsymbol ({\thinlinefont .})
\plot 22.187 10.372 22.092 10.181 /
%
%
\linethickness= 0.500pt
\setplotsymbol ({\thinlinefont .})
\plot 22.451 10.181 22.335 10.001 /
%
%
\linethickness= 0.500pt
\setplotsymbol ({\thinlinefont .})
\plot 22.684  9.991 22.547  9.821 /
%
%
\linethickness= 0.500pt
\setplotsymbol ({\thinlinefont .})
\plot 23.002  9.578 22.790  9.451 /
%
%
\linethickness= 0.500pt
\setplotsymbol ({\thinlinefont .})
\plot 23.108  9.377 22.885  9.260 /
%
%
\linethickness= 0.500pt
\setplotsymbol ({\thinlinefont .})
\plot 23.277  8.731 23.055  8.710 /
%
%
\linethickness= 0.500pt
\setplotsymbol ({\thinlinefont .})
\putrule from 23.510 -3.222 to 23.309 -3.222
%
%
\linethickness= 0.500pt
\setplotsymbol ({\thinlinefont .})
\plot 23.425 -3.603 23.256 -3.592 /
%
%
\linethickness= 0.500pt
\setplotsymbol ({\thinlinefont .})
\plot 23.203 -4.132 23.012 -4.068 /
%
%
\linethickness= 0.500pt
\setplotsymbol ({\thinlinefont .})
\plot 23.298 -3.909 23.118 -3.857 /
%
%
\linethickness= 0.500pt
\setplotsymbol ({\thinlinefont .})
\plot 23.002 -4.365 22.832 -4.280 /
%
%
\linethickness= 0.500pt
\setplotsymbol ({\thinlinefont .})
\plot 22.758 -4.597 22.642 -4.481 /
%
%
\linethickness= 0.500pt
\setplotsymbol ({\thinlinefont .})
\plot 22.494 -4.777 22.409 -4.650 /
%
%
\linethickness= 0.500pt
\setplotsymbol ({\thinlinefont .})
\plot 22.229 -4.936 22.155 -4.788 /
%
%
\linethickness= 0.500pt
\setplotsymbol ({\thinlinefont .})
\plot 21.954 -5.063 21.901 -4.873 /
%
%
\linethickness= 0.500pt
\setplotsymbol ({\thinlinefont .})
\plot 21.721 -5.116 21.700 -4.925 /
%
%
\linethickness= 0.500pt
\setplotsymbol ({\thinlinefont .})
\putrule from 21.531 -5.127 to 21.531 -4.915
%
%
\linethickness= 0.500pt
\setplotsymbol ({\thinlinefont .})
\putrule from 20.970 -5.127 to 20.970 -4.883
%
%
\linethickness= 0.500pt
\setplotsymbol ({\thinlinefont .})
\plot 20.673 -5.063 20.705 -4.851 /
%
%
\linethickness= 0.500pt
\setplotsymbol ({\thinlinefont .})
\plot 20.377 -4.957 20.462 -4.777 /
%
%
\linethickness= 0.500pt
\setplotsymbol ({\thinlinefont .})
\plot 19.075 -3.465 19.075 -3.465 /
%
%
\linethickness= 0.500pt
\setplotsymbol ({\thinlinefont .})
\plot 19.086 -3.518 19.308 -3.444 /
%
%
\linethickness= 0.500pt
\setplotsymbol ({\thinlinefont .})
\plot 19.467 -4.290 19.657 -4.142 /
%
%
\linethickness= 0.500pt
\setplotsymbol ({\thinlinefont .})
\plot 20.144 -4.862 20.240 -4.671 /
%
%
\linethickness= 0.500pt
\setplotsymbol ({\thinlinefont .})
\plot 19.880 -4.671 19.996 -4.492 /
%
%
\linethickness= 0.500pt
\setplotsymbol ({\thinlinefont .})
\plot 19.647 -4.481 19.784 -4.312 /
%
%
\linethickness= 0.500pt
\setplotsymbol ({\thinlinefont .})
\plot 19.329 -4.068 19.541 -3.941 /
%
%
\linethickness= 0.500pt
\setplotsymbol ({\thinlinefont .})
\plot 19.224 -3.867 19.446 -3.751 /
%
%
\linethickness= 0.500pt
\setplotsymbol ({\thinlinefont .})
\plot 19.054 -3.222 19.276 -3.200 /
%
%
\linethickness= 0.500pt
\setplotsymbol ({\thinlinefont .})
\plot 34.506 -3.528 34.612 -3.116 /
%
%
\linethickness= 0.500pt
\setplotsymbol ({\thinlinefont .})
\plot 34.294 -3.592 34.485 -2.946 /
%
%
\linethickness= 0.500pt
\setplotsymbol ({\thinlinefont .})
\plot 34.051 -3.518 34.262 -2.862 /
%
%
\linethickness= 0.500pt
\setplotsymbol ({\thinlinefont .})
\plot 33.934 -3.317 34.019 -2.989 /
\linethickness=1pt
\setplotsymbol ({\makebox(0,0)[l]{\tencirc\symbol{'160}}})
%
%
%
\plot	14.516 14.446 14.516 16.589
 	14.516 16.656
	14.516 16.721
	14.516 16.786
	14.516 16.849
	14.515 16.973
	14.515 17.094
	14.514 17.210
	14.513 17.322
	14.512 17.431
	14.511 17.535
	14.510 17.635
	14.508 17.732
	14.507 17.825
	14.505 17.913
	14.503 17.998
	14.501 18.079
	14.499 18.156
	14.496 18.228
	14.494 18.297
	14.491 18.362
	14.485 18.481
	14.479 18.583
	14.471 18.670
	14.455 18.796
	14.437 18.860
	14.338 18.969
	14.273 19.010
	14.199 19.043
	14.098 19.066
	14.033 19.075
	13.957 19.081
	13.870 19.085
	13.773 19.086
	13.666 19.086
	13.548 19.083
	13.428 19.080
	13.316 19.081
	13.211 19.084
	13.113 19.091
	13.023 19.100
	12.939 19.113
	12.863 19.128
	12.794 19.146
	12.676 19.198
	12.582 19.275
	12.513 19.377
	12.469 19.503
	12.460 19.576
	12.463 19.656
	12.479 19.743
	12.507 19.838
	12.548 19.940
	12.601 20.049
	12.668 20.165
	12.747 20.288
	12.833 20.411
	12.920 20.524
	13.009 20.628
	13.100 20.723
	13.192 20.808
	13.285 20.885
	13.380 20.952
	13.477 21.010
	13.571 21.060
	13.660 21.104
	13.742 21.140
	13.818 21.170
	13.889 21.194
	13.953 21.210
	14.064 21.224
	 /
\plot 14.064 21.224 14.262 21.224 /
\linethickness=1pt
\setplotsymbol ({\makebox(0,0)[l]{\tencirc\symbol{'160}}})
%
%
%
\plot	14.802 21.224 15.024 21.192
 	15.146 21.162
	15.215 21.136
	15.290 21.103
	15.370 21.063
	15.456 21.015
	15.547 20.961
	15.643 20.899
	15.740 20.830
	15.833 20.755
	15.921 20.674
	16.004 20.586
	16.084 20.492
	16.158 20.391
	16.228 20.283
	16.294 20.170
	16.353 20.057
	16.403 19.953
	16.443 19.857
	16.475 19.769
	16.497 19.690
	16.510 19.619
	16.509 19.503
	16.476 19.410
	16.418 19.330
	16.334 19.263
	16.223 19.209
	16.158 19.186
	16.088 19.166
	16.010 19.147
	15.927 19.131
	15.838 19.118
	15.742 19.106
	15.640 19.097
	15.532 19.090
	15.425 19.085
	15.327 19.082
	15.237 19.082
	15.155 19.083
	15.082 19.086
	15.017 19.090
	14.913 19.106
	14.832 19.139
	14.762 19.197
	14.705 19.282
	14.659 19.392
	14.640 19.481
	14.632 19.546
	14.624 19.624
	14.617 19.716
	14.611 19.822
	14.605 19.942
	14.602 20.006
	14.599 20.075
	14.597 20.146
	14.595 20.221
	14.593 20.300
	14.591 20.382
	14.589 20.467
	14.587 20.556
	14.586 20.648
	14.585 20.743
	14.583 20.842
	14.582 20.945
	14.582 21.050
	14.581 21.160
	14.580 21.272
	14.580 21.388
	14.580 21.507
	14.580 21.630
	 /
\plot 14.580 21.630 14.580 23.622 /
\linethickness=1pt
\setplotsymbol ({\makebox(0,0)[l]{\tencirc\symbol{'160}}})
%
%
%
\plot	21.107 14.478 21.107 16.621
 	21.107 16.688
	21.107 16.753
	21.107 16.818
	21.107 16.881
	21.107 17.005
	21.106 17.125
	21.105 17.242
	21.105 17.354
	21.104 17.462
	21.102 17.567
	21.101 17.667
	21.100 17.764
	21.098 17.856
	21.096 17.945
	21.094 18.030
	21.092 18.110
	21.090 18.187
	21.088 18.260
	21.085 18.329
	21.082 18.394
	21.076 18.512
	21.070 18.615
	21.063 18.702
	21.047 18.828
	21.028 18.891
	20.929 19.001
	20.864 19.042
	20.790 19.074
	20.690 19.098
	20.624 19.106
	20.548 19.112
	20.461 19.116
	20.364 19.118
	20.257 19.117
	20.139 19.115
	20.020 19.112
	19.907 19.113
	19.802 19.116
	19.705 19.122
	19.614 19.132
	19.531 19.144
	19.454 19.160
	19.385 19.178
	19.267 19.230
	19.174 19.307
	19.105 19.408
	19.061 19.535
	19.051 19.608
	19.054 19.688
	19.070 19.775
	19.098 19.870
	19.139 19.971
	19.193 20.080
	19.259 20.197
	19.338 20.320
	19.424 20.442
	19.511 20.556
	19.600 20.660
	19.691 20.754
	19.783 20.840
	19.876 20.917
	19.972 20.984
	20.068 21.042
	20.162 21.092
	20.251 21.135
	20.333 21.172
	20.410 21.202
	20.480 21.226
	20.545 21.242
	20.655 21.256
	 /
\plot 20.655 21.256 20.853 21.256 /
\linethickness=1pt
\setplotsymbol ({\makebox(0,0)[l]{\tencirc\symbol{'160}}})
%
%
%
\plot	21.393 21.256 21.615 21.224
 	21.737 21.194
	21.807 21.168
	21.881 21.135
	21.961 21.094
	22.047 21.047
	22.138 20.992
	22.235 20.931
	22.332 20.862
	22.424 20.787
	22.512 20.706
	22.596 20.618
	22.675 20.523
	22.749 20.422
	22.820 20.315
	22.885 20.201
	22.944 20.089
	22.994 19.984
	23.035 19.888
	23.066 19.801
	23.088 19.722
	23.101 19.651
	23.100 19.535
	23.068 19.442
	23.009 19.362
	22.925 19.295
	22.814 19.240
	22.750 19.218
	22.679 19.197
	22.602 19.179
	22.518 19.163
	22.429 19.150
	22.333 19.138
	22.231 19.129
	22.123 19.122
	22.016 19.117
	21.918 19.114
	21.828 19.113
	21.746 19.114
	21.673 19.117
	21.608 19.122
	21.504 19.138
	21.423 19.170
	21.353 19.229
	21.296 19.313
	21.250 19.424
	21.232 19.513
	21.223 19.577
	21.216 19.656
	21.208 19.748
	21.202 19.854
	21.196 19.973
	21.193 20.038
	21.191 20.106
	21.188 20.178
	21.186 20.253
	21.184 20.332
	21.182 20.414
	21.180 20.499
	21.179 20.588
	21.177 20.680
	21.176 20.775
	21.175 20.874
	21.174 20.976
	21.173 21.082
	21.172 21.191
	21.172 21.304
	21.171 21.420
	21.171 21.539
	21.171 21.662
	 /
\plot 21.171 21.662 21.171 23.654 /
\linethickness=1pt
\setplotsymbol ({\makebox(0,0)[l]{\tencirc\symbol{'160}}})
%
%
%
\plot	 8.587 14.351  8.587 16.494
 	 8.587 16.561
	 8.587 16.626
	 8.587 16.691
	 8.587 16.754
	 8.587 16.878
	 8.586 16.998
	 8.585 17.115
	 8.585 17.227
	 8.584 17.335
	 8.582 17.440
	 8.581 17.540
	 8.580 17.637
	 8.578 17.729
	 8.576 17.818
	 8.574 17.903
	 8.572 17.983
	 8.570 18.060
	 8.567 18.133
	 8.565 18.202
	 8.562 18.267
	 8.556 18.385
	 8.550 18.488
	 8.543 18.575
	 8.527 18.701
	 8.508 18.764
	 8.409 18.874
	 8.344 18.915
	 8.270 18.947
	 8.170 18.971
	 8.104 18.979
	 8.028 18.985
	 7.941 18.989
	 7.844 18.991
	 7.737 18.990
	 7.619 18.988
	 7.500 18.985
	 7.387 18.986
	 7.282 18.989
	 7.184 18.995
	 7.094 19.005
	 7.011 19.017
	 6.934 19.033
	 6.865 19.051
	 6.747 19.103
	 6.653 19.180
	 6.585 19.281
	 6.541 19.408
	 6.531 19.481
	 6.534 19.561
	 6.550 19.648
	 6.578 19.743
	 6.619 19.844
	 6.673 19.953
	 6.739 20.070
	 6.818 20.193
	 6.904 20.315
	 6.991 20.429
	 7.080 20.533
	 7.171 20.627
	 7.263 20.713
	 7.356 20.790
	 7.451 20.857
	 7.548 20.915
	 7.642 20.965
	 7.731 21.008
	 7.813 21.045
	 7.890 21.075
	 7.960 21.099
	 8.024 21.115
	 8.135 21.129
	 /
\plot  8.135 21.129  8.333 21.129 /
\linethickness=1pt
\setplotsymbol ({\makebox(0,0)[l]{\tencirc\symbol{'160}}})
%
%
%
\plot	 8.873 21.129  9.095 21.097
 	 9.217 21.067
	 9.287 21.041
	 9.361 21.008
	 9.441 20.967
	 9.527 20.920
	 9.618 20.865
	 9.714 20.804
	 9.811 20.735
	 9.904 20.660
	 9.992 20.579
	10.076 20.491
	10.155 20.396
	10.229 20.295
	10.300 20.188
	10.365 20.074
	10.424 19.962
	10.474 19.857
	10.515 19.761
	10.546 19.674
	10.568 19.595
	10.581 19.524
	10.580 19.408
	10.548 19.315
	10.489 19.235
	10.405 19.168
	10.294 19.114
	10.230 19.091
	10.159 19.070
	10.082 19.052
	 9.998 19.036
	 9.909 19.023
	 9.813 19.011
	 9.711 19.002
	 9.603 18.995
	 9.496 18.990
	 9.398 18.987
	 9.308 18.986
	 9.226 18.987
	 9.153 18.990
	 9.088 18.995
	 8.984 19.011
	 8.903 19.043
	 8.833 19.102
	 8.776 19.186
	 8.730 19.297
	 8.712 19.386
	 8.703 19.450
	 8.695 19.529
	 8.688 19.621
	 8.682 19.727
	 8.676 19.846
	 8.673 19.911
	 8.671 19.979
	 8.668 20.051
	 8.666 20.126
	 8.664 20.205
	 8.662 20.287
	 8.660 20.372
	 8.659 20.461
	 8.657 20.553
	 8.656 20.648
	 8.655 20.747
	 8.654 20.849
	 8.653 20.955
	 8.652 21.064
	 8.652 21.177
	 8.651 21.293
	 8.651 21.412
	 8.651 21.535
	 /
\plot  8.651 21.535  8.651 23.527 /
\linethickness=1pt
\setplotsymbol ({\makebox(0,0)[l]{\tencirc\symbol{'160}}})
%
%
%
\plot	27.318 23.654 27.318 21.511
 	27.318 21.444
	27.318 21.379
	27.318 21.314
	27.318 21.251
	27.318 21.127
	27.319 21.006
	27.320 20.890
	27.320 20.778
	27.321 20.669
	27.323 20.565
	27.324 20.465
	27.325 20.368
	27.327 20.275
	27.329 20.187
	27.331 20.102
	27.333 20.021
	27.335 19.944
	27.338 19.872
	27.340 19.803
	27.343 19.738
	27.349 19.619
	27.355 19.517
	27.362 19.430
	27.378 19.304
	27.397 19.241
	27.496 19.131
	27.561 19.090
	27.635 19.057
	27.735 19.034
	27.801 19.025
	27.877 19.019
	27.964 19.015
	28.061 19.014
	28.168 19.014
	28.286 19.017
	28.406 19.020
	28.518 19.019
	28.623 19.016
	28.721 19.009
	28.811 19.000
	28.894 18.987
	28.971 18.972
	29.040 18.954
	29.158 18.902
	29.252 18.825
	29.320 18.723
	29.365 18.597
	29.374 18.524
	29.371 18.444
	29.355 18.357
	29.327 18.262
	29.286 18.160
	29.232 18.051
	29.166 17.935
	29.087 17.812
	29.001 17.689
	28.914 17.576
	28.825 17.472
	28.734 17.377
	28.642 17.292
	28.549 17.215
	28.454 17.148
	28.357 17.090
	28.263 17.040
	28.174 16.996
	28.092 16.960
	28.015 16.930
	27.945 16.906
	27.881 16.890
	27.770 16.876
	 /
\plot 27.770 16.876 27.572 16.876 /
\linethickness=1pt
\setplotsymbol ({\makebox(0,0)[l]{\tencirc\symbol{'160}}})
%
%
%
\plot	27.032 16.876 26.810 16.908
 	26.688 16.938
	26.618 16.964
	26.544 16.997
	26.464 17.037
	26.378 17.085
	26.287 17.139
	26.191 17.201
	26.094 17.270
	26.001 17.345
	25.913 17.426
	25.829 17.514
	25.750 17.608
	25.676 17.709
	25.605 17.817
	25.540 17.930
	25.481 18.043
	25.431 18.147
	25.390 18.243
	25.359 18.331
	25.337 18.410
	25.324 18.481
	25.325 18.597
	25.357 18.690
	25.416 18.770
	25.500 18.837
	25.611 18.891
	25.675 18.914
	25.746 18.934
	25.823 18.953
	25.907 18.969
	25.996 18.982
	26.092 18.994
	26.194 19.003
	26.302 19.010
	26.409 19.015
	26.507 19.018
	26.597 19.018
	26.679 19.017
	26.752 19.014
	26.817 19.010
	26.921 18.994
	27.002 18.961
	27.072 18.903
	27.129 18.818
	27.175 18.708
	27.193 18.619
	27.202 18.554
	27.210 18.476
	27.217 18.384
	27.223 18.278
	27.229 18.158
	27.232 18.094
	27.234 18.025
	27.237 17.954
	27.239 17.879
	27.241 17.800
	27.243 17.718
	27.245 17.633
	27.246 17.544
	27.248 17.452
	27.249 17.357
	27.250 17.258
	27.251 17.155
	27.252 17.050
	27.253 16.940
	27.254 16.828
	27.254 16.712
	27.254 16.593
	27.254 16.470
	 /
\plot 27.254 16.470 27.254 14.478 /
\linethickness=1pt
\setplotsymbol ({\makebox(0,0)[l]{\tencirc\symbol{'160}}})
%
%
%
\plot	34.148 23.685 34.148 21.542
 	34.148 21.476
	34.148 21.410
	34.148 21.346
	34.148 21.282
	34.149 21.158
	34.149 21.038
	34.150 20.922
	34.151 20.810
	34.152 20.701
	34.153 20.597
	34.154 20.496
	34.156 20.400
	34.158 20.307
	34.159 20.219
	34.161 20.134
	34.163 20.053
	34.166 19.976
	34.168 19.903
	34.171 19.834
	34.173 19.769
	34.179 19.651
	34.186 19.549
	34.193 19.462
	34.209 19.336
	34.228 19.272
	34.327 19.163
	34.391 19.122
	34.466 19.089
	34.566 19.066
	34.632 19.057
	34.708 19.051
	34.794 19.047
	34.891 19.046
	34.999 19.046
	35.117 19.049
	35.236 19.051
	35.348 19.051
	35.453 19.047
	35.551 19.041
	35.642 19.032
	35.725 19.019
	35.801 19.004
	35.870 18.985
	35.988 18.933
	36.082 18.857
	36.151 18.755
	36.195 18.629
	36.205 18.556
	36.201 18.476
	36.186 18.388
	36.157 18.294
	36.116 18.192
	36.063 18.083
	35.997 17.967
	35.918 17.844
	35.832 17.721
	35.744 17.608
	35.655 17.504
	35.565 17.409
	35.473 17.323
	35.379 17.247
	35.284 17.180
	35.187 17.122
	35.093 17.072
	35.005 17.028
	34.922 16.991
	34.846 16.961
	34.775 16.938
	34.711 16.921
	34.600 16.908
	 /
\plot 34.600 16.908 34.402 16.908 /
\linethickness=1pt
\setplotsymbol ({\makebox(0,0)[l]{\tencirc\symbol{'160}}})
%
%
%
\plot	33.862 16.908 33.640 16.940
 	33.518 16.970
	33.449 16.996
	33.374 17.029
	33.294 17.069
	33.209 17.116
	33.118 17.171
	33.021 17.233
	32.924 17.301
	32.832 17.376
	32.743 17.458
	32.660 17.546
	32.581 17.640
	32.506 17.741
	32.436 17.848
	32.370 17.962
	32.311 18.075
	32.262 18.179
	32.221 18.275
	32.189 18.363
	32.167 18.442
	32.154 18.512
	32.155 18.629
	32.188 18.722
	32.246 18.802
	32.331 18.869
	32.441 18.923
	32.506 18.946
	32.577 18.966
	32.654 18.984
	32.737 19.000
	32.827 19.014
	32.922 19.025
	33.024 19.035
	33.132 19.042
	33.239 19.046
	33.338 19.049
	33.428 19.050
	33.509 19.049
	33.582 19.046
	33.647 19.041
	33.751 19.026
	33.833 18.993
	33.902 18.935
	33.960 18.850
	34.005 18.740
	34.024 18.651
	34.032 18.586
	34.040 18.508
	34.047 18.415
	34.054 18.310
	34.060 18.190
	34.062 18.125
	34.065 18.057
	34.067 17.985
	34.069 17.910
	34.072 17.832
	34.074 17.750
	34.075 17.665
	34.077 17.576
	34.078 17.484
	34.080 17.388
	34.081 17.289
	34.082 17.187
	34.083 17.081
	34.083 16.972
	34.084 16.860
	34.084 16.744
	34.085 16.624
	34.085 16.502
	 /
\plot 34.085 16.502 34.085 14.510 /
\linethickness=1pt
\setplotsymbol ({\makebox(0,0)[l]{\tencirc\symbol{'160}}})
%
%
%
\plot	21.048  3.651 21.048  5.794
 	21.048  5.861
	21.048  5.926
	21.048  5.991
	21.048  6.054
	21.047  6.178
	21.047  6.299
	21.046  6.415
	21.045  6.527
	21.044  6.636
	21.043  6.740
	21.042  6.840
	21.040  6.937
	21.039  7.030
	21.037  7.118
	21.035  7.203
	21.033  7.284
	21.031  7.361
	21.028  7.433
	21.026  7.502
	21.023  7.567
	21.017  7.686
	21.011  7.788
	21.003  7.875
	20.987  8.001
	20.969  8.064
	20.870  8.174
	20.805  8.215
	20.731  8.248
	20.630  8.271
	20.565  8.280
	20.489  8.286
	20.402  8.290
	20.305  8.291
	20.198  8.291
	20.080  8.288
	19.960  8.285
	19.848  8.286
	19.743  8.289
	19.645  8.296
	19.555  8.305
	19.471  8.318
	19.395  8.333
	19.326  8.351
	19.208  8.403
	19.114  8.480
	19.045  8.582
	19.001  8.708
	18.992  8.781
	18.995  8.861
	19.011  8.948
	19.039  9.043
	19.080  9.145
	19.133  9.254
	19.200  9.370
	19.279  9.493
	19.365  9.616
	19.452  9.729
	19.541  9.833
	19.632  9.928
	19.724 10.013
	19.817 10.090
	19.912 10.157
	20.009 10.215
	20.103 10.265
	20.192 10.309
	20.274 10.345
	20.350 10.375
	20.421 10.399
	20.485 10.415
	20.596 10.429
	 /
\plot 20.596 10.429 20.794 10.429 /
\linethickness=1pt
\setplotsymbol ({\makebox(0,0)[l]{\tencirc\symbol{'160}}})
%
%
%
\plot	21.334 10.429 21.556 10.397
 	21.678 10.367
	21.747 10.341
	21.822 10.308
	21.902 10.268
	21.988 10.220
	22.079 10.166
	22.175 10.104
	22.272 10.035
	22.365  9.960
	22.453  9.879
	22.536  9.791
	22.616  9.697
	22.690  9.596
	22.760  9.488
	22.826  9.375
	22.885  9.262
	22.935  9.158
	22.975  9.062
	23.007  8.974
	23.029  8.895
	23.042  8.824
	23.041  8.708
	23.008  8.615
	22.950  8.535
	22.866  8.468
	22.755  8.414
	22.690  8.391
	22.620  8.371
	22.542  8.352
	22.459  8.336
	22.370  8.323
	22.274  8.311
	22.172  8.302
	22.064  8.295
	21.957  8.290
	21.859  8.287
	21.769  8.287
	21.687  8.288
	21.614  8.291
	21.549  8.295
	21.445  8.311
	21.364  8.344
	21.294  8.402
	21.237  8.487
	21.191  8.597
	21.172  8.686
	21.164  8.751
	21.156  8.829
	21.149  8.921
	21.143  9.027
	21.137  9.147
	21.134  9.211
	21.131  9.280
	21.129  9.351
	21.127  9.426
	21.125  9.505
	21.123  9.587
	21.121  9.672
	21.119  9.761
	21.118  9.853
	21.117  9.948
	21.115 10.047
	21.114 10.150
	21.114 10.255
	21.113 10.365
	21.112 10.477
	21.112 10.593
	21.112 10.712
	21.112 10.835
	 /
\plot 21.112 10.835 21.112 12.827 /
\linethickness=1pt
\setplotsymbol ({\makebox(0,0)[l]{\tencirc\symbol{'160}}})
%
%
%
\plot	27.874  1.890 27.874 -0.253
 	27.874 -0.319
	27.874 -0.385
	27.875 -0.449
	27.875 -0.513
	27.875 -0.637
	27.876 -0.757
	27.876 -0.873
	27.877 -0.986
	27.878 -1.094
	27.879 -1.198
	27.881 -1.299
	27.882 -1.396
	27.884 -1.488
	27.886 -1.577
	27.887 -1.661
	27.890 -1.742
	27.892 -1.819
	27.894 -1.892
	27.897 -1.961
	27.899 -2.026
	27.905 -2.144
	27.912 -2.247
	27.919 -2.334
	27.935 -2.460
	27.954 -2.523
	28.053 -2.632
	28.117 -2.674
	28.192 -2.706
	28.292 -2.730
	28.358 -2.738
	28.434 -2.744
	28.521 -2.748
	28.618 -2.750
	28.725 -2.749
	28.843 -2.746
	28.962 -2.744
	29.074 -2.744
	29.179 -2.748
	29.277 -2.754
	29.368 -2.764
	29.451 -2.776
	29.527 -2.792
	29.596 -2.810
	29.715 -2.862
	29.808 -2.939
	29.877 -3.040
	29.921 -3.167
	29.931 -3.239
	29.928 -3.320
	29.912 -3.407
	29.884 -3.501
	29.843 -3.603
	29.789 -3.712
	29.723 -3.828
	29.644 -3.952
	29.558 -4.074
	29.470 -4.187
	29.381 -4.291
	29.291 -4.386
	29.199 -4.472
	29.105 -4.548
	29.010 -4.616
	28.914 -4.674
	28.819 -4.724
	28.731 -4.767
	28.649 -4.804
	28.572 -4.834
	28.502 -4.857
	28.437 -4.874
	28.326 -4.887
	 /
\plot 28.326 -4.887 28.128 -4.887 /
\linethickness=1pt
\setplotsymbol ({\makebox(0,0)[l]{\tencirc\symbol{'160}}})
%
%
%
\plot	21.283  1.858 21.283 -0.285
 	21.283 -0.351
	21.283 -0.417
	21.283 -0.481
	21.283 -0.545
	21.284 -0.669
	21.284 -0.789
	21.285 -0.905
	21.286 -1.017
	21.287 -1.126
	21.288 -1.230
	21.289 -1.331
	21.291 -1.427
	21.292 -1.520
	21.294 -1.609
	21.296 -1.693
	21.298 -1.774
	21.301 -1.851
	21.303 -1.924
	21.305 -1.993
	21.308 -2.058
	21.314 -2.176
	21.321 -2.279
	21.328 -2.365
	21.344 -2.492
	21.362 -2.555
	21.462 -2.664
	21.526 -2.705
	21.601 -2.738
	21.701 -2.761
	21.767 -2.770
	21.843 -2.776
	21.929 -2.780
	22.026 -2.782
	22.134 -2.781
	22.251 -2.778
	22.371 -2.776
	22.483 -2.776
	22.588 -2.780
	22.686 -2.786
	22.776 -2.795
	22.860 -2.808
	22.936 -2.823
	23.005 -2.842
	23.123 -2.894
	23.217 -2.970
	23.286 -3.072
	23.330 -3.198
	23.339 -3.271
	23.336 -3.351
	23.321 -3.439
	23.292 -3.533
	23.251 -3.635
	23.198 -3.744
	23.131 -3.860
	23.053 -3.984
	22.967 -4.106
	22.879 -4.219
	22.790 -4.323
	22.700 -4.418
	22.608 -4.504
	22.514 -4.580
	22.419 -4.647
	22.322 -4.705
	22.228 -4.755
	22.140 -4.799
	22.057 -4.836
	21.981 -4.866
	21.910 -4.889
	21.846 -4.906
	21.735 -4.919
	 /
\plot 21.735 -4.919 21.537 -4.919 /
\linethickness=1pt
\setplotsymbol ({\makebox(0,0)[l]{\tencirc\symbol{'160}}})
%
%
%
\plot	33.803  1.985 33.803 -0.158
 	33.803 -0.224
	33.803 -0.290
	33.803 -0.354
	33.803 -0.418
	33.804 -0.542
	33.804 -0.662
	33.805 -0.778
	33.806 -0.890
	33.807 -0.999
	33.808 -1.103
	33.809 -1.204
	33.811 -1.300
	33.813 -1.393
	33.814 -1.482
	33.816 -1.566
	33.818 -1.647
	33.821 -1.724
	33.823 -1.797
	33.826 -1.866
	33.828 -1.931
	33.834 -2.049
	33.841 -2.152
	33.848 -2.238
	33.864 -2.365
	33.883 -2.428
	33.982 -2.537
	34.046 -2.578
	34.121 -2.611
	34.221 -2.634
	34.287 -2.643
	34.363 -2.649
	34.449 -2.653
	34.546 -2.655
	34.654 -2.654
	34.772 -2.651
	34.891 -2.649
	35.003 -2.649
	35.108 -2.653
	35.206 -2.659
	35.297 -2.668
	35.380 -2.681
	35.456 -2.696
	35.525 -2.715
	35.643 -2.767
	35.737 -2.843
	35.806 -2.945
	35.850 -3.071
	35.860 -3.144
	35.856 -3.224
	35.841 -3.312
	35.812 -3.406
	35.771 -3.508
	35.718 -3.617
	35.652 -3.733
	35.573 -3.857
	35.487 -3.979
	35.399 -4.092
	35.310 -4.196
	35.220 -4.291
	35.128 -4.377
	35.034 -4.453
	34.939 -4.520
	34.842 -4.578
	34.748 -4.628
	34.660 -4.672
	34.577 -4.709
	34.501 -4.739
	34.430 -4.762
	34.366 -4.779
	34.255 -4.792
	 /
\plot 34.255 -4.792 34.057 -4.792 /
\linethickness=1pt
\setplotsymbol ({\makebox(0,0)[l]{\tencirc\symbol{'160}}})
%
%
%
\plot	15.359 -0.540 15.581 -0.572
 	15.703 -0.602
	15.772 -0.628
	15.847 -0.661
	15.927 -0.701
	16.012 -0.748
	16.103 -0.803
	16.200 -0.865
	16.297 -0.933
	16.389 -1.008
	16.477 -1.090
	16.561 -1.178
	16.640 -1.272
	16.715 -1.373
	16.785 -1.480
	16.851 -1.594
	16.910 -1.707
	16.959 -1.811
	17.000 -1.907
	17.031 -1.994
	17.054 -2.074
	17.067 -2.144
	17.066 -2.261
	17.033 -2.354
	16.975 -2.434
	16.890 -2.501
	16.780 -2.555
	16.715 -2.578
	16.644 -2.598
	16.567 -2.616
	16.484 -2.632
	16.394 -2.646
	16.299 -2.657
	16.197 -2.666
	16.089 -2.673
	15.982 -2.678
	15.883 -2.681
	15.793 -2.682
	15.712 -2.681
	15.639 -2.678
	15.574 -2.673
	15.470 -2.657
	15.388 -2.625
	15.319 -2.566
	15.261 -2.482
	15.216 -2.372
	15.197 -2.283
	15.189 -2.218
	15.181 -2.139
	15.174 -2.047
	15.167 -1.941
	15.161 -1.822
	15.159 -1.757
	15.156 -1.689
	15.154 -1.617
	15.151 -1.542
	15.149 -1.464
	15.147 -1.382
	15.146 -1.296
	15.144 -1.208
	15.143 -1.116
	15.141 -1.020
	15.140 -0.921
	15.139 -0.819
	15.138 -0.713
	15.138 -0.604
	15.137 -0.491
	15.137 -0.376
	15.136 -0.256
	15.136 -0.133
	 /
\plot 15.136 -0.133 15.136  1.858 /
\linethickness=1pt
\setplotsymbol ({\makebox(0,0)[l]{\tencirc\symbol{'160}}})
%
%
%
\plot	 8.528 -0.571  8.750 -0.603
 	 8.872 -0.633
	 8.942 -0.659
	 9.016 -0.692
	 9.096 -0.733
	 9.182 -0.780
	 9.273 -0.835
	 9.369 -0.896
	 9.466 -0.965
	 9.559 -1.040
	 9.647 -1.121
	 9.731 -1.209
	 9.810 -1.304
	 9.884 -1.405
	 9.955 -1.512
	10.020 -1.626
	10.079 -1.738
	10.129 -1.843
	10.170 -1.939
	10.201 -2.026
	10.223 -2.105
	10.236 -2.176
	10.235 -2.292
	10.203 -2.385
	10.144 -2.465
	10.060 -2.532
	 9.949 -2.587
	 9.885 -2.609
	 9.814 -2.630
	 9.737 -2.648
	 9.653 -2.664
	 9.564 -2.678
	 9.468 -2.689
	 9.366 -2.698
	 9.258 -2.705
	 9.151 -2.710
	 9.053 -2.713
	 8.963 -2.714
	 8.881 -2.713
	 8.808 -2.710
	 8.743 -2.705
	 8.639 -2.689
	 8.558 -2.657
	 8.488 -2.598
	 8.431 -2.514
	 8.385 -2.403
	 8.367 -2.315
	 8.358 -2.250
	 8.350 -2.171
	 8.343 -2.079
	 8.337 -1.973
	 8.331 -1.854
	 8.328 -1.789
	 8.326 -1.721
	 8.323 -1.649
	 8.321 -1.574
	 8.319 -1.495
	 8.317 -1.414
	 8.315 -1.328
	 8.314 -1.240
	 8.312 -1.147
	 8.311 -1.052
	 8.310 -0.953
	 8.309 -0.851
	 8.308 -0.745
	 8.307 -0.636
	 8.306 -0.523
	 8.306 -0.407
	 8.306 -0.288
	 8.306 -0.165
	 /
\plot  8.306 -0.165  8.306  1.827 /
\linethickness=1pt
\setplotsymbol ({\makebox(0,0)[l]{\tencirc\symbol{'160}}})
%
%
%
\plot	27.589 -4.887 27.366 -4.856
 	27.244 -4.825
	27.175 -4.800
	27.100 -4.766
	27.020 -4.726
	26.935 -4.679
	26.844 -4.624
	26.747 -4.562
	26.650 -4.494
	26.558 -4.419
	26.470 -4.337
	26.386 -4.249
	26.307 -4.155
	26.232 -4.054
	26.162 -3.947
	26.096 -3.833
	26.037 -3.721
	25.988 -3.616
	25.947 -3.520
	25.916 -3.433
	25.893 -3.354
	25.880 -3.283
	25.882 -3.167
	25.914 -3.074
	25.973 -2.993
	26.057 -2.926
	26.167 -2.872
	26.232 -2.850
	26.303 -2.829
	26.380 -2.811
	26.463 -2.795
	26.553 -2.781
	26.649 -2.770
	26.750 -2.761
	26.858 -2.754
	26.965 -2.749
	27.064 -2.746
	27.154 -2.745
	27.235 -2.746
	27.309 -2.749
	27.373 -2.754
	27.478 -2.770
	27.559 -2.802
	27.628 -2.861
	27.686 -2.945
	27.732 -3.055
	27.750 -3.144
	27.758 -3.209
	27.766 -3.288
	27.773 -3.380
	27.780 -3.486
	27.786 -3.605
	27.788 -3.670
	27.791 -3.738
	27.793 -3.810
	27.796 -3.885
	27.798 -3.963
	27.800 -4.045
	27.802 -4.131
	27.803 -4.219
	27.805 -4.311
	27.806 -4.407
	27.807 -4.506
	27.808 -4.608
	27.809 -4.714
	27.810 -4.823
	27.810 -4.936
	27.811 -5.052
	27.811 -5.171
	27.811 -5.294
	 /
\plot 27.811 -5.294 27.811 -7.286 /
\linethickness=1pt
\setplotsymbol ({\makebox(0,0)[l]{\tencirc\symbol{'160}}})
%
%
%
\plot	20.997 -4.919 20.775 -4.887
 	20.653 -4.857
	20.584 -4.831
	20.509 -4.798
	20.429 -4.758
	20.343 -4.711
	20.252 -4.656
	20.156 -4.594
	20.059 -4.526
	19.966 -4.451
	19.878 -4.369
	19.795 -4.281
	19.716 -4.187
	19.641 -4.086
	19.571 -3.979
	19.505 -3.865
	19.446 -3.752
	19.396 -3.648
	19.356 -3.552
	19.324 -3.464
	19.302 -3.385
	19.289 -3.315
	19.290 -3.198
	19.323 -3.105
	19.381 -3.025
	19.466 -2.958
	19.576 -2.904
	19.641 -2.881
	19.712 -2.861
	19.789 -2.843
	19.872 -2.827
	19.962 -2.813
	20.057 -2.802
	20.159 -2.792
	20.267 -2.786
	20.374 -2.781
	20.472 -2.778
	20.563 -2.777
	20.644 -2.778
	20.717 -2.781
	20.782 -2.786
	20.886 -2.801
	20.968 -2.834
	21.037 -2.892
	21.095 -2.977
	21.140 -3.087
	21.159 -3.176
	21.167 -3.241
	21.175 -3.320
	21.182 -3.412
	21.189 -3.518
	21.194 -3.637
	21.197 -3.702
	21.200 -3.770
	21.202 -3.842
	21.204 -3.917
	21.206 -3.995
	21.208 -4.077
	21.210 -4.162
	21.212 -4.251
	21.213 -4.343
	21.215 -4.439
	21.216 -4.538
	21.217 -4.640
	21.218 -4.746
	21.218 -4.855
	21.219 -4.967
	21.219 -5.083
	21.220 -5.203
	21.220 -5.326
	 /
\plot 21.220 -5.326 21.220 -7.317 /
\linethickness=1pt
\setplotsymbol ({\makebox(0,0)[l]{\tencirc\symbol{'160}}})
%
%
%
\plot	33.517 -4.792 33.295 -4.760
 	33.173 -4.730
	33.104 -4.704
	33.029 -4.671
	32.949 -4.631
	32.864 -4.584
	32.773 -4.529
	32.676 -4.467
	32.579 -4.399
	32.487 -4.324
	32.398 -4.242
	32.315 -4.154
	32.236 -4.060
	32.161 -3.959
	32.091 -3.852
	32.025 -3.738
	31.966 -3.625
	31.916 -3.521
	31.876 -3.425
	31.844 -3.337
	31.822 -3.258
	31.809 -3.188
	31.810 -3.071
	31.843 -2.978
	31.901 -2.898
	31.986 -2.831
	32.096 -2.777
	32.161 -2.754
	32.232 -2.734
	32.309 -2.716
	32.392 -2.700
	32.482 -2.686
	32.577 -2.675
	32.679 -2.665
	32.787 -2.659
	32.894 -2.654
	32.993 -2.651
	33.083 -2.650
	33.164 -2.651
	33.237 -2.654
	33.302 -2.659
	33.406 -2.674
	33.488 -2.707
	33.557 -2.765
	33.615 -2.850
	33.660 -2.960
	33.679 -3.049
	33.687 -3.114
	33.695 -3.193
	33.702 -3.285
	33.709 -3.391
	33.715 -3.510
	33.717 -3.575
	33.720 -3.643
	33.722 -3.715
	33.724 -3.790
	33.727 -3.868
	33.729 -3.950
	33.730 -4.035
	33.732 -4.124
	33.733 -4.216
	33.735 -4.312
	33.736 -4.411
	33.737 -4.513
	33.738 -4.619
	33.738 -4.728
	33.739 -4.840
	33.739 -4.956
	33.740 -5.076
	33.740 -5.199
	 /
\plot 33.740 -5.199 33.740 -7.190 /
\linethickness=1pt
\setplotsymbol ({\makebox(0,0)[l]{\tencirc\symbol{'160}}})
%
%
%
\plot	15.073 -7.317 15.073 -5.174
 	15.073 -5.108
	15.073 -5.042
	15.073 -4.978
	15.072 -4.914
	15.072 -4.790
	15.072 -4.670
	15.071 -4.554
	15.070 -4.441
	15.069 -4.333
	15.068 -4.229
	15.067 -4.128
	15.065 -4.032
	15.063 -3.939
	15.062 -3.850
	15.060 -3.766
	15.058 -3.685
	15.055 -3.608
	15.053 -3.535
	15.050 -3.466
	15.048 -3.401
	15.042 -3.283
	15.035 -3.180
	15.028 -3.094
	15.012 -2.967
	14.993 -2.904
	14.894 -2.795
	14.830 -2.753
	14.755 -2.721
	14.655 -2.697
	14.589 -2.689
	14.513 -2.683
	14.427 -2.679
	14.330 -2.677
	14.222 -2.678
	14.104 -2.681
	13.985 -2.683
	13.873 -2.683
	13.768 -2.679
	13.670 -2.673
	13.579 -2.663
	13.496 -2.651
	13.420 -2.636
	13.351 -2.617
	13.233 -2.565
	13.139 -2.488
	13.070 -2.387
	13.026 -2.261
	13.016 -2.188
	13.019 -2.108
	13.035 -2.020
	13.064 -1.926
	13.105 -1.824
	13.158 -1.715
	13.224 -1.599
	13.303 -1.475
	13.389 -1.353
	13.477 -1.240
	13.566 -1.136
	13.656 -1.041
	13.748 -0.955
	13.842 -0.879
	13.937 -0.812
	14.034 -0.754
	14.128 -0.703
	14.216 -0.660
	14.299 -0.623
	14.375 -0.593
	14.446 -0.570
	14.510 -0.553
	14.621 -0.540
	 /
\plot 14.621 -0.540 14.819 -0.540 /
\linethickness=1pt
\setplotsymbol ({\makebox(0,0)[l]{\tencirc\symbol{'160}}})
%
%
%
\plot	 8.242 -7.349  8.242 -5.206
 	 8.242 -5.139
	 8.242 -5.074
	 8.242 -5.009
	 8.242 -4.946
	 8.242 -4.822
	 8.241 -4.702
	 8.240 -4.585
	 8.240 -4.473
	 8.239 -4.365
	 8.237 -4.260
	 8.236 -4.160
	 8.235 -4.063
	 8.233 -3.971
	 8.231 -3.882
	 8.229 -3.797
	 8.227 -3.717
	 8.225 -3.640
	 8.222 -3.567
	 8.220 -3.498
	 8.217 -3.433
	 8.211 -3.315
	 8.205 -3.212
	 8.198 -3.125
	 8.182 -2.999
	 8.163 -2.936
	 8.064 -2.827
	 7.999 -2.785
	 7.925 -2.753
	 7.825 -2.729
	 7.759 -2.721
	 7.683 -2.715
	 7.596 -2.711
	 7.499 -2.709
	 7.392 -2.710
	 7.274 -2.713
	 7.154 -2.715
	 7.042 -2.714
	 6.937 -2.711
	 6.839 -2.705
	 6.749 -2.695
	 6.666 -2.683
	 6.589 -2.667
	 6.520 -2.649
	 6.402 -2.597
	 6.308 -2.520
	 6.240 -2.419
	 6.195 -2.292
	 6.186 -2.219
	 6.189 -2.139
	 6.205 -2.052
	 6.233 -1.957
	 6.274 -1.856
	 6.328 -1.747
	 6.394 -1.630
	 6.473 -1.507
	 6.559 -1.385
	 6.646 -1.271
	 6.735 -1.167
	 6.826 -1.073
	 6.918 -0.987
	 7.011 -0.911
	 7.106 -0.843
	 7.203 -0.785
	 7.297 -0.735
	 7.386 -0.692
	 7.468 -0.655
	 7.545 -0.625
	 7.615 -0.602
	 7.679 -0.585
	 7.790 -0.571
	 /
\plot  7.790 -0.571  7.988 -0.571 /
\linethickness=1pt
\setplotsymbol ({\makebox(0,0)[l]{\tencirc\symbol{'160}}})
%
%
%
\plot	27.970 12.939 27.970 10.796
 	27.970 10.730
	27.970 10.664
	27.970 10.600
	27.970 10.536
	27.970 10.412
	27.971 10.292
	27.972 10.176
	27.972 10.063
	27.973  9.955
	27.975  9.851
	27.976  9.750
	27.977  9.653
	27.979  9.561
	27.981  9.472
	27.983  9.388
	27.985  9.307
	27.987  9.230
	27.989  9.157
	27.992  9.088
	27.995  9.023
	28.001  8.905
	28.007  8.802
	28.014  8.715
	28.030  8.589
	28.049  8.526
	28.148  8.417
	28.213  8.375
	28.287  8.343
	28.387  8.319
	28.453  8.311
	28.529  8.305
	28.616  8.301
	28.713  8.299
	28.820  8.300
	28.938  8.303
	29.057  8.305
	29.170  8.305
	29.275  8.301
	29.372  8.295
	29.463  8.285
	29.546  8.273
	29.623  8.257
	29.692  8.239
	29.810  8.187
	29.903  8.110
	29.972  8.009
	30.016  7.882
	30.026  7.810
	30.023  7.729
	30.007  7.642
	29.979  7.548
	29.938  7.446
	29.884  7.337
	29.818  7.221
	29.739  7.097
	29.653  6.975
	29.566  6.862
	29.477  6.758
	29.386  6.663
	29.294  6.577
	29.201  6.501
	29.106  6.433
	29.009  6.375
	28.915  6.325
	28.826  6.282
	28.744  6.245
	28.667  6.215
	28.597  6.192
	28.532  6.175
	28.422  6.162
	 /
\plot 28.422  6.162 28.224  6.162 /
\linethickness=1pt
\setplotsymbol ({\makebox(0,0)[l]{\tencirc\symbol{'160}}})
%
%
%
\plot	27.684  6.162 27.462  6.193
 	27.340  6.224
	27.270  6.249
	27.196  6.283
	27.116  6.323
	27.030  6.370
	26.939  6.425
	26.843  6.487
	26.746  6.555
	26.653  6.630
	26.565  6.712
	26.481  6.800
	26.402  6.894
	26.328  6.995
	26.257  7.102
	26.192  7.216
	26.133  7.328
	26.083  7.433
	26.042  7.529
	26.011  7.616
	25.989  7.695
	25.976  7.766
	25.977  7.882
	26.009  7.975
	26.068  8.056
	26.152  8.123
	26.263  8.177
	26.327  8.199
	26.398  8.220
	26.475  8.238
	26.559  8.254
	26.648  8.268
	26.744  8.279
	26.846  8.288
	26.954  8.295
	27.061  8.300
	27.159  8.303
	27.249  8.304
	27.331  8.303
	27.404  8.300
	27.469  8.295
	27.573  8.279
	27.654  8.247
	27.724  8.188
	27.781  8.104
	27.827  7.994
	27.845  7.905
	27.854  7.840
	27.861  7.761
	27.869  7.669
	27.875  7.563
	27.881  7.444
	27.884  7.379
	27.886  7.311
	27.889  7.239
	27.891  7.164
	27.893  7.086
	27.895  7.004
	27.897  6.918
	27.898  6.830
	27.900  6.738
	27.901  6.642
	27.902  6.543
	27.903  6.441
	27.904  6.335
	27.905  6.226
	27.905  6.113
	27.906  5.997
	27.906  5.878
	27.906  5.755
	 /
\plot 27.906  5.755 27.906  3.763 /
\linethickness=1pt
\setplotsymbol ({\makebox(0,0)[l]{\tencirc\symbol{'160}}})
%
%
%
\plot	33.898 13.034 33.898 10.891
 	33.898 10.825
	33.898 10.759
	33.899 10.695
	33.899 10.631
	33.899 10.507
	33.900 10.387
	33.900 10.271
	33.901 10.159
	33.902 10.050
	33.903  9.946
	33.905  9.845
	33.906  9.749
	33.908  9.656
	33.910  9.567
	33.912  9.483
	33.914  9.402
	33.916  9.325
	33.918  9.252
	33.921  9.183
	33.924  9.118
	33.929  9.000
	33.936  8.897
	33.943  8.811
	33.959  8.684
	33.978  8.621
	34.077  8.512
	34.142  8.471
	34.216  8.438
	34.316  8.415
	34.382  8.406
	34.458  8.400
	34.545  8.396
	34.642  8.394
	34.749  8.395
	34.867  8.398
	34.986  8.400
	35.098  8.400
	35.203  8.396
	35.301  8.390
	35.392  8.381
	35.475  8.368
	35.551  8.353
	35.620  8.334
	35.739  8.282
	35.832  8.206
	35.901  8.104
	35.945  7.978
	35.955  7.905
	35.952  7.825
	35.936  7.737
	35.908  7.643
	35.867  7.541
	35.813  7.432
	35.747  7.316
	35.668  7.192
	35.582  7.070
	35.495  6.957
	35.406  6.853
	35.315  6.758
	35.223  6.672
	35.129  6.596
	35.034  6.529
	34.938  6.471
	34.843  6.421
	34.755  6.377
	34.673  6.340
	34.596  6.310
	34.526  6.287
	34.461  6.270
	34.350  6.257
	 /
\plot 34.350  6.257 34.152  6.257 /
\linethickness=1pt
\setplotsymbol ({\makebox(0,0)[l]{\tencirc\symbol{'160}}})
%
%
%
\plot	33.613  6.257 33.390  6.289
 	33.268  6.319
	33.199  6.345
	33.125  6.378
	33.044  6.418
	32.959  6.465
	32.868  6.520
	32.771  6.582
	32.674  6.650
	32.582  6.725
	32.494  6.807
	32.410  6.895
	32.331  6.989
	32.256  7.090
	32.186  7.197
	32.120  7.311
	32.061  7.424
	32.012  7.528
	31.971  7.624
	31.940  7.712
	31.917  7.791
	31.904  7.861
	31.906  7.978
	31.938  8.071
	31.997  8.151
	32.081  8.218
	32.191  8.272
	32.256  8.295
	32.327  8.315
	32.404  8.333
	32.487  8.349
	32.577  8.363
	32.673  8.374
	32.774  8.384
	32.882  8.390
	32.989  8.395
	33.088  8.398
	33.178  8.399
	33.259  8.398
	33.333  8.395
	33.397  8.390
	33.502  8.375
	33.583  8.342
	33.652  8.284
	33.710  8.199
	33.756  8.089
	33.774  8.000
	33.783  7.935
	33.790  7.856
	33.797  7.764
	33.804  7.658
	33.810  7.539
	33.813  7.474
	33.815  7.406
	33.817  7.334
	33.820  7.259
	33.822  7.181
	33.824  7.099
	33.826  7.014
	33.827  6.925
	33.829  6.833
	33.830  6.737
	33.831  6.638
	33.832  6.536
	33.833  6.430
	33.834  6.321
	33.834  6.209
	33.835  6.093
	33.835  5.973
	33.835  5.850
	 /
\plot 33.835  5.850 33.835  3.859 /
\linethickness=0pt
\putrectangle corners at  6.132 23.717 and 36.258 -7.381
\endpicture
\endfig

This figure should be understood as follows:  Time passes down the
page.  At times before 0 the picture is static and is the same
as figure \mtwist.  Just before time 1 (the top row of the figure)
a small 2--sphere has appeared (a 0--handle).  At time 1 this sphere
has grown a little and now appears as the three little circles
in the middle three slices (think of the sphere as a circle which grows
from a point and then shrinks back down).  At time 2 (the middle
row of the figure) the sphere has grown to the point where it
encloses the finger.  At this point a 1--handle bridges between
the finger and the sphere and this has the effect of flipping
the finger over from a left finger to a right finger as shown in
the bottom row (time 3).

It should now be clear how to continue this sequence of constructions
to construct embeddings of $n$--space in $(n+2)$--space for all $n$.

\sh{Higher codimensions}

We shall now describe the end result of the straightening process
for the simplest situation in higher codimensions, namely an
inclined plane with the analogue of the twist field of figure
\twist.   More precisely we start with an embedding of
$\re^c$ in $\re^{2c+1}$ with a perpendicular field which points
up outside a disc, down at the centre of the disc and at other
points of the disc has direction given by identifying the
disc (rel boundary) with the normal sphere to $\re^c$ in $\re^{2c+1}$.
We shall describe the result of applying the compression theorem. 
The move which achieves this result is a similar twist to the 
codimension 2 twist of figures \vectp\ to \twist\ above.
We shall describe an immersion of
$\re^c$ in $\re^{2c}$ with a single double point covered by
an embedding in $\re^{2c+1}$.

The case  $c=2$ can be described as follows:

Take a straight line
in $\re^3$.  Put a twist in it.  Pass the line through itself to
get the opposite twist.  Pull straight again.  This moving picture
of immersions of $\re^1$ in $\re^3$ defines an
immersion of $\re^2$ in $\re^4$.  Lying above it is a moving picture
of $\re^1$ in $\re^4$, where at the critical stage we have the $\re^1$ 
embedded in a 3-dimensional subspace (and have a twist like figure
\twist\ there).

We construct the general picture inductively on $c$.
Suppose it is constructed for $c$.  So we have an
embedding of $\re^c$ in $\re^{2c+1}$ lying over an immersion in
$\re^{2c}$ with a single double point.  Now consider this immersion
as being in $\re^{2c+1}$ then we can undo
the immersion in two different ways: lift the double point off
upwards and then pull flat and ditto downwards.  Performing the reverse
of one of these undoings followed by the other describes
a moving picture of $\re^c$ in $\re^{2c+1}$, 
ie, an immersion of $\re^{c+1}$ in $\re^{2c+2}$
with a single double point.  This is the projection of the
next standard picture.  To get the lift, just use the extra dimension
to lift off the double point.

The move which produces the immersions just described is very
similar to the move described by Koschorke and Sanderson
in [\KSa; Theorem 5.2], see in particular the picture on [\KSa; page
216].

\sh{Removal of a Whitney umbrella}

We now return to the initial example described in figures \vectp\
to \twist, but from a different point of view.  The compression
isotopy is a sequence of embeddings of a line in $\re^3$ with normal
vector fields.  The whole sequence defines an embedding of a 
plane in $\re^4$ with a normal
vector field.  The projection on $\re^3$ of this plane has a
singularity---in fact a Whitney umbrella.  To see this,
we have reproduced this sequence in figure \figkey\Whit\ below.  We
have changed the isotopy to an equivalent one which shows the
singularity clearly.  

\fig{\Whit: Surface above a Whitney umbrella}
\beginpicture
\setcoordinatesystem units <0.50000cm,0.50000cm>
\unitlength=1.00000cm
\linethickness=1pt
\setplotsymbol ({\makebox(0,0)[l]{\tencirc\symbol{'160}}})
\setshadesymbol ({\thinlinefont .})
\setlinear
%
%
\linethickness= 0.500pt
\setplotsymbol ({\thinlinefont .})
\putrule from 18.540 24.481 to 18.540 25.370
%
%
\plot 18.603 25.116 18.540 25.370 18.476 25.116 /
%
%
%
\linethickness= 0.500pt
\setplotsymbol ({\thinlinefont .})
\plot 21.969 24.037 21.960 25.021 /
%
%
\plot 22.026 24.768 21.960 25.021 21.899 24.767 /
%
%
%
\linethickness= 0.500pt
\setplotsymbol ({\thinlinefont .})
\putrule from 20.292 24.354 to 20.292 25.229
%
%
\plot 20.356 24.975 20.292 25.229 20.229 24.975 /
%
%
%
\linethickness= 0.500pt
\setplotsymbol ({\thinlinefont .})
\putrule from 23.404 23.783 to 23.404 24.752
%
%
\plot 23.467 24.498 23.404 24.752 23.340 24.498 /
%
%
%
\linethickness= 0.500pt
\setplotsymbol ({\thinlinefont .})
\putrule from 23.848 22.593 to 23.848 23.370
%
%
\plot 23.912 23.116 23.848 23.370 23.785 23.116 /
%
%
%
\linethickness= 0.500pt
\setplotsymbol ({\thinlinefont .})
\plot 23.975 20.720 23.992 21.814 /
%
%
\plot 24.052 21.559 23.992 21.814 23.925 21.561 /
%
%
%
\linethickness= 0.500pt
\setplotsymbol ({\thinlinefont .})
\putrule from 26.198 20.085 to 26.198 21.084
%
%
\plot 26.261 20.830 26.198 21.084 26.134 20.830 /
%
%
%
\linethickness= 0.500pt
\setplotsymbol ({\thinlinefont .})
\putrule from 28.452 19.560 to 28.452 20.544
%
%
\plot 28.516 20.290 28.452 20.544 28.389 20.290 /
%
%
%
\linethickness= 0.500pt
\setplotsymbol ({\thinlinefont .})
\putrule from 18.656 17.693 to 18.656 18.709
%
%
\plot 18.720 18.455 18.656 18.709 18.593 18.455 /
%
%
%
\linethickness= 0.500pt
\setplotsymbol ({\thinlinefont .})
\plot 23.719 15.265 24.816 15.280 /
%
%
\plot 24.563 15.213 24.816 15.280 24.561 15.340 /
%
%
%
\linethickness= 0.500pt
\setplotsymbol ({\thinlinefont .})
\putrule from 25.148 13.106 to 25.148 13.995
%
%
\plot 25.212 13.741 25.148 13.995 25.085 13.741 /
%
%
%
\linethickness= 0.500pt
\setplotsymbol ({\thinlinefont .})
\putrule from 27.673 12.931 to 27.673 13.805
%
%
\plot 27.737 13.551 27.673 13.805 27.610 13.551 /
%
%
%
\linethickness= 0.500pt
\setplotsymbol ({\thinlinefont .})
\putrule from 20.576 17.560 to 20.576 18.481
%
%
\plot 20.640 18.227 20.576 18.481 20.513 18.227 /
%
%
%
\linethickness= 0.500pt
\setplotsymbol ({\thinlinefont .})
\plot 22.608 17.115 23.084 17.973 /
%
%
\plot 23.017 17.720 23.084 17.973 22.906 17.781 /
%
%
%
\linethickness= 0.500pt
\setplotsymbol ({\thinlinefont .})
\plot 23.349 16.305 24.731 16.622 /
%
%
\plot 24.498 16.503 24.731 16.622 24.469 16.627 /
%
%
%
\linethickness=1pt
\setplotsymbol ({\makebox(0,0)[l]{\tencirc\symbol{'160}}})
\plot 22.515  8.291 29.180  5.908 /
%
%
\linethickness= 0.500pt
\setplotsymbol ({\thinlinefont .})
\plot 18.768  9.701 19.135 10.604 /
%
%
\plot 19.098 10.345 19.135 10.604 18.980 10.393 /
%
%
%
\linethickness= 0.500pt
\setplotsymbol ({\thinlinefont .})
\plot 20.388  9.049 20.786 10.001 /
%
%
\plot 20.746  9.742 20.786 10.001 20.629  9.791 /
%
%
%
\linethickness= 0.500pt
\setplotsymbol ({\thinlinefont .})
\plot 23.213  8.001 22.928  7.256 /
%
%
\plot 22.959  7.516 22.928  7.256 23.078  7.470 /
%
%
%
\linethickness= 0.500pt
\setplotsymbol ({\thinlinefont .})
\plot 24.642  7.493 24.382  6.716 /
%
%
\plot 24.402  6.977 24.382  6.716 24.523  6.937 /
%
%
%
\linethickness= 0.500pt
\setplotsymbol ({\thinlinefont .})
\plot 26.945  6.716 27.392  7.694 /
%
%
\plot 27.344  7.437 27.392  7.694 27.229  7.489 /
%
%
%
\linethickness= 0.500pt
\setplotsymbol ({\thinlinefont .})
\plot 28.723  6.113 29.199  7.192 /
%
%
\plot 29.155  6.934 29.199  7.192 29.039  6.986 /
%
%
%
\linethickness=1pt
\setplotsymbol ({\makebox(0,0)[l]{\tencirc\symbol{'160}}})
\plot 18.385  9.749 22.020  8.469 /
%
%
\linethickness=1pt
\setplotsymbol ({\makebox(0,0)[l]{\tencirc\symbol{'160}}})
\plot  5.954 15.636 16.161 15.623 /
%
%
\linethickness=1pt
\setplotsymbol ({\makebox(0,0)[l]{\tencirc\symbol{'160}}})
\plot  9.373 15.900  9.749 15.636 /
\plot  9.749 15.636  9.356 15.361 /
%
%
\linethickness=1pt
\setplotsymbol ({\makebox(0,0)[l]{\tencirc\symbol{'160}}})
\plot  6.001  8.748 16.239  8.750 /
%
%
\linethickness=1pt
\setplotsymbol ({\makebox(0,0)[l]{\tencirc\symbol{'160}}})
\plot  9.881  9.091 10.293  8.742 /
\plot 10.293  8.742  9.864  8.471 /
\linethickness=1pt
\setplotsymbol ({\makebox(0,0)[l]{\tencirc\symbol{'160}}})
%
%
%
\plot	18.127 24.513 19.842 24.354
 	19.948 24.344
	20.054 24.334
	20.158 24.324
	20.261 24.313
	20.363 24.302
	20.464 24.291
	20.564 24.280
	20.663 24.268
	20.761 24.256
	20.857 24.244
	20.953 24.232
	21.047 24.220
	21.140 24.207
	21.232 24.194
	21.323 24.181
	21.413 24.168
	21.502 24.154
	21.590 24.141
	21.676 24.127
	21.762 24.113
	21.846 24.098
	21.929 24.084
	22.011 24.069
	22.092 24.054
	22.172 24.038
	22.250 24.023
	22.328 24.007
	22.404 23.991
	22.480 23.975
	22.554 23.959
	22.627 23.943
	22.699 23.926
	 /
\plot 22.699 23.926 23.842 23.656 /
\linethickness=1pt
\setplotsymbol ({\makebox(0,0)[l]{\tencirc\symbol{'160}}})
%
%
%
\plot	24.382 23.497 24.718 23.354
 	24.753 23.230
	24.671 23.129
	24.605 23.068
	24.522 23.001
	24.422 22.927
	24.306 22.847
	24.241 22.804
	24.173 22.760
	24.100 22.714
	24.023 22.667
	23.942 22.617
	23.857 22.567
	23.767 22.514
	23.674 22.460
	23.576 22.404
	23.474 22.347
	23.368 22.288
	23.258 22.227
	23.151 22.166
	23.055 22.104
	22.970 22.043
	22.896 21.982
	22.780 21.861
	22.708 21.740
	22.679 21.620
	22.694 21.500
	22.753 21.381
	22.855 21.263
	22.922 21.204
	23.000 21.145
	23.089 21.086
	23.189 21.028
	23.300 20.969
	23.422 20.911
	23.487 20.882
	23.555 20.853
	23.625 20.824
	23.698 20.795
	23.774 20.767
	23.853 20.738
	23.934 20.709
	24.018 20.680
	24.105 20.651
	24.195 20.623
	24.287 20.594
	24.382 20.565
	24.479 20.537
	24.580 20.508
	24.683 20.480
	24.789 20.451
	24.897 20.423
	25.008 20.395
	25.122 20.366
	25.239 20.338
	 /
\plot 25.239 20.338 29.017 19.433 /
\linethickness= 0.500pt
\setplotsymbol ({\thinlinefont .})
%
%
%
\plot	18.159 25.466 19.969 25.323
 	20.082 25.314
	20.195 25.304
	20.309 25.294
	20.422 25.283
	20.536 25.272
	20.650 25.260
	20.763 25.247
	20.877 25.234
	20.992 25.221
	21.106 25.207
	21.220 25.193
	21.335 25.178
	21.450 25.162
	21.564 25.146
	21.679 25.129
	21.794 25.112
	21.909 25.095
	22.025 25.077
	22.140 25.058
	22.256 25.039
	22.371 25.019
	22.487 24.999
	22.603 24.978
	22.719 24.957
	22.835 24.935
	22.951 24.912
	23.068 24.890
	23.184 24.866
	23.301 24.842
	23.418 24.818
	23.535 24.793
	23.652 24.767
	23.765 24.741
	23.873 24.714
	23.973 24.686
	24.067 24.657
	24.154 24.628
	24.235 24.597
	24.309 24.566
	24.376 24.534
	24.491 24.467
	24.579 24.397
	24.640 24.324
	24.676 24.247
	24.684 24.167
	24.666 24.084
	24.622 23.997
	24.551 23.907
	24.453 23.813
	24.329 23.717
	24.257 23.667
	24.178 23.616
	24.093 23.565
	24.001 23.513
	23.907 23.461
	23.817 23.410
	23.731 23.361
	23.648 23.313
	23.569 23.266
	23.493 23.221
	23.421 23.177
	23.352 23.134
	23.225 23.052
	23.112 22.976
	23.014 22.904
	22.929 22.838
	22.803 22.722
	22.733 22.626
	22.719 22.551
	22.763 22.497
	 /
\plot 22.763 22.497 23.048 22.322 /
\linethickness= 0.500pt
\setplotsymbol ({\thinlinefont .})
%
%
%
\plot	23.683 21.973 25.207 21.465
 	25.302 21.434
	25.397 21.403
	25.490 21.372
	25.583 21.343
	25.676 21.313
	25.768 21.285
	25.859 21.257
	25.950 21.229
	26.040 21.202
	26.129 21.176
	26.218 21.150
	26.306 21.124
	26.393 21.100
	26.480 21.075
	26.566 21.052
	26.652 21.029
	26.737 21.006
	26.821 20.984
	26.905 20.963
	26.988 20.942
	27.071 20.921
	27.153 20.902
	27.234 20.882
	27.315 20.864
	27.395 20.846
	27.474 20.828
	27.553 20.811
	27.631 20.795
	27.709 20.779
	27.786 20.764
	27.862 20.749
	27.938 20.735
	 /
\plot 27.938 20.735 29.144 20.513 /
\linethickness= 0.500pt
\setplotsymbol ({\thinlinefont .})
%
%
%
\plot	23.910 14.201 25.227 14.066
 	25.310 14.058
	25.391 14.050
	25.472 14.042
	25.553 14.034
	25.633 14.027
	25.712 14.019
	25.791 14.012
	25.870 14.005
	25.948 13.998
	26.025 13.992
	26.102 13.985
	26.179 13.979
	26.255 13.973
	26.330 13.967
	26.405 13.961
	26.480 13.955
	26.554 13.950
	26.627 13.945
	26.700 13.940
	26.772 13.935
	26.844 13.930
	26.916 13.926
	26.987 13.922
	27.057 13.918
	27.127 13.914
	27.196 13.910
	27.265 13.907
	27.333 13.903
	27.401 13.900
	27.468 13.897
	27.535 13.894
	27.601 13.892
	 /
\plot 27.601 13.892 28.658 13.851 /
\linethickness=1pt
\setplotsymbol ({\makebox(0,0)[l]{\tencirc\symbol{'160}}})
%
%
%
\plot	 5.874 24.122  7.493 23.725
 	 7.594 23.699
	 7.696 23.673
	 7.797 23.646
	 7.899 23.618
	 8.000 23.589
	 8.102 23.560
	 8.204 23.529
	 8.306 23.497
	 8.408 23.465
	 8.510 23.432
	 8.612 23.397
	 8.714 23.362
	 8.817 23.326
	 8.919 23.289
	 9.022 23.251
	 9.124 23.213
	 9.227 23.173
	 9.330 23.133
	 9.433 23.091
	 9.536 23.049
	 9.639 23.005
	 9.742 22.961
	 9.845 22.916
	 9.949 22.870
	10.052 22.824
	10.156 22.776
	10.259 22.727
	10.363 22.678
	10.467 22.627
	10.571 22.576
	10.675 22.524
	10.779 22.471
	10.881 22.417
	10.978 22.365
	11.071 22.314
	11.158 22.264
	11.241 22.214
	11.320 22.166
	11.394 22.118
	11.463 22.072
	11.587 21.981
	11.692 21.895
	11.779 21.812
	11.847 21.732
	11.896 21.657
	11.926 21.585
	11.938 21.517
	11.931 21.453
	11.861 21.335
	11.798 21.282
	11.716 21.232
	11.625 21.187
	11.537 21.147
	11.450 21.112
	11.365 21.082
	11.281 21.058
	11.199 21.038
	11.119 21.023
	11.041 21.014
	10.965 21.010
	10.890 21.010
	10.817 21.016
	10.745 21.027
	10.676 21.043
	10.608 21.064
	10.542 21.090
	10.478 21.121
	10.365 21.191
	10.282 21.265
	10.228 21.344
	10.204 21.427
	10.208 21.514
	10.242 21.606
	10.306 21.703
	10.398 21.804
	 /
\plot 10.398 21.804 10.827 22.217 /
\linethickness=1pt
\setplotsymbol ({\makebox(0,0)[l]{\tencirc\symbol{'160}}})
%
%
%
\plot	11.239 22.566 12.129 23.169
 	12.242 23.243
	12.359 23.312
	12.481 23.378
	12.608 23.440
	12.673 23.469
	12.738 23.498
	12.805 23.526
	12.873 23.552
	12.943 23.578
	13.013 23.603
	13.084 23.626
	13.156 23.649
	13.230 23.671
	13.304 23.692
	13.380 23.712
	13.457 23.731
	13.535 23.749
	13.613 23.766
	13.693 23.782
	13.775 23.797
	13.857 23.811
	13.940 23.824
	14.024 23.837
	14.110 23.848
	14.196 23.858
	14.284 23.868
	14.372 23.876
	14.462 23.883
	 /
\plot 14.462 23.883 15.907 23.995 /
\linethickness= 0.500pt
\setplotsymbol ({\thinlinefont .})
%
%
%
\plot	17.990 18.773 20.013 18.643
 	20.139 18.634
	20.263 18.625
	20.386 18.614
	20.507 18.603
	20.627 18.591
	20.746 18.578
	20.863 18.564
	20.979 18.550
	21.093 18.534
	21.206 18.518
	21.317 18.501
	21.427 18.483
	21.535 18.464
	21.642 18.444
	21.748 18.424
	21.852 18.402
	21.955 18.380
	22.056 18.356
	22.156 18.332
	22.254 18.307
	22.351 18.282
	22.447 18.255
	22.541 18.227
	22.633 18.199
	22.725 18.170
	22.814 18.140
	22.903 18.109
	22.989 18.077
	23.075 18.044
	23.159 18.011
	23.241 17.976
	23.322 17.941
	23.402 17.904
	23.479 17.867
	23.555 17.828
	23.628 17.788
	23.699 17.746
	23.768 17.703
	23.836 17.659
	23.901 17.613
	24.025 17.518
	24.142 17.418
	24.250 17.312
	24.350 17.202
	24.443 17.086
	24.527 16.964
	24.603 16.838
	24.638 16.773
	24.671 16.706
	24.702 16.638
	24.731 16.569
	24.758 16.499
	24.783 16.427
	24.806 16.354
	24.827 16.280
	24.845 16.204
	24.862 16.127
	24.877 16.050
	24.887 15.975
	24.895 15.902
	24.898 15.831
	24.899 15.762
	24.896 15.695
	24.889 15.629
	24.879 15.566
	24.849 15.444
	24.805 15.331
	24.747 15.224
	24.676 15.125
	24.591 15.034
	24.492 14.950
	24.379 14.874
	24.252 14.805
	24.184 14.774
	24.112 14.744
	24.037 14.717
	23.958 14.691
	23.876 14.667
	23.790 14.645
	23.701 14.624
	23.608 14.606
	23.516 14.589
	23.428 14.572
	23.345 14.556
	23.265 14.540
	23.190 14.524
	23.119 14.509
	23.052 14.495
	22.990 14.480
	22.877 14.453
	22.782 14.428
	22.703 14.404
	22.642 14.381
	22.759 14.264
	 /
\plot 22.759 14.264 23.419 14.201 /
\linethickness=1pt
\setplotsymbol ({\makebox(0,0)[l]{\tencirc\symbol{'160}}})
%
%
%
\plot	18.004 17.678 20.513 17.476
 	20.590 17.470
	20.667 17.463
	20.743 17.456
	20.818 17.448
	20.891 17.440
	20.964 17.431
	21.035 17.422
	21.106 17.413
	21.175 17.403
	21.244 17.393
	21.311 17.383
	21.377 17.372
	21.442 17.361
	21.507 17.349
	21.632 17.324
	21.753 17.298
	21.869 17.270
	21.982 17.241
	22.090 17.210
	22.194 17.178
	22.294 17.143
	22.389 17.108
	22.481 17.070
	22.568 17.031
	22.651 16.991
	22.730 16.949
	22.804 16.905
	22.874 16.860
	22.940 16.813
	23.060 16.715
	23.163 16.610
	23.248 16.499
	23.317 16.381
	23.369 16.257
	23.411 16.133
	23.451 16.015
	23.487 15.904
	23.521 15.798
	23.553 15.699
	23.581 15.607
	23.607 15.520
	23.630 15.440
	23.651 15.366
	23.669 15.298
	23.696 15.181
	23.713 15.089
	23.718 15.022
	 /
\plot 23.718 15.022 23.719 14.804 /
\linethickness=1pt
\setplotsymbol ({\makebox(0,0)[l]{\tencirc\symbol{'160}}})
%
%
%
\plot	23.705 14.391 23.695 14.255
 	23.694 14.165
	23.695 14.103
	23.698 14.031
	23.703 13.947
	23.710 13.853
	23.718 13.748
	23.728 13.631
	23.776 13.516
	23.898 13.412
	23.987 13.365
	24.095 13.321
	24.221 13.280
	24.292 13.260
	24.367 13.242
	24.446 13.224
	24.530 13.207
	24.619 13.190
	24.713 13.175
	24.811 13.160
	24.914 13.146
	25.021 13.132
	25.133 13.120
	25.250 13.108
	25.371 13.097
	25.498 13.086
	25.562 13.081
	25.628 13.077
	25.695 13.072
	25.764 13.068
	25.833 13.064
	25.904 13.060
	25.976 13.056
	26.049 13.053
	26.123 13.049
	26.198 13.046
	 /
\plot 26.198 13.046 28.626 12.948 /
\linethickness= 0.500pt
\setplotsymbol ({\thinlinefont .})
%
%
%
\plot	 9.878 15.873  9.887 16.224
 	 9.898 16.310
	 9.927 16.391
	 9.975 16.468
	10.040 16.541
	10.124 16.609
	10.225 16.673
	10.345 16.732
	10.412 16.760
	10.483 16.787
	10.556 16.812
	10.628 16.834
	10.699 16.853
	10.770 16.869
	10.840 16.882
	10.908 16.892
	10.976 16.900
	11.044 16.904
	11.110 16.905
	11.176 16.903
	11.240 16.898
	11.304 16.890
	11.430 16.866
	11.552 16.830
	11.668 16.782
	11.778 16.724
	11.880 16.656
	11.975 16.577
	12.063 16.488
	12.143 16.388
	12.216 16.278
	12.282 16.158
	12.344 16.040
	12.405 15.937
	12.466 15.850
	12.527 15.779
	12.646 15.683
	12.764 15.650
	 /
\plot 12.764 15.650 12.996 15.644 /
\linethickness= 0.500pt
\setplotsymbol ({\thinlinefont .})
%
%
%
\plot	18.730 10.742 20.419 10.145
 	20.524 10.107
	20.625 10.069
	20.723 10.029
	20.819  9.988
	20.911  9.946
	21.001  9.903
	21.088  9.859
	21.172  9.814
	21.253  9.768
	21.331  9.721
	21.406  9.673
	21.478  9.623
	21.548  9.573
	21.614  9.522
	21.739  9.417
	21.852  9.307
	21.953  9.193
	22.043  9.075
	22.121  8.953
	22.188  8.827
	22.217  8.763
	22.243  8.697
	22.266  8.630
	22.286  8.563
	22.304  8.494
	22.318  8.424
	22.332  8.355
	22.347  8.287
	22.365  8.220
	22.383  8.155
	22.425  8.029
	22.473  7.910
	22.528  7.796
	22.589  7.688
	22.656  7.586
	22.729  7.489
	22.809  7.399
	22.895  7.315
	22.987  7.236
	23.086  7.163
	23.191  7.097
	23.302  7.036
	23.420  6.981
	23.544  6.932
	23.669  6.888
	23.791  6.849
	23.909  6.814
	24.024  6.784
	24.135  6.758
	24.243  6.736
	24.348  6.719
	24.449  6.707
	24.546  6.699
	24.640  6.695
	24.731  6.696
	24.818  6.701
	24.901  6.711
	24.981  6.725
	25.058  6.744
	25.131  6.767
	 /
\plot 25.131  6.767 25.703  6.970 /
\linethickness= 0.500pt
\setplotsymbol ({\thinlinefont .})
%
%
%
\plot	25.995  7.250 26.611  7.586
 	26.690  7.625
	26.773  7.659
	26.860  7.686
	26.951  7.708
	27.047  7.723
	27.146  7.733
	27.249  7.737
	27.357  7.735
	27.468  7.728
	27.584  7.714
	27.704  7.695
	27.828  7.669
	27.891  7.655
	27.955  7.638
	28.021  7.621
	28.087  7.601
	28.155  7.581
	28.223  7.559
	28.293  7.535
	28.363  7.510
	 /
\plot 28.363  7.510 29.500  7.097 /
\linethickness= 0.500pt
\setplotsymbol ({\thinlinefont .})
%
%
%
\plot	 8.871  8.769  9.009  8.764
 	 9.120  8.717
	 9.207  8.661
	 9.315  8.583
	 9.378  8.535
	 9.445  8.483
	 9.518  8.424
	 9.596  8.360
	 9.679  8.291
	 9.768  8.216
	 9.862  8.136
	 9.961  8.050
	10.062  7.966
	10.160  7.893
	10.255  7.831
	10.349  7.779
	10.439  7.738
	10.528  7.708
	10.613  7.688
	10.697  7.678
	10.777  7.679
	10.856  7.690
	10.931  7.712
	11.005  7.745
	11.076  7.788
	11.144  7.842
	11.210  7.906
	11.273  7.981
	 /
\plot 11.273  7.981 11.771  8.621 /
\linethickness= 0.500pt
\setplotsymbol ({\thinlinefont .})
%
%
%
\plot	11.968  8.894 12.460  9.513
 	12.522  9.586
	12.584  9.648
	12.648  9.699
	12.712  9.741
	12.777  9.772
	12.842  9.792
	12.909  9.803
	12.976  9.803
	13.044  9.793
	13.112  9.773
	13.181  9.742
	13.252  9.701
	13.322  9.650
	13.394  9.588
	13.466  9.516
	13.539  9.434
	13.611  9.349
	13.681  9.270
	13.747  9.197
	13.811  9.129
	13.931  9.009
	14.041  8.911
	14.139  8.836
	14.227  8.782
	14.304  8.750
	14.370  8.741
	 /
\plot 14.370  8.741 14.613  8.746 /
\linethickness= 0.500pt
\setplotsymbol ({\thinlinefont .})
%
%
%
\plot	10.668 15.613 10.811 15.607
 	10.920 15.521
	10.950 15.418
	10.958 15.351
	10.962 15.274
	10.960 15.195
	10.950 15.124
	10.903 15.001
	10.824 14.906
	10.710 14.838
	10.645 14.814
	10.578 14.797
	10.511 14.785
	10.442 14.781
	10.372 14.782
	10.302 14.790
	10.230 14.804
	10.158 14.825
	10.089 14.852
	10.029 14.884
	 9.936 14.965
	 9.877 15.067
	 9.854 15.192
	 /
\plot  9.854 15.192  9.842 15.464 /
%
%
\put{\SetFigFont{12}{14.4}{rm}\bullet} [lB] <0pt, -3pt> at 23.656  7.882
%
%
\put{\SetFigFont{12}{14.4}{rm}D} [lB] at 23.624  8.289
%
%
\put{\SetFigFont{12}{14.4}{rm}\bullet} [lB] <0pt, -2.5pt> at 11.191  8.750
%
%
\put{\SetFigFont{12}{14.4}{rm}D} [lB] at 11.096  9.100
%
%
\put{\SetFigFont{12}{14.4}{rm}H} [lB] at 22.735 14.994
%
%
\put{\SetFigFont{12}{14.4}{rm}H} [lB] <0pt, -2pt> at 11.034 15.081
%
%
\put{\SetFigFont{12}{14.4}{rm}\bullet} [lB] <-0.5pt, 0pt> at 23.497 15.075
%
%
\put{\SetFigFont{12}{14.4}{rm}\bullet} [lB] <0pt, -2pt> at 11.157 15.621
\linethickness=0pt
\putrectangle corners at  5.842 25.483 and 29.517  5.861
\endpicture
\endfig

The left hand sequence in figure \Whit\ is the view
from the top, with the curve drawn out by the tip of the vector 
shown (in most places the vector is pointing up so this curve 
coincides with the submanifold), and the right hand sequence is
the view from the side.  The arrows in the bottom two pictures
of the left sequence indicate slope.  Arrows point {\it downhill}.
Reading up the page gives an equivalent
compression isotopy to that pictured in figures \vectp\ to \twist.
The whole sequence defines an embedding of $\re^2$ in $\re^4$
with normal vector field which projects to a map of $\re^2$ to $\re^3$.
The image has a line of double points at the top terminating in
a singular point (the image of $H$) which can be seen to be a 
Whitney umbrella.  The surface is flat at the bottom and has
a `ripple' at the top (ie, a region of the form (curve shaped like the
letter $\alpha$)$\times I$).   The ripple shrinks to a point
at the umbrella.

If we now apply the compression theorem to this embedding with
vector field, then the image changes into a sheet with a continuous
ripple (a twist is created in each slice from the $H$--slice down).
Thus the Whitney umbrella is desingularised by having a new ripple
spliced into it.

\sh{Cross-cap to Boy's surface}

We finish by describing a complete compression.  In this example
we apply the compression theorem to an immersion (rather than an
embedding).  We shall see in the next section that the compression
theorem can be applied to any immersion with a normal vector field
by working locally.  Here we anticipate this.  Consider the
well-known non-immersion of $P^2$ in $\re^3$ with a line of 
double points and two singularities (Whitney umbrellas) 
--- a 2--sphere with cross-cap --- illustrated in figure \figkey\CCap.

\fig{\CCap: non-immersion of $P^2$ in $\re^3$}
\beginpicture
\setcoordinatesystem units <0.80000cm,0.80000cm>
\unitlength=1.00000cm
\linethickness=1pt
\setplotsymbol ({\makebox(0,0)[l]{\tencirc\symbol{'160}}})
\setshadesymbol ({\thinlinefont .})
\setlinear
%
%
\linethickness= 0.500pt
\setplotsymbol ({\thinlinefont .})
\circulararc 95.709 degrees from  4.047 25.559 center at  4.938 26.410
%
%
\linethickness= 0.500pt
\setplotsymbol ({\thinlinefont .})
\circulararc 95.940 degrees from  5.906 25.590 center at  6.725 26.373
%
%
\linethickness= 0.500pt
\setplotsymbol ({\thinlinefont .})
\circulararc 74.957 degrees from  2.904 24.479 center at  4.428 26.401
%
%
\linethickness= 0.500pt
\setplotsymbol ({\thinlinefont .})
\circulararc 80.607 degrees from  5.920 24.431 center at  7.273 26.197
%
%
\linethickness= 0.500pt
\setplotsymbol ({\thinlinefont .})
\circulararc 91.749 degrees from  2.366 23.099 center at  4.166 24.736
%
%
\linethickness= 0.500pt
\setplotsymbol ({\thinlinefont .})
\circulararc 92.846 degrees from  5.889 23.004 center at  7.529 24.730
%
%
\linethickness= 0.500pt
\setplotsymbol ({\thinlinefont .})
\circulararc 79.470 degrees from  2.381 21.876 center at  4.192 23.938
%
%
\linethickness= 0.500pt
\setplotsymbol ({\thinlinefont .})
\circulararc 81.254 degrees from  5.874 21.795 center at  7.547 23.878
%
%
\linethickness= 0.500pt
\setplotsymbol ({\thinlinefont .})
\circulararc 39.867 degrees from  4.064 19.511 center at  5.841 24.884
%
%
\linethickness= 0.500pt
\setplotsymbol ({\thinlinefont .})
\setdots < 0.0953cm>
\circulararc 52.980 degrees from  5.857 21.795 center at  4.079 18.350
%
%
\linethickness= 0.500pt
\setplotsymbol ({\thinlinefont .})
\circulararc 44.251 degrees from  9.351 21.891 center at  7.755 17.577
%
%
\linethickness= 0.500pt
\setplotsymbol ({\thinlinefont .})
\circulararc 60.428 degrees from  5.874 23.019 center at  4.094 20.059
%
%
\linethickness= 0.500pt
\setplotsymbol ({\thinlinefont .})
\circulararc 67.087 degrees from  9.351 23.194 center at  7.740 20.503
%
%
\linethickness= 0.500pt
\setplotsymbol ({\thinlinefont .})
\circulararc 78.784 degrees from  5.906 24.448 center at  4.375 22.662
%
%
\linethickness= 0.500pt
\setplotsymbol ({\thinlinefont .})
\circulararc 91.124 degrees from  8.795 24.575 center at  7.412 23.078
%
%
\linethickness= 0.500pt
\setplotsymbol ({\thinlinefont .})
\circulararc 67.637 degrees from  5.889 25.654 center at  5.047 24.295
%
%
\linethickness= 0.500pt
\setplotsymbol ({\thinlinefont .})
\circulararc 56.090 degrees from  7.588 25.639 center at  6.768 24.028
%
%
\linethickness= 0.500pt
\setplotsymbol ({\thinlinefont .})
\circulararc 41.113 degrees from  7.889 19.607 center at  6.103 14.459
%
%
\linethickness= 0.500pt
\setplotsymbol ({\thinlinefont .})
\setsolid
\ellipticalarc axes ratio  3.539:3.539  360 degrees 
	from  9.413 22.543 center at  5.874 22.543
%
%
\linethickness= 0.500pt
\setplotsymbol ({\thinlinefont .})
\putrule from  5.874 26.099 to  5.874 21.780
\linethickness= 0.500pt
\setplotsymbol ({\thinlinefont .})
%
%
%
\plot	 2.587 21.209  3.412 20.749
 	 3.515 20.695
	 3.619 20.649
	 3.723 20.610
	 3.827 20.579
	 3.931 20.556
	 4.036 20.541
	 4.140 20.533
	 4.246 20.532
	 4.351 20.540
	 4.456 20.555
	 4.562 20.577
	 4.668 20.608
	 4.774 20.645
	 4.881 20.691
	 4.988 20.744
	 5.095 20.805
	 5.202 20.865
	 5.309 20.919
	 5.417 20.964
	 5.524 21.002
	 5.632 21.032
	 5.740 21.055
	 5.848 21.070
	 5.956 21.077
	 6.064 21.077
	 6.173 21.069
	 6.281 21.053
	 6.390 21.030
	 6.498 20.999
	 6.607 20.961
	 6.716 20.915
	 6.825 20.861
	 6.934 20.807
	 7.041 20.759
	 7.148 20.718
	 7.253 20.684
	 7.358 20.657
	 7.461 20.636
	 7.563 20.622
	 7.665 20.615
	 7.765 20.614
	 7.864 20.620
	 7.962 20.633
	 8.059 20.652
	 8.155 20.679
	 8.250 20.711
	 8.344 20.751
	 8.437 20.797
	 /
\plot  8.437 20.797  9.176 21.194 /
\linethickness= 0.500pt
\setplotsymbol ({\thinlinefont .})
%
%
%
\plot	 2.999 20.447  3.905 20.043
 	 4.017 19.995
	 4.127 19.951
	 4.235 19.913
	 4.340 19.879
	 4.443 19.850
	 4.543 19.825
	 4.641 19.806
	 4.737 19.791
	 4.831 19.781
	 4.922 19.776
	 5.011 19.775
	 5.098 19.780
	 5.182 19.789
	 5.265 19.803
	 5.344 19.821
	 5.422 19.845
	 5.498 19.869
	 5.575 19.890
	 5.651 19.909
	 5.727 19.924
	 5.804 19.937
	 5.880 19.946
	 5.956 19.952
	 6.032 19.956
	 6.109 19.956
	 6.185 19.954
	 6.261 19.948
	 6.338 19.940
	 6.414 19.929
	 6.490 19.914
	 6.567 19.897
	 6.643 19.877
	 6.720 19.857
	 6.799 19.842
	 6.879 19.833
	 6.960 19.827
	 7.043 19.827
	 7.127 19.831
	 7.213 19.841
	 7.300 19.855
	 7.388 19.873
	 7.478 19.897
	 7.570 19.925
	 7.663 19.958
	 7.757 19.996
	 7.853 20.038
	 7.950 20.086
	 8.049 20.138
	 /
\plot  8.049 20.138  8.843 20.574 /
\linethickness= 0.500pt
\setplotsymbol ({\thinlinefont .})
\setdots < 0.0953cm>
%
%
%
\plot	 2.555 21.258  3.484 21.551
 	 3.600 21.585
	 3.714 21.615
	 3.827 21.640
	 3.940 21.660
	 4.051 21.676
	 4.161 21.687
	 4.270 21.694
	 4.377 21.695
	 4.484 21.692
	 4.589 21.685
	 4.694 21.673
	 4.797 21.656
	 4.899 21.635
	 5.000 21.608
	 5.100 21.578
	 5.199 21.542
	 5.297 21.507
	 5.397 21.477
	 5.498 21.452
	 5.601 21.431
	 5.704 21.416
	 5.809 21.405
	 5.915 21.400
	 6.022 21.400
	 6.130 21.404
	 6.239 21.413
	 6.350 21.428
	 6.462 21.447
	 6.575 21.471
	 6.689 21.501
	 6.804 21.535
	 6.920 21.574
	 7.036 21.613
	 7.150 21.647
	 7.261 21.675
	 7.370 21.697
	 7.477 21.714
	 7.582 21.726
	 7.684 21.732
	 7.784 21.733
	 7.881 21.728
	 7.977 21.718
	 8.070 21.702
	 8.160 21.681
	 8.249 21.655
	 8.335 21.623
	 8.419 21.585
	 8.501 21.542
	 /
\plot  8.501 21.542  9.144 21.177 /
\linethickness= 0.500pt
\setplotsymbol ({\thinlinefont .})
%
%
%
\plot	 2.999 20.447  3.722 20.701
 	 3.813 20.731
	 3.904 20.757
	 3.996 20.780
	 4.088 20.800
	 4.182 20.816
	 4.275 20.828
	 4.370 20.837
	 4.465 20.842
	 4.560 20.844
	 4.656 20.842
	 4.753 20.837
	 4.850 20.828
	 4.948 20.815
	 5.047 20.800
	 5.146 20.780
	 5.246 20.757
	 5.346 20.734
	 5.446 20.714
	 5.546 20.697
	 5.645 20.684
	 5.745 20.674
	 5.844 20.667
	 5.943 20.663
	 6.041 20.662
	 6.140 20.665
	 6.239 20.671
	 6.337 20.680
	 6.435 20.692
	 6.533 20.707
	 6.630 20.726
	 6.728 20.748
	 6.825 20.773
	 6.922 20.798
	 7.017 20.819
	 7.112 20.838
	 7.205 20.852
	 7.297 20.863
	 7.387 20.871
	 7.477 20.875
	 7.565 20.876
	 7.653 20.873
	 7.739 20.867
	 7.824 20.857
	 7.907 20.844
	 7.990 20.828
	 8.071 20.808
	 8.152 20.784
	 8.231 20.757
	 /
\plot  8.231 20.757  8.858 20.527 /
\linethickness=0pt
\putrectangle corners at  2.318 26.124 and  9.430 18.989
\endpicture
\endfig

This is covered by an embedding in $\re^4$ which we can think of
as a sequence of 3--dimensional slices, which start with a
small circle in the shape of a `figure 8' (corresponding
to a level just below the top singular point), move to a genuine
`figure 8' at the bottom singular point and then resolve into
a circle.   
Now it is an old result of Whitney that no embedding in $\re^4$ admits
a normal vector field, so we modify this embedding to an immersion
with a single double point in order to construct a normal vector
field.  The result, in a more general projection, together with the 
vector field is illustrated in figure \figkey\PVect\ 
below.

\fig{\PVect: covering immersion in $\re^4$}
\beginpicture
\setcoordinatesystem units <0.50000cm,0.50000cm>
\unitlength=1.00000cm
\linethickness=1pt
\setplotsymbol ({\makebox(0,0)[l]{\tencirc\symbol{'160}}})
\setshadesymbol ({\thinlinefont .})
\setlinear
%
%
\linethickness= 0.500pt
\setplotsymbol ({\thinlinefont .})
\ellipticalarc axes ratio  1.937:0.667  360 degrees 
	from 14.743 24.132 center at 12.806 24.132
%
%
\linethickness= 0.500pt
\setplotsymbol ({\thinlinefont .})
\ellipticalarc axes ratio  0.667:0.318  360 degrees 
	from 13.494 25.559 center at 12.827 25.559
%
%
\linethickness= 0.500pt
\setplotsymbol ({\thinlinefont .})
\ellipticalarc axes ratio  2.635:0.953  360 degrees 
	from 15.102  3.827 center at 12.467  3.827
%
%
\linethickness= 0.500pt
\setplotsymbol ({\thinlinefont .})
\ellipticalarc axes ratio  1.016:0.444  360 degrees 
	from 13.610  1.827 center at 12.594  1.827
%
%
\linethickness= 0.500pt
\setplotsymbol ({\thinlinefont .})
\ellipticalarc axes ratio  2.517:0.851  360 degrees 
	from 15.238 22.011 center at 12.721 22.011
%
%
\linethickness= 0.500pt
\setplotsymbol ({\thinlinefont .})
\ellipticalarc axes ratio  2.961:1.094  360 degrees 
	from 15.507  6.555 center at 12.545  6.555
%
%
\linethickness= 0.500pt
\setplotsymbol ({\thinlinefont .})
\plot 12.150 23.592 12.006 23.529 /
\plot 12.006 23.529 12.133 23.434 /
%
%
\linethickness= 0.500pt
\setplotsymbol ({\thinlinefont .})
\plot 13.561 23.609 13.705 23.546 /
\plot 13.705 23.546 13.593 23.451 /
%
%
\linethickness= 0.500pt
\setplotsymbol ({\thinlinefont .})
\plot 11.451 18.633 11.301 18.555 /
\plot 11.301 18.555 11.402 18.411 /
%
%
\linethickness= 0.500pt
\setplotsymbol ({\thinlinefont .})
\setdots < 0.0953cm>
\plot 11.995 14.374 12.979 14.406 /
%
%
\linethickness= 0.500pt
\setplotsymbol ({\thinlinefont .})
\setsolid
\plot 10.615 12.025 10.829 11.968 /
\plot 10.829 11.968 10.679 11.849 /
%
%
\linethickness= 0.500pt
\setplotsymbol ({\thinlinefont .})
\plot 11.496  3.086 11.646  2.913 /
\plot 11.646  2.913 11.481  2.849 /
%
%
\linethickness= 0.500pt
\setplotsymbol ({\thinlinefont .})
\plot 14.575  3.397 14.497  3.213 /
\plot 14.497  3.213 14.662  3.190 /
%
%
\linethickness= 0.500pt
\setplotsymbol ({\thinlinefont .})
\plot 11.151  8.352 11.278  8.242 /
\plot 11.278  8.242 11.159  8.153 /
%
%
\linethickness= 0.500pt
\setplotsymbol ({\thinlinefont .})
\plot 11.942 21.296 11.769 21.241 /
\plot 11.769 21.241 11.879 21.105 /
%
%
\linethickness= 0.500pt
\setplotsymbol ({\thinlinefont .})
\plot 13.379 21.258 13.555 21.209 /
\plot 13.555 21.209 13.428 21.105 /
%
%
\linethickness= 0.500pt
\setplotsymbol ({\thinlinefont .})
\plot 11.212  5.717 11.331  5.550 /
\plot 11.331  5.550 11.148  5.503 /
%
%
\linethickness= 0.500pt
\setplotsymbol ({\thinlinefont .})
\plot 14.618  5.925 14.522  5.734 /
\putrule from 14.522  5.734 to 14.730  5.734
\linethickness= 0.500pt
\setplotsymbol ({\thinlinefont .})
%
%
%
\plot	 9.808 17.670 
 	 9.876 17.718
	 9.946 17.760
	10.017 17.797
	10.090 17.828
	10.164 17.854
	10.239 17.874
	10.316 17.889
	10.394 17.899
	10.474 17.903
	10.555 17.901
	10.637 17.894
	10.721 17.882
	10.807 17.864
	10.893 17.840
	10.982 17.811
	11.071 17.777
	11.162 17.737
	11.255 17.691
	11.348 17.640
	11.444 17.584
	11.540 17.522
	11.638 17.455
	11.738 17.382
	11.839 17.303
	11.941 17.220
	12.045 17.130
	12.150 17.036
	12.256 16.935
	12.364 16.829
	12.474 16.718
	12.585 16.602
	12.697 16.479
	12.809 16.357
	12.919 16.242
	13.027 16.132
	13.134 16.028
	13.239 15.930
	13.342 15.839
	13.443 15.754
	13.543 15.674
	13.640 15.601
	13.736 15.534
	13.830 15.472
	13.923 15.417
	14.013 15.368
	14.102 15.325
	14.189 15.289
	14.274 15.258
	14.357 15.233
	14.439 15.214
	14.519 15.202
	14.597 15.195
	14.673 15.195
	14.747 15.201
	14.820 15.212
	14.891 15.230
	14.960 15.254
	15.027 15.284
	15.092 15.320
	15.156 15.362
	15.278 15.465
	15.392 15.591
	15.446 15.661
	15.496 15.730
	15.542 15.798
	15.585 15.867
	15.623 15.935
	15.658 16.003
	15.689 16.071
	15.717 16.138
	15.741 16.205
	15.760 16.272
	15.777 16.338
	15.789 16.405
	15.798 16.471
	15.803 16.536
	15.804 16.602
	15.801 16.667
	15.795 16.732
	15.785 16.796
	15.771 16.861
	15.753 16.925
	15.732 16.988
	15.707 17.052
	15.645 17.178
	15.569 17.303
	15.477 17.427
	15.371 17.549
	15.250 17.671
	15.184 17.729
	15.117 17.781
	15.049 17.827
	14.979 17.867
	14.907 17.902
	14.834 17.930
	14.759 17.953
	14.682 17.971
	14.604 17.982
	14.524 17.988
	14.443 17.988
	14.360 17.982
	14.275 17.970
	14.189 17.953
	14.101 17.930
	14.011 17.901
	13.920 17.867
	13.827 17.826
	13.733 17.780
	13.637 17.728
	13.540 17.670
	13.441 17.607
	13.340 17.538
	13.237 17.463
	13.133 17.382
	13.028 17.295
	12.921 17.203
	12.812 17.105
	12.701 17.001
	12.589 16.892
	12.476 16.776
	12.360 16.655
	12.245 16.534
	12.132 16.419
	12.021 16.310
	11.911 16.207
	11.803 16.110
	11.698 16.019
	11.594 15.934
	11.492 15.855
	11.392 15.782
	11.294 15.716
	11.198 15.655
	11.103 15.600
	11.011 15.551
	10.920 15.508
	10.832 15.471
	10.745 15.441
	10.660 15.416
	10.577 15.397
	10.496 15.384
	10.417 15.378
	10.340 15.377
	10.264 15.382
	10.191 15.393
	10.119 15.411
	10.049 15.434
	 9.981 15.463
	 9.915 15.499
	 9.851 15.540
	 9.729 15.641
	 9.614 15.766
	 9.560 15.834
	 9.510 15.902
	 9.465 15.969
	 9.422 16.035
	 9.384 16.101
	 9.350 16.167
	 9.319 16.231
	 9.293 16.296
	 9.251 16.423
	 9.224 16.547
	 9.213 16.670
	 9.218 16.790
	 9.238 16.908
	 9.273 17.023
	 9.323 17.137
	 9.389 17.248
	 9.471 17.357
	 9.568 17.463
	 9.680 17.568
	 9.742 17.619
	 9.808 17.670
	/
\linethickness= 0.500pt
\setplotsymbol ({\thinlinefont .})
%
%
%
\plot	12.531 13.155 12.959 13.545
 	13.054 13.642
	13.125 13.738
	13.171 13.833
	13.193 13.928
	13.191 14.022
	13.165 14.115
	13.114 14.207
	13.039 14.299
	12.976 14.383
	12.961 14.451
	13.078 14.541
	13.138 14.554
	13.210 14.563
	13.294 14.569
	13.390 14.570
	13.498 14.568
	13.619 14.561
	13.683 14.557
	13.751 14.551
	13.822 14.545
	13.896 14.537
	13.971 14.529
	14.045 14.519
	14.119 14.509
	14.192 14.497
	14.264 14.484
	14.336 14.470
	14.406 14.455
	14.476 14.439
	14.546 14.422
	14.614 14.404
	14.682 14.384
	14.749 14.364
	14.816 14.342
	14.882 14.320
	14.947 14.296
	15.011 14.271
	15.138 14.219
	15.261 14.161
	15.382 14.100
	15.500 14.034
	15.615 13.964
	15.728 13.890
	15.837 13.811
	15.944 13.728
	16.041 13.642
	16.122 13.554
	16.187 13.466
	16.235 13.375
	16.268 13.284
	16.285 13.191
	16.285 13.097
	16.269 13.001
	16.237 12.904
	16.189 12.806
	16.125 12.706
	16.045 12.605
	15.949 12.503
	15.836 12.399
	15.774 12.347
	15.708 12.294
	15.637 12.241
	15.563 12.188
	15.486 12.136
	15.409 12.088
	15.332 12.045
	15.255 12.005
	15.177 11.969
	15.099 11.937
	15.020 11.908
	14.942 11.884
	14.863 11.864
	14.783 11.847
	14.704 11.835
	14.624 11.826
	14.543 11.821
	14.463 11.821
	14.382 11.824
	14.301 11.831
	14.219 11.841
	14.137 11.856
	14.055 11.875
	13.973 11.898
	13.890 11.924
	13.807 11.954
	13.724 11.989
	13.640 12.027
	13.556 12.069
	13.472 12.115
	13.387 12.165
	13.302 12.219
	13.217 12.277
	13.131 12.338
	13.045 12.404
	12.959 12.474
	12.875 12.544
	12.793 12.614
	12.716 12.682
	12.641 12.749
	12.570 12.815
	12.502 12.879
	12.377 13.004
	12.265 13.124
	12.167 13.239
	12.082 13.349
	12.011 13.454
	11.953 13.553
	11.909 13.648
	11.878 13.738
	11.861 13.822
	11.857 13.901
	11.867 13.975
	11.891 14.045
	11.927 14.109
	11.993 14.226
	12.016 14.334
	11.995 14.431
	11.931 14.517
	11.824 14.594
	11.755 14.628
	11.674 14.659
	11.583 14.688
	11.481 14.715
	11.368 14.738
	11.245 14.760
	11.180 14.768
	11.116 14.775
	11.051 14.780
	10.986 14.783
	10.921 14.784
	10.856 14.783
	10.790 14.780
	10.725 14.775
	10.659 14.768
	10.594 14.760
	10.528 14.749
	10.462 14.736
	10.396 14.721
	10.330 14.704
	10.263 14.685
	10.197 14.664
	10.131 14.641
	10.064 14.617
	 9.997 14.590
	 9.930 14.561
	 9.863 14.530
	 9.796 14.498
	 9.729 14.463
	 9.661 14.426
	 9.594 14.387
	 9.526 14.347
	 9.458 14.304
	 9.390 14.259
	 9.322 14.213
	 9.254 14.164
	 9.186 14.114
	 9.118 14.061
	 9.051 14.008
	 8.989 13.954
	 8.877 13.848
	 8.782 13.743
	 8.704 13.639
	 8.642 13.536
	 8.597 13.434
	 8.569 13.332
	 8.558 13.232
	 8.563 13.132
	 8.585 13.033
	 8.624 12.935
	 8.680 12.838
	 8.752 12.741
	 8.841 12.646
	 8.947 12.551
	 9.070 12.458
	 9.135 12.412
	 9.200 12.368
	 9.264 12.327
	 9.327 12.287
	 9.452 12.214
	 9.574 12.148
	 9.694 12.091
	 9.811 12.041
	 9.926 11.999
	10.038 11.966
	10.148 11.940
	10.256 11.922
	10.361 11.912
	10.463 11.910
	10.563 11.916
	10.660 11.930
	10.755 11.952
	10.848 11.981
	10.938 12.015
	11.024 12.049
	11.107 12.082
	11.186 12.116
	11.262 12.149
	11.335 12.183
	11.405 12.216
	11.471 12.249
	11.593 12.315
	11.702 12.381
	11.798 12.447
	11.880 12.513
	 /
\plot 11.880 12.513 12.181 12.774 /
\linethickness= 0.500pt
\setplotsymbol ({\thinlinefont .})
%
%
%
\plot	12.785  9.028 13.040  9.272
 	13.131  9.392
	13.152  9.507
	13.100  9.617
	12.977  9.723
	12.897  9.767
	12.817  9.797
	12.735  9.813
	12.652  9.814
	12.568  9.801
	12.483  9.774
	12.397  9.732
	12.311  9.675
	12.244  9.603
	12.217  9.514
	12.232  9.408
	12.286  9.285
	12.329  9.217
	12.382  9.144
	12.445  9.068
	12.518  8.987
	12.602  8.902
	12.695  8.813
	12.799  8.719
	12.913  8.621
	13.033  8.528
	13.156  8.450
	13.281  8.387
	13.344  8.361
	13.408  8.338
	13.473  8.319
	13.538  8.304
	13.604  8.292
	13.671  8.284
	13.738  8.280
	13.806  8.279
	13.875  8.282
	13.944  8.289
	14.014  8.299
	14.084  8.313
	14.155  8.331
	14.226  8.352
	14.299  8.377
	14.372  8.406
	14.445  8.438
	14.519  8.474
	14.594  8.514
	14.669  8.557
	14.745  8.604
	14.822  8.654
	14.899  8.709
	14.977  8.766
	15.055  8.828
	15.134  8.893
	15.210  8.960
	15.279  9.027
	15.341  9.093
	15.397  9.160
	15.445  9.226
	15.487  9.292
	15.521  9.357
	15.549  9.423
	15.570  9.488
	15.584  9.553
	15.591  9.618
	15.591  9.683
	15.584  9.747
	15.571  9.811
	15.550  9.875
	15.523  9.939
	15.489 10.003
	15.448 10.066
	15.345 10.193
	15.283 10.256
	15.215 10.318
	15.139 10.381
	15.057 10.443
	14.968 10.505
	14.871 10.567
	14.768 10.628
	14.658 10.690
	14.542 10.751
	14.418 10.812
	14.354 10.843
	14.287 10.873
	14.220 10.903
	14.150 10.934
	14.080 10.963
	14.010 10.992
	13.940 11.020
	13.871 11.047
	13.802 11.073
	13.733 11.098
	13.665 11.122
	13.597 11.145
	13.530 11.167
	13.463 11.188
	13.396 11.208
	13.329 11.227
	13.263 11.245
	13.197 11.262
	13.132 11.279
	13.066 11.294
	13.002 11.308
	12.937 11.321
	12.873 11.334
	12.809 11.345
	12.683 11.365
	12.557 11.381
	12.434 11.393
	12.311 11.401
	12.190 11.405
	12.070 11.406
	11.952 11.402
	11.835 11.395
	11.719 11.384
	11.605 11.369
	11.492 11.350
	11.380 11.327
	11.270 11.301
	11.161 11.270
	11.054 11.236
	10.948 11.197
	10.843 11.155
	10.740 11.109
	10.638 11.059
	10.537 11.005
	10.438 10.948
	10.340 10.886
	10.245 10.823
	10.156 10.759
	10.071 10.696
	 9.992 10.633
	 9.917 10.570
	 9.848 10.508
	 9.725 10.382
	 9.622 10.258
	 9.539 10.134
	 9.477 10.011
	 9.435  9.888
	 9.413  9.765
	 9.412  9.644
	 9.431  9.523
	 9.471  9.402
	 9.531  9.282
	 9.611  9.162
	 9.711  9.043
	 9.832  8.925
	 9.897  8.867
	 9.962  8.812
	10.026  8.759
	10.090  8.709
	10.215  8.616
	10.339  8.533
	10.460  8.461
	10.580  8.398
	10.697  8.346
	10.812  8.303
	10.925  8.271
	11.036  8.249
	11.145  8.237
	11.252  8.235
	11.356  8.244
	11.459  8.262
	11.559  8.290
	11.658  8.329
	 /
\plot 11.658  8.329 12.435  8.678 /
\linethickness= 0.500pt
\setplotsymbol ({\thinlinefont .})
%
%
%
\plot	11.148 23.789 11.220 23.753
 	11.331 23.776
	11.415 23.831
	11.519 23.916
	11.622 23.988
	11.706 24.004
	11.771 23.964
	11.816 23.869
	 /
\plot 11.816 23.869 11.887 23.622 /
\linethickness= 0.500pt
\setplotsymbol ({\thinlinefont .})
%
%
%
\plot	11.927 23.440 12.075 23.218
 	12.149 23.136
	12.226 23.112
	12.306 23.146
	12.388 23.238
	12.464 23.343
	12.527 23.417
	12.613 23.472
	 /
\plot 12.613 23.472 12.672 23.463 /
\linethickness= 0.500pt
\setplotsymbol ({\thinlinefont .})
%
%
%
\plot	13.053 23.463 13.121 23.472
 	13.243 23.408
	13.345 23.323
	13.407 23.267
	13.475 23.202
	13.543 23.141
	13.604 23.098
	13.708 23.066
	13.785 23.104
	13.836 23.213
	 /
\plot 13.836 23.213 13.911 23.503 /
\linethickness= 0.500pt
\setplotsymbol ({\thinlinefont .})
%
%
%
\plot	13.928 23.671 13.980 23.944
 	14.018 24.051
	14.082 24.100
	14.171 24.091
	14.285 24.023
	14.400 23.941
	14.489 23.888
	14.590 23.869
	 /
\plot 14.590 23.869 14.641 23.908 /
\linethickness= 0.500pt
\setplotsymbol ({\thinlinefont .})
%
%
%
\plot	10.894 21.440 10.994 21.404
 	11.052 21.404
	11.127 21.440
	11.219 21.513
	11.327 21.622
	11.433 21.716
	11.516 21.747
	11.576 21.712
	11.613 21.613
	 /
\plot 11.613 21.613 11.665 21.351 /
\linethickness= 0.500pt
\setplotsymbol ({\thinlinefont .})
%
%
%
\plot	11.688 21.145 11.752 20.903
 	11.849 20.786
	11.923 20.821
	12.014 20.919
	12.105 21.032
	12.178 21.113
	12.272 21.177
	 /
\plot 12.272 21.177 12.332 21.177 /
\linethickness= 0.500pt
\setplotsymbol ({\thinlinefont .})
%
%
%
\plot	13.974 21.264 14.062 21.284
 	14.164 21.269
	14.225 21.236
	14.293 21.185
	14.339 21.123
	14.334 21.056
	14.277 20.986
	14.170 20.911
	14.050 20.855
	13.955 20.843
	13.841 20.951
	 /
\plot 13.841 20.951 13.777 21.145 /
\linethickness= 0.500pt
\setplotsymbol ({\thinlinefont .})
%
%
%
\plot	13.737 21.336 13.741 21.475
 	13.800 21.585
	13.870 21.618
	13.968 21.638
	14.077 21.641
	14.183 21.627
	14.285 21.596
	14.384 21.547
	14.473 21.498
	14.545 21.468
	14.641 21.463
	 /
\plot 14.641 21.463 14.704 21.495 /
\linethickness= 0.500pt
\setplotsymbol ({\thinlinefont .})
%
%
%
\plot	10.427 18.959 10.494 18.916
 	10.544 18.910
	10.624 18.935
	10.734 18.993
	10.801 19.033
	10.875 19.082
	10.948 19.126
	11.011 19.155
	11.104 19.164
	11.154 19.109
	11.161 18.990
	 /
\plot 11.161 18.990 11.134 18.688 /
\linethickness= 0.500pt
\setplotsymbol ({\thinlinefont .})
%
%
%
\plot	14.355 18.563 14.530 18.357
 	14.621 18.269
	14.720 18.214
	14.828 18.190
	14.943 18.198
	15.040 18.239
	15.093 18.313
	15.103 18.421
	15.092 18.488
	15.070 18.563
	15.026 18.635
	14.950 18.691
	14.841 18.732
	14.775 18.747
	14.701 18.758
	14.618 18.765
	14.528 18.768
	14.429 18.767
	14.322 18.763
	14.208 18.754
	14.085 18.742
	14.021 18.734
	13.954 18.726
	13.886 18.716
	13.815 18.706
	13.744 18.695
	13.674 18.685
	13.604 18.675
	13.534 18.665
	13.466 18.655
	13.397 18.646
	13.330 18.636
	13.263 18.628
	13.196 18.619
	13.131 18.611
	13.065 18.603
	13.001 18.595
	12.937 18.588
	12.873 18.581
	12.748 18.567
	12.625 18.555
	12.505 18.544
	12.386 18.533
	12.271 18.524
	12.157 18.517
	12.046 18.510
	11.937 18.504
	11.831 18.500
	11.726 18.498
	11.622 18.502
	11.518 18.511
	11.415 18.525
	11.313 18.545
	11.211 18.570
	11.110 18.600
	11.010 18.635
	10.910 18.675
	10.810 18.720
	10.712 18.771
	10.614 18.827
	10.516 18.888
	10.419 18.955
	10.323 19.026
	10.228 19.103
	10.139 19.182
	10.065 19.259
	10.004 19.335
	 9.957 19.410
	 9.924 19.484
	 9.904 19.557
	 9.898 19.628
	 9.906 19.698
	 9.928 19.767
	 9.964 19.835
	10.013 19.901
	10.076 19.966
	10.153 20.030
	10.243 20.093
	10.348 20.154
	10.466 20.214
	10.590 20.272
	10.711 20.326
	10.830 20.377
	10.946 20.424
	11.060 20.467
	11.171 20.507
	11.280 20.543
	11.387 20.575
	11.491 20.604
	11.592 20.630
	11.691 20.651
	11.787 20.670
	11.881 20.684
	11.973 20.695
	12.062 20.703
	12.149 20.706
	12.234 20.706
	12.319 20.703
	12.404 20.696
	12.489 20.685
	12.574 20.671
	12.658 20.653
	12.743 20.632
	12.827 20.607
	12.912 20.578
	12.996 20.546
	13.080 20.511
	13.164 20.471
	13.247 20.428
	13.331 20.382
	13.415 20.331
	13.498 20.278
	13.579 20.222
	13.655 20.168
	13.726 20.114
	13.793 20.061
	13.911 19.958
	14.010 19.857
	14.090 19.759
	14.150 19.664
	14.191 19.573
	14.212 19.484
	 /
\plot 14.212 19.484 14.260 19.135 /
\linethickness= 0.500pt
\setplotsymbol ({\thinlinefont .})
%
%
%
\plot	11.117 18.506 11.145 18.168
 	11.180 18.042
	11.256 18.000
	11.373 18.042
	11.447 18.095
	11.531 18.168
	11.615 18.248
	11.690 18.316
	11.812 18.423
	11.896 18.487
	11.943 18.510
	 /
\plot 11.943 18.510 11.999 18.514 /
\linethickness= 0.500pt
\setplotsymbol ({\thinlinefont .})
%
%
%
\plot	 9.682 15.710  9.745 15.651
 	 9.851 15.648
	 9.936 15.689
	10.043 15.758
	10.147 15.816
	10.226 15.827
	10.308 15.701
	 /
\plot 10.308 15.701 10.340 15.479 /
\linethickness= 0.500pt
\setplotsymbol ({\thinlinefont .})
%
%
%
\plot	10.348 15.306 10.408 15.107
 	10.507 15.033
	10.586 15.090
	10.685 15.209
	10.737 15.280
	10.784 15.343
	10.863 15.441
	10.962 15.531
	 /
\plot 10.962 15.531 11.021 15.551 /
\linethickness= 0.500pt
\setplotsymbol ({\thinlinefont .})
%
%
%
\plot	 9.665 12.103  9.740 12.067
 	 9.851 12.105
	 9.934 12.179
	10.034 12.289
	10.137 12.386
	10.228 12.416
	10.306 12.381
	10.372 12.281
	 /
\plot 10.372 12.281 10.490 12.014 /
\linethickness= 0.500pt
\setplotsymbol ({\thinlinefont .})
%
%
%
\plot	10.530 11.824 10.665 11.570
 	10.735 11.484
	10.812 11.479
	10.896 11.557
	10.941 11.627
	10.987 11.717
	11.031 11.811
	11.072 11.894
	11.110 11.965
	11.143 12.024
	11.237 12.142
	 /
\plot 11.237 12.142 11.301 12.167 /
\linethickness= 0.500pt
\setplotsymbol ({\thinlinefont .})
%
%
%
\plot	10.181  8.649 10.272  8.585
 	10.328  8.569
	10.405  8.583
	10.503  8.629
	10.621  8.705
	10.738  8.768
	10.828  8.776
	10.893  8.729
	10.932  8.625
	 /
\plot 10.932  8.625 10.983  8.363 /
\linethickness= 0.500pt
\setplotsymbol ({\thinlinefont .})
%
%
%
\plot	10.998  8.166 11.062  7.892
 	11.172  7.775
	11.262  7.834
	11.315  7.893
	11.375  7.971
	11.434  8.055
	11.487  8.127
	11.574  8.241
	11.672  8.344
	 /
\plot 11.672  8.344 11.720  8.363 /
\linethickness= 0.500pt
\setplotsymbol ({\thinlinefont .})
%
%
%
\plot	10.340  5.838 10.392  5.806
 	10.498  5.829
	10.592  5.881
	10.714  5.961
	10.833  6.028
	10.921  6.043
	11.002  5.917
	 /
\plot 11.002  5.917 11.021  5.687 /
\linethickness= 0.500pt
\setplotsymbol ({\thinlinefont .})
%
%
%
\plot	11.030  5.529 11.057  5.259
 	11.088  5.157
	11.154  5.122
	11.253  5.153
	11.315  5.193
	11.387  5.250
	11.458  5.311
	11.522  5.364
	11.626  5.444
	11.740  5.500
	 /
\plot 11.740  5.500 11.792  5.489 /
\linethickness= 0.500pt
\setplotsymbol ({\thinlinefont .})
%
%
%
\plot	10.530  3.194 10.614  3.158
 	10.670  3.155
	10.753  3.183
	10.864  3.240
	10.930  3.281
	11.002  3.328
	11.074  3.373
	11.137  3.404
	11.237  3.424
	11.332  3.297
	 /
\plot 11.332  3.297 11.356  3.059 /
\linethickness= 0.500pt
\setplotsymbol ({\thinlinefont .})
%
%
%
\plot	11.364  2.877 11.396  2.611
 	11.428  2.512
	11.492  2.479
	11.588  2.514
	11.649  2.556
	11.717  2.615
	11.785  2.678
	11.847  2.733
	11.949  2.816
	12.066  2.881
	 /
\plot 12.066  2.881 12.126  2.877 /
\linethickness= 0.500pt
\setplotsymbol ({\thinlinefont .})
%
%
%
\plot	13.625  2.963 13.733  2.988
 	13.796  2.984
	13.877  2.947
	13.977  2.878
	14.095  2.777
	14.207  2.691
	14.291  2.668
	14.346  2.709
	14.372  2.812
	 /
\plot 14.372  2.812 14.395  3.082 /
\linethickness= 0.500pt
\setplotsymbol ({\thinlinefont .})
%
%
%
\plot	14.387  3.298 14.407  3.571
 	14.434  3.683
	14.496  3.745
	14.591  3.758
	14.652  3.746
	14.721  3.721
	14.791  3.693
	14.852  3.670
	14.950  3.641
	15.047  3.646
	 /
\plot 15.047  3.646 15.079  3.694 /
\linethickness= 0.500pt
\setplotsymbol ({\thinlinefont .})
%
%
%
\plot	14.254  5.662 14.150  5.639
 	14.032  5.651
	13.963  5.684
	13.887  5.734
	13.831  5.797
	13.819  5.867
	13.852  5.943
	13.931  6.027
	14.026  6.090
	14.108  6.104
	14.179  6.068
	14.237  5.984
	 /
\plot 14.237  5.984 14.340  5.766 /
\linethickness= 0.500pt
\setplotsymbol ({\thinlinefont .})
%
%
%
\plot	14.395  5.599 14.455  5.420
 	14.459  5.335
	14.414  5.260
	14.319  5.194
	14.252  5.165
	14.173  5.138
	14.090  5.118
	14.010  5.108
	13.934  5.109
	13.862  5.120
	13.794  5.141
	13.729  5.174
	13.610  5.269
	13.508  5.371
	13.425  5.442
	13.313  5.492
	 /
\plot 13.313  5.492 13.238  5.480 /
%
%
\put{\SetFigFont{12}{14.4}{rm}\bullet} <0pt, 2pt> at 12.764 26.162
%
%
\put{\SetFigFont{12}{14.4}{rm}\bullet} at 12.435 16.745
%
%
\put{\SetFigFont{12}{14.4}{rm}X} [lB] <-2pt, 2pt> at 10.141 15.875
%
%
\put{\SetFigFont{12}{14.4}{rm}Q} [lB] <-2pt, 0pt> at 12.300 17.145
%
%
\put{\SetFigFont{12}{14.4}{rm}{H_1}} [lB] <0pt, -2pt> at 14.459 21.082
%
%
\put{\SetFigFont{12}{14.4}{rm}\bullet}  <0pt, 4pt> at 12.594  0.366
%
%
\put{\SetFigFont{12}{14.4}{rm}{H_2}} [lB] <-2pt, 0pt> at 14.893  5.357
%
%
\put{\SetFigFont{12}{14.4}{rm}B} [lB] <-2pt, 0pt> at 12.321 14.696
\linethickness=0pt
\putrectangle corners at  8.528 26.187 and 16.383  0.360
\endpicture
\endfig

In figure \PVect\ we have sliced the sphere-with-cross-cap by
planes which are roughly horizontal but tilted a little into
general position and then lifted to an immersion in $\re^4$.
The vector information is given using the conventions defined 
in the left-hand pictures in figure \Whit.  Ignoring
the vectors for the time being, the sequence starts with a
0--handle producing a small circle which passes through 
the top Whitney umbrella at $H_1$ (think of a plane tilted to
the left slicing the sphere-with-cross-cap).  The double point
occurs at $Q$ and then a 1--handle (a bridge) at $B$ and the
bottom Whitney umbrella at $H_2$ (think of the slicing plane titled 
backwards here).  Then the resulting circle shrinks to a
2--handle at the bottom.  

We now explain the vector field.  At the top it is up (towards
the eye) and then two opposite twist fields appear (think of a
twist on a `U' shaped curve).  The right-hand twist field passes
through $H_1$ to form a twist in the curve.  The local sequence
here is the same as figure \Whit\ read upwards.  The left hand
twist field continues downwards, but notice that at $X$ the
slope (indicated by the arrows) changes.  The field near $X$ is
the same as that used to create the figure \mtwist\ sequence
described earlier.  Another figure \Whit\
sequence (reflected this time) changes the twist resulting from 
the bridge into an opposite twist field and the two twist fields
cancel.  We now draw the result of the compression.

\fig{\figkey\PImm: immersion of $P^2$ in $\re^3$}
\beginpicture
\setcoordinatesystem units <0.80000cm,0.80000cm>
\unitlength=1.00000cm
\linethickness=1pt
\setplotsymbol ({\makebox(0,0)[l]{\tencirc\symbol{'160}}})
\setshadesymbol ({\thinlinefont .})
\setlinear
%
%
\linethickness= 0.500pt
\setplotsymbol ({\thinlinefont .})
\circulararc 95.709 degrees from  4.047 25.559 center at  4.938 26.410
%
%
\linethickness= 0.500pt
\setplotsymbol ({\thinlinefont .})
\circulararc 95.940 degrees from  5.906 25.590 center at  6.725 26.373
%
%
\linethickness= 0.500pt
\setplotsymbol ({\thinlinefont .})
\circulararc 74.957 degrees from  2.904 24.479 center at  4.428 26.401
%
%
\linethickness= 0.500pt
\setplotsymbol ({\thinlinefont .})
\circulararc 80.607 degrees from  5.920 24.431 center at  7.273 26.197
%
%
\linethickness= 0.500pt
\setplotsymbol ({\thinlinefont .})
\circulararc 91.749 degrees from  2.366 23.099 center at  4.166 24.736
%
%
\linethickness= 0.500pt
\setplotsymbol ({\thinlinefont .})
\circulararc 92.846 degrees from  5.889 23.004 center at  7.529 24.730
%
%
\linethickness= 0.500pt
\setplotsymbol ({\thinlinefont .})
\circulararc 79.470 degrees from  2.381 21.876 center at  4.192 23.938
%
%
\linethickness= 0.500pt
\setplotsymbol ({\thinlinefont .})
\circulararc 81.254 degrees from  5.874 21.795 center at  7.547 23.878
%
%
\linethickness= 0.500pt
\setplotsymbol ({\thinlinefont .})
\circulararc 39.867 degrees from  4.064 19.511 center at  5.841 24.884
%
%
\linethickness= 0.500pt
\setplotsymbol ({\thinlinefont .})
\setdots < 0.0953cm>
\circulararc 52.980 degrees from  5.857 21.795 center at  4.079 18.350
%
%
\linethickness= 0.500pt
\setplotsymbol ({\thinlinefont .})
\circulararc 44.251 degrees from  9.351 21.891 center at  7.755 17.577
%
%
\linethickness= 0.500pt
\setplotsymbol ({\thinlinefont .})
\circulararc 60.428 degrees from  5.874 23.019 center at  4.094 20.059
%
%
\linethickness= 0.500pt
\setplotsymbol ({\thinlinefont .})
\circulararc 67.087 degrees from  9.351 23.194 center at  7.740 20.503
%
%
\linethickness= 0.500pt
\setplotsymbol ({\thinlinefont .})
\circulararc 78.784 degrees from  5.906 24.448 center at  4.375 22.662
%
%
\linethickness= 0.500pt
\setplotsymbol ({\thinlinefont .})
\circulararc 91.124 degrees from  8.795 24.575 center at  7.412 23.078
%
%
\linethickness= 0.500pt
\setplotsymbol ({\thinlinefont .})
\circulararc 67.637 degrees from  5.889 25.654 center at  5.047 24.295
%
%
\linethickness= 0.500pt
\setplotsymbol ({\thinlinefont .})
\circulararc 56.090 degrees from  7.588 25.639 center at  6.768 24.028
%
%
\linethickness= 0.500pt
\setplotsymbol ({\thinlinefont .})
\circulararc 41.113 degrees from  7.889 19.607 center at  6.103 14.459
%
%
\linethickness= 0.500pt
\setplotsymbol ({\thinlinefont .})
\setsolid
\ellipticalarc axes ratio  3.539:3.539  360 degrees 
	from  9.413 22.543 center at  5.874 22.543
%
%
\linethickness= 0.500pt
\setplotsymbol ({\thinlinefont .})
\ellipticalarc axes ratio  0.243:0.078  360 degrees 
	from 12.234 26.242 center at 11.991 26.242
%
%
\linethickness= 0.500pt
\setplotsymbol ({\thinlinefont .})
\ellipticalarc axes ratio  1.270:0.222  360 degrees 
	from 13.346 19.008 center at 12.076 19.008
%
%
\linethickness= 0.500pt
\setplotsymbol ({\thinlinefont .})
\ellipticalarc axes ratio  0.296:0.127  360 degrees 
	from 12.287 18.510 center at 11.991 18.510
%
%
\linethickness= 0.500pt
\setplotsymbol ({\thinlinefont .})
\ellipticalarc axes ratio  0.150:0.150  360 degrees 
	from  5.929 15.113 center at  5.779 15.113
%
%
\linethickness= 0.500pt
\setplotsymbol ({\thinlinefont .})
\putrule from  5.874 26.099 to  5.874 21.780
%
%
\linethickness=1pt
\setplotsymbol ({\makebox(0,0)[l]{\tencirc\symbol{'160}}})
\putrule from  3.698 23.063 to  4.443 23.063
\plot  4.443 23.063  4.445 22.475 /
\putrule from  4.445 22.475 to  3.698 22.475
\plot  3.698 22.475  3.700 23.063 /
%
%
\linethickness= 0.500pt
\setplotsymbol ({\thinlinefont .})
\setdots < 0.051cm>
\plot 11.991 21.548 12.319 21.516 /
%
%
\linethickness= 0.500pt
\setplotsymbol ({\thinlinefont .})
\setsolid
\putrule from 10.725 23.133 to 10.725 22.689
\putrule from 10.725 22.689 to 11.254 22.689
\putrule from 11.254 22.689 to 11.254 23.133
\putrule from 11.254 23.133 to 10.725 23.133
%
%
\linethickness= 0.500pt
\setplotsymbol ({\thinlinefont .})
\setdots < 0.051cm>
\plot  5.726 16.912  6.054 16.806 /
%
%
\linethickness= 0.500pt
\setplotsymbol ({\thinlinefont .})
\setsolid
\putrule from  2.434 16.703 to  2.434 16.258
\putrule from  2.434 16.258 to  2.963 16.258
\putrule from  2.963 16.258 to  2.963 16.703
\putrule from  2.963 16.703 to  2.434 16.703
%
%
\linethickness= 0.500pt
\setplotsymbol ({\thinlinefont .})
\putrule from 11.007 19.748 to 11.007 19.304
\putrule from 11.007 19.304 to 11.536 19.304
\putrule from 11.536 19.304 to 11.536 19.748
\putrule from 11.536 19.748 to 11.007 19.748
%
%
\linethickness= 0.500pt
\setplotsymbol ({\thinlinefont .})
\putrule from  9.576 16.654 to  9.576 16.209
\putrule from  9.576 16.209 to 10.105 16.209
\putrule from 10.105 16.209 to 10.105 16.654
\putrule from 10.105 16.654 to  9.576 16.654
\linethickness= 0.500pt
\setplotsymbol ({\thinlinefont .})
%
%
%
\plot	 2.587 21.209  3.412 20.749
 	 3.515 20.695
	 3.619 20.649
	 3.723 20.610
	 3.827 20.579
	 3.931 20.556
	 4.036 20.541
	 4.140 20.533
	 4.246 20.532
	 4.351 20.540
	 4.456 20.555
	 4.562 20.577
	 4.668 20.608
	 4.774 20.645
	 4.881 20.691
	 4.988 20.744
	 5.095 20.805
	 5.202 20.865
	 5.309 20.919
	 5.417 20.964
	 5.524 21.002
	 5.632 21.032
	 5.740 21.055
	 5.848 21.070
	 5.956 21.077
	 6.064 21.077
	 6.173 21.069
	 6.281 21.053
	 6.390 21.030
	 6.498 20.999
	 6.607 20.961
	 6.716 20.915
	 6.825 20.861
	 6.934 20.807
	 7.041 20.759
	 7.148 20.718
	 7.253 20.684
	 7.358 20.657
	 7.461 20.636
	 7.563 20.622
	 7.665 20.615
	 7.765 20.614
	 7.864 20.620
	 7.962 20.633
	 8.059 20.652
	 8.155 20.679
	 8.250 20.711
	 8.344 20.751
	 8.437 20.797
	 /
\plot  8.437 20.797  9.176 21.194 /
\linethickness= 0.500pt
\setplotsymbol ({\thinlinefont .})
%
%
%
\plot	 2.999 20.447  3.905 20.043
 	 4.017 19.995
	 4.127 19.951
	 4.235 19.913
	 4.340 19.879
	 4.443 19.850
	 4.543 19.825
	 4.641 19.806
	 4.737 19.791
	 4.831 19.781
	 4.922 19.776
	 5.011 19.775
	 5.098 19.780
	 5.182 19.789
	 5.265 19.803
	 5.344 19.821
	 5.422 19.845
	 5.498 19.869
	 5.575 19.890
	 5.651 19.909
	 5.727 19.924
	 5.804 19.937
	 5.880 19.946
	 5.956 19.952
	 6.032 19.956
	 6.109 19.956
	 6.185 19.954
	 6.261 19.948
	 6.338 19.940
	 6.414 19.929
	 6.490 19.914
	 6.567 19.897
	 6.643 19.877
	 6.720 19.857
	 6.799 19.842
	 6.879 19.833
	 6.960 19.827
	 7.043 19.827
	 7.127 19.831
	 7.213 19.841
	 7.300 19.855
	 7.388 19.873
	 7.478 19.897
	 7.570 19.925
	 7.663 19.958
	 7.757 19.996
	 7.853 20.038
	 7.950 20.086
	 8.049 20.138
	 /
\plot  8.049 20.138  8.843 20.574 /
\linethickness= 0.500pt
\setplotsymbol ({\thinlinefont .})
\setdots < 0.0953cm>
%
%
%
\plot	 2.555 21.258  3.484 21.551
 	 3.600 21.585
	 3.714 21.615
	 3.827 21.640
	 3.940 21.660
	 4.051 21.676
	 4.161 21.687
	 4.270 21.694
	 4.377 21.695
	 4.484 21.692
	 4.589 21.685
	 4.694 21.673
	 4.797 21.656
	 4.899 21.635
	 5.000 21.608
	 5.100 21.578
	 5.199 21.542
	 5.297 21.507
	 5.397 21.477
	 5.498 21.452
	 5.601 21.431
	 5.704 21.416
	 5.809 21.405
	 5.915 21.400
	 6.022 21.400
	 6.130 21.404
	 6.239 21.413
	 6.350 21.428
	 6.462 21.447
	 6.575 21.471
	 6.689 21.501
	 6.804 21.535
	 6.920 21.574
	 7.036 21.613
	 7.150 21.647
	 7.261 21.675
	 7.370 21.697
	 7.477 21.714
	 7.582 21.726
	 7.684 21.732
	 7.784 21.733
	 7.881 21.728
	 7.977 21.718
	 8.070 21.702
	 8.160 21.681
	 8.249 21.655
	 8.335 21.623
	 8.419 21.585
	 8.501 21.542
	 /
\plot  8.501 21.542  9.144 21.177 /
\linethickness= 0.500pt
\setplotsymbol ({\thinlinefont .})
%
%
%
\plot	 2.999 20.447  3.722 20.701
 	 3.813 20.731
	 3.904 20.757
	 3.996 20.780
	 4.088 20.800
	 4.182 20.816
	 4.275 20.828
	 4.370 20.837
	 4.465 20.842
	 4.560 20.844
	 4.656 20.842
	 4.753 20.837
	 4.850 20.828
	 4.948 20.815
	 5.047 20.800
	 5.146 20.780
	 5.246 20.757
	 5.346 20.734
	 5.446 20.714
	 5.546 20.697
	 5.645 20.684
	 5.745 20.674
	 5.844 20.667
	 5.943 20.663
	 6.041 20.662
	 6.140 20.665
	 6.239 20.671
	 6.337 20.680
	 6.435 20.692
	 6.533 20.707
	 6.630 20.726
	 6.728 20.748
	 6.825 20.773
	 6.922 20.798
	 7.017 20.819
	 7.112 20.838
	 7.205 20.852
	 7.297 20.863
	 7.387 20.871
	 7.477 20.875
	 7.565 20.876
	 7.653 20.873
	 7.739 20.867
	 7.824 20.857
	 7.907 20.844
	 7.990 20.828
	 8.071 20.808
	 8.152 20.784
	 8.231 20.757
	 /
\plot  8.231 20.757  8.858 20.527 /
\linethickness=2pt
\setplotsymbol ({\makebox(0,0)[l]{\tencirc\symbol{'161}}})
\setsolid
%
%
%
\plot	 5.889 26.020  5.754 26.043
 	 5.666 26.034
	 5.606 26.013
	 5.536 25.982
	 5.456 25.939
	 5.365 25.887
	 5.264 25.823
	 5.152 25.749
	 5.038 25.662
	 4.929 25.560
	 4.826 25.443
	 4.776 25.379
	 4.728 25.311
	 4.681 25.239
	 4.635 25.163
	 4.591 25.084
	 4.548 25.001
	 4.506 24.914
	 4.466 24.823
	 4.427 24.729
	 4.390 24.631
	 4.354 24.532
	 4.321 24.437
	 4.289 24.345
	 4.259 24.257
	 4.232 24.172
	 4.206 24.091
	 4.183 24.013
	 4.161 23.938
	 4.142 23.867
	 4.125 23.799
	 4.096 23.674
	 4.076 23.562
	 4.064 23.464
	 /
\plot  4.064 23.464  4.032 23.099 /
\linethickness=2pt
\setplotsymbol ({\makebox(0,0)[l]{\tencirc\symbol{'161}}})
\setdashes < 0.17180cm>
%
%
%
\plot	 5.874 21.732  5.858 21.574
 	 5.829 21.491
	 5.759 21.400
	 5.647 21.301
	 5.575 21.249
	 5.493 21.194
	 5.404 21.146
	 5.312 21.111
	 5.216 21.090
	 5.118 21.083
	 5.016 21.090
	 4.911 21.111
	 4.802 21.146
	 4.691 21.194
	 4.583 21.254
	 4.487 21.323
	 4.403 21.401
	 4.331 21.488
	 4.270 21.584
	 4.221 21.688
	 4.184 21.802
	 4.158 21.924
	 /
\plot  4.158 21.924  4.079 22.432 /
\linethickness= 0.500pt
\setplotsymbol ({\thinlinefont .})
\setsolid
%
%
%
\plot	12.634 23.789 
 	12.761 23.732
	12.885 23.693
	13.006 23.672
	13.124 23.668
	13.239 23.681
	13.352 23.713
	13.461 23.762
	13.568 23.828
	13.659 23.902
	13.723 23.973
	13.760 24.040
	13.769 24.103
	13.704 24.221
	13.630 24.274
	13.529 24.325
	13.404 24.363
	13.335 24.374
	13.261 24.379
	13.182 24.379
	13.098 24.374
	13.010 24.364
	12.916 24.348
	12.819 24.326
	12.716 24.300
	12.609 24.268
	12.496 24.230
	12.380 24.188
	12.258 24.140
	12.132 24.086
	12.067 24.057
	12.000 24.027
	11.935 23.997
	11.871 23.967
	11.749 23.911
	11.636 23.857
	11.530 23.807
	11.433 23.759
	11.343 23.715
	11.261 23.674
	11.188 23.636
	11.064 23.570
	10.973 23.517
	10.885 23.447
	10.867 23.370
	10.908 23.312
	10.980 23.289
	11.055 23.316
	11.097 23.373
	11.071 23.440
	11.023 23.489
	10.940 23.562
	10.821 23.661
	10.748 23.720
	10.666 23.785
	10.592 23.853
	10.544 23.922
	10.521 23.990
	10.523 24.058
	10.551 24.127
	10.604 24.195
	10.682 24.264
	10.785 24.332
	10.901 24.392
	11.016 24.435
	11.130 24.460
	11.243 24.468
	11.355 24.459
	11.466 24.433
	11.577 24.390
	11.686 24.329
	11.796 24.260
	11.909 24.191
	12.024 24.123
	12.141 24.055
	12.261 23.988
	12.383 23.921
	12.508 23.855
	12.634 23.789
	/
\linethickness= 0.500pt
\setplotsymbol ({\thinlinefont .})
%
%
%
\plot	11.018 24.987 
 	10.940 25.053
	10.936 25.118
	11.006 25.184
	11.069 25.216
	11.150 25.249
	11.244 25.274
	11.347 25.284
	11.459 25.279
	11.580 25.260
	11.644 25.244
	11.709 25.225
	11.777 25.203
	11.847 25.176
	11.920 25.146
	11.994 25.112
	12.071 25.074
	12.150 25.033
	12.229 24.992
	12.305 24.955
	12.380 24.922
	12.453 24.894
	12.523 24.870
	12.592 24.850
	12.658 24.834
	12.722 24.823
	12.845 24.813
	12.959 24.820
	13.065 24.844
	13.162 24.885
	13.245 24.933
	13.308 24.980
	13.369 25.070
	13.347 25.153
	13.242 25.230
	13.161 25.260
	13.066 25.277
	12.956 25.281
	12.831 25.272
	12.764 25.262
	12.693 25.250
	12.618 25.234
	12.539 25.215
	12.457 25.192
	12.371 25.166
	12.282 25.137
	12.189 25.105
	12.096 25.072
	12.009 25.040
	11.926 25.010
	11.848 24.981
	11.775 24.955
	11.706 24.929
	11.584 24.884
	11.481 24.845
	11.398 24.812
	11.288 24.765
	11.213 24.698
	11.269 24.634
	11.385 24.608
	11.486 24.657
	11.515 24.736
	11.414 24.800
	11.327 24.833
	11.232 24.875
	11.129 24.926
	11.018 24.987
	/
\linethickness= 0.500pt
\setplotsymbol ({\thinlinefont .})
%
%
%
\plot	11.310 25.706 
 	11.402 25.719
	11.497 25.702
	11.596 25.655
	11.699 25.579
	11.800 25.509
	11.892 25.481
	11.977 25.497
	12.053 25.555
	12.139 25.622
	12.253 25.664
	12.320 25.676
	12.393 25.681
	12.474 25.681
	12.561 25.674
	12.642 25.663
	12.701 25.649
	12.757 25.612
	12.730 25.564
	12.620 25.504
	12.545 25.476
	12.474 25.456
	12.406 25.447
	12.342 25.446
	12.224 25.472
	12.120 25.535
	12.022 25.596
	11.921 25.616
	11.818 25.594
	11.713 25.531
	11.608 25.462
	11.506 25.422
	11.407 25.410
	11.311 25.428
	11.229 25.460
	11.171 25.495
	11.127 25.570
	11.175 25.644
	11.232 25.676
	11.310 25.706
	/
\linethickness= 0.500pt
\setplotsymbol ({\thinlinefont .})
%
%
%
\plot	11.610 26.029 11.521 25.989
 	11.485 25.939
	11.556 25.867
	11.621 25.833
	11.694 25.815
	11.773 25.813
	11.860 25.827
	11.952 25.842
	12.048 25.843
	12.149 25.831
	12.253 25.806
	12.352 25.786
	12.436 25.789
	12.555 25.862
	12.567 25.968
	12.517 26.010
	12.428 26.046
	12.326 26.068
	12.235 26.073
	12.155 26.060
	12.086 26.030
	12.022 25.998
	11.957 25.982
	11.891 25.982
	11.825 25.998
	11.710 26.028
	11.634 26.021
	 /
\plot 11.634 26.021 11.578 25.997 /
\linethickness= 0.500pt
\setplotsymbol ({\thinlinefont .})
%
%
%
\plot	12.664 21.920 
 	12.791 21.863
	12.914 21.824
	13.035 21.803
	13.153 21.799
	13.269 21.812
	13.381 21.844
	13.491 21.893
	13.597 21.959
	13.689 22.033
	13.753 22.103
	13.789 22.171
	13.798 22.234
	13.733 22.352
	13.660 22.405
	13.558 22.456
	13.435 22.493
	13.367 22.504
	13.295 22.508
	13.219 22.507
	13.139 22.501
	13.055 22.489
	12.966 22.471
	12.874 22.448
	12.777 22.419
	12.676 22.385
	12.571 22.345
	12.462 22.299
	12.349 22.248
	12.232 22.191
	12.111 22.129
	11.990 22.066
	11.877 22.009
	11.771 21.958
	11.671 21.912
	11.578 21.872
	11.492 21.838
	11.412 21.808
	11.340 21.785
	11.214 21.755
	11.116 21.747
	11.001 21.797
	10.954 21.885
	10.935 21.956
	10.980 22.051
	11.069 22.073
	11.134 22.009
	11.142 21.898
	11.060 21.776
	10.986 21.738
	10.892 21.748
	10.777 21.809
	10.712 21.858
	10.642 21.919
	10.581 21.987
	10.543 22.054
	10.530 22.122
	10.539 22.190
	10.573 22.258
	10.630 22.326
	10.711 22.395
	10.815 22.463
	10.931 22.523
	11.046 22.566
	11.160 22.591
	11.273 22.599
	11.385 22.590
	11.496 22.564
	11.606 22.521
	11.716 22.460
	11.826 22.391
	11.939 22.322
	12.053 22.254
	12.171 22.186
	12.291 22.119
	12.413 22.052
	12.537 21.986
	12.664 21.920
	/
\linethickness= 0.500pt
\setplotsymbol ({\thinlinefont .})
%
%
%
\plot	12.649 20.832 
 	12.776 20.775
	12.900 20.736
	13.021 20.715
	13.139 20.711
	13.254 20.724
	13.366 20.756
	13.476 20.805
	13.583 20.871
	13.674 20.945
	13.738 21.016
	13.775 21.083
	13.783 21.146
	13.719 21.264
	13.645 21.317
	13.543 21.368
	13.429 21.413
	13.315 21.453
	13.202 21.486
	13.089 21.513
	12.978 21.534
	12.867 21.549
	12.757 21.558
	12.648 21.560
	12.548 21.559
	12.464 21.554
	12.348 21.538
	12.319 21.474
	12.324 21.409
	12.234 21.301
	12.153 21.231
	12.048 21.150
	11.987 21.105
	11.919 21.057
	11.846 21.007
	11.767 20.954
	11.686 20.902
	11.610 20.855
	11.539 20.813
	11.471 20.776
	11.348 20.716
	11.242 20.675
	11.153 20.655
	11.081 20.653
	10.987 20.709
	10.940 20.797
	10.920 20.868
	10.965 20.963
	11.054 20.985
	11.119 20.921
	11.127 20.810
	11.045 20.688
	10.971 20.650
	10.877 20.661
	10.762 20.721
	10.697 20.770
	10.627 20.831
	10.566 20.899
	10.528 20.966
	10.515 21.034
	10.525 21.102
	10.558 21.170
	10.615 21.238
	10.696 21.307
	10.800 21.375
	10.915 21.439
	11.029 21.493
	11.140 21.538
	11.249 21.573
	11.356 21.599
	11.461 21.615
	11.564 21.622
	11.665 21.620
	11.757 21.611
	11.835 21.602
	11.945 21.578
	11.986 21.511
	11.955 21.409
	11.976 21.342
	12.022 21.263
	12.104 21.173
	12.164 21.124
	12.236 21.071
	12.321 21.016
	12.418 20.958
	12.527 20.896
	12.649 20.832
	/
\linethickness= 0.500pt
\setplotsymbol ({\thinlinefont .})
%
%
%
\plot	11.393 20.611 
 	11.495 20.628
	11.598 20.641
	11.700 20.653
	11.802 20.661
	11.904 20.667
	12.005 20.670
	12.107 20.670
	12.208 20.668
	12.309 20.663
	12.410 20.655
	12.511 20.645
	12.611 20.632
	12.712 20.616
	12.812 20.598
	12.912 20.577
	13.012 20.553
	13.109 20.528
	13.200 20.504
	13.285 20.480
	13.363 20.457
	13.436 20.434
	13.503 20.412
	13.618 20.369
	13.709 20.328
	13.776 20.290
	13.838 20.219
	13.823 20.149
	13.754 20.070
	13.630 19.982
	13.548 19.935
	13.451 19.886
	13.348 19.840
	13.242 19.804
	13.136 19.776
	13.028 19.758
	12.919 19.748
	12.808 19.748
	12.697 19.757
	12.584 19.775
	12.475 19.798
	12.376 19.821
	12.286 19.845
	12.207 19.870
	12.137 19.896
	12.077 19.922
	11.986 19.976
	11.904 20.085
	11.943 20.188
	12.053 20.243
	12.115 20.237
	12.181 20.209
	12.229 20.110
	12.185 20.045
	12.097 19.971
	11.971 19.904
	11.898 19.880
	11.818 19.862
	11.730 19.851
	11.635 19.846
	11.534 19.847
	11.425 19.854
	11.319 19.866
	11.228 19.881
	11.151 19.898
	11.087 19.918
	10.975 20.024
	11.005 20.132
	11.097 20.204
	11.195 20.222
	11.245 20.172
	11.188 20.078
	11.097 20.024
	11.036 19.995
	10.964 19.965
	10.892 19.942
	10.827 19.932
	10.724 19.952
	10.653 20.027
	10.630 20.085
	10.615 20.156
	10.618 20.231
	10.648 20.301
	10.705 20.366
	10.788 20.426
	10.899 20.480
	10.965 20.505
	11.037 20.529
	11.116 20.552
	11.201 20.573
	11.294 20.593
	11.393 20.611
	/
\linethickness= 0.500pt
\setplotsymbol ({\thinlinefont .})
%
%
%
\plot	 4.085 17.018  4.503 16.838
 	 4.608 16.797
	 4.715 16.766
	 4.822 16.743
	 4.931 16.728
	 5.040 16.723
	 5.151 16.726
	 5.263 16.738
	 5.376 16.759
	 5.480 16.784
	 5.565 16.810
	 5.674 16.862
	 5.657 16.970
	 5.563 17.013
	 5.457 17.029
	 5.338 17.018
	 5.274 17.003
	 5.207 16.981
	 5.108 16.922
	 5.085 16.847
	 5.140 16.756
	 5.195 16.704
	 5.271 16.648
	 5.356 16.593
	 5.444 16.547
	 5.535 16.509
	 5.628 16.478
	 5.723 16.456
	 5.820 16.441
	 5.920 16.435
	 6.022 16.436
	 6.120 16.443
	 6.209 16.455
	 6.288 16.470
	 6.357 16.490
	 6.465 16.542
	 6.535 16.611
	 6.554 16.681
	 6.509 16.740
	 6.400 16.787
	 6.322 16.806
	 6.228 16.822
	 6.138 16.832
	 6.071 16.830
	 6.006 16.797
	 6.034 16.722
	 6.154 16.605
	 6.238 16.540
	 6.324 16.484
	 6.410 16.435
	 6.498 16.394
	 6.587 16.361
	 6.678 16.335
	 6.770 16.318
	 6.863 16.309
	 /
\plot  6.863 16.309  7.239 16.288 /
\linethickness= 0.500pt
\setplotsymbol ({\thinlinefont .})
%
%
%
\plot	 5.493 16.145 
 	 5.559 16.153
	 5.625 16.159
	 5.751 16.161
	 5.872 16.152
	 5.986 16.130
	 6.095 16.097
	 6.198 16.053
	 6.296 15.996
	 6.387 15.928
	 6.461 15.856
	 6.507 15.788
	 6.514 15.663
	 6.407 15.555
	 6.311 15.507
	 6.186 15.462
	 6.117 15.443
	 6.050 15.429
	 5.985 15.418
	 5.922 15.412
	 5.800 15.413
	 5.685 15.432
	 5.576 15.468
	 5.475 15.521
	 5.380 15.591
	 5.292 15.679
	 5.221 15.771
	 5.177 15.853
	 5.161 15.926
	 5.173 15.989
	 5.278 16.086
	 5.372 16.120
	 5.493 16.145
	/
\linethickness= 0.500pt
\setplotsymbol ({\thinlinefont .})
%
%
%
\plot	 4.072 16.487  4.877 16.137
 	 4.977 16.095
	 5.077 16.054
	 5.177 16.015
	 5.277 15.979
	 5.376 15.944
	 5.476 15.911
	 5.575 15.880
	 5.673 15.852
	 5.772 15.825
	 5.870 15.800
	 5.968 15.777
	 6.065 15.756
	 6.163 15.738
	 6.260 15.721
	 6.357 15.706
	 6.454 15.693
	 /
\plot  6.454 15.693  7.226 15.598 /
\linethickness= 0.500pt
\setplotsymbol ({\thinlinefont .})
%
%
%
\plot	 4.032 15.505  4.837 15.155
 	 4.937 15.113
	 5.037 15.072
	 5.137 15.033
	 5.237 14.997
	 5.336 14.962
	 5.435 14.929
	 5.534 14.898
	 5.633 14.870
	 5.731 14.843
	 5.830 14.818
	 5.928 14.795
	 6.025 14.774
	 6.123 14.755
	 6.220 14.739
	 6.317 14.724
	 6.413 14.711
	 /
\plot  6.413 14.711  7.186 14.616 /
\linethickness= 0.500pt
\setplotsymbol ({\thinlinefont .})
%
%
%
\plot	 4.032 14.796  4.837 14.446
 	 4.937 14.404
	 5.037 14.363
	 5.137 14.324
	 5.237 14.287
	 5.336 14.253
	 5.435 14.220
	 5.534 14.189
	 5.633 14.160
	 5.731 14.134
	 5.830 14.109
	 5.928 14.086
	 6.025 14.065
	 6.123 14.046
	 6.220 14.030
	 6.317 14.015
	 6.413 14.002
	 /
\plot  6.413 14.002  7.186 13.906 /
\linethickness= 0.500pt
\setplotsymbol ({\thinlinefont .})
%
%
%
\plot	 4.096 17.759  4.466 17.584
 	 4.556 17.544
	 4.639 17.512
	 4.716 17.488
	 4.786 17.472
	 4.909 17.462
	 5.006 17.484
	 5.122 17.578
	 5.141 17.639
	 5.133 17.711
	 5.064 17.801
	 4.942 17.764
	 4.890 17.706
	 4.882 17.633
	 4.917 17.545
	 4.995 17.441
	 5.112 17.341
	 5.182 17.300
	 5.260 17.264
	 5.346 17.234
	 5.440 17.209
	 5.542 17.190
	 5.652 17.177
	 5.760 17.169
	 5.859 17.166
	 5.948 17.167
	 6.027 17.174
	 6.096 17.186
	 6.156 17.203
	 6.244 17.251
	 6.337 17.370
	 6.339 17.494
	 6.275 17.568
	 6.165 17.537
	 6.121 17.416
	 6.166 17.329
	 6.255 17.224
	 6.377 17.123
	 6.447 17.081
	 6.522 17.044
	 6.603 17.014
	 6.689 16.990
	 6.781 16.972
	 6.879 16.960
	 /
\plot  6.879 16.960  7.281 16.923 /
\linethickness= 0.500pt
\setplotsymbol ({\thinlinefont .})
%
%
%
\plot	11.096 17.541 11.556 17.356
 	11.664 17.310
	11.759 17.264
	11.841 17.219
	11.909 17.174
	12.005 17.087
	12.048 17.001
	12.053 16.924
	12.037 16.863
	11.942 16.789
	11.830 16.787
	11.768 16.863
	11.801 16.975
	11.870 17.026
	11.974 17.075
	12.038 17.095
	12.109 17.109
	12.187 17.115
	12.272 17.115
	12.363 17.107
	12.461 17.093
	12.566 17.072
	12.678 17.043
	12.788 17.014
	12.887 16.991
	12.975 16.973
	13.053 16.961
	13.178 16.954
	13.260 16.969
	13.352 17.046
	13.387 17.170
	13.355 17.263
	13.249 17.244
	13.175 17.131
	13.190 17.044
	13.239 16.938
	13.319 16.832
	13.428 16.750
	13.493 16.717
	13.565 16.690
	13.644 16.668
	13.731 16.652
	 /
\plot 13.731 16.652 14.091 16.599 /
\linethickness= 0.500pt
\setplotsymbol ({\thinlinefont .})
%
%
%
\plot	11.096 16.662 11.556 16.477
 	11.664 16.431
	11.759 16.386
	11.841 16.341
	11.909 16.296
	12.005 16.208
	12.048 16.123
	12.053 16.046
	12.037 15.985
	11.942 15.911
	11.830 15.908
	11.768 15.985
	11.801 16.096
	11.870 16.148
	11.974 16.197
	12.097 16.232
	12.223 16.242
	12.287 16.237
	12.351 16.226
	12.416 16.209
	12.482 16.186
	12.547 16.163
	12.611 16.147
	12.735 16.135
	12.852 16.149
	12.964 16.191
	13.052 16.249
	13.101 16.309
	13.110 16.373
	13.080 16.440
	12.986 16.495
	12.907 16.513
	12.806 16.523
	12.683 16.527
	12.614 16.527
	12.539 16.524
	12.459 16.520
	12.373 16.515
	12.282 16.507
	12.186 16.498
	12.093 16.487
	12.012 16.472
	11.889 16.434
	11.794 16.321
	11.822 16.246
	11.900 16.159
	11.958 16.111
	12.028 16.060
	12.111 16.005
	12.207 15.948
	12.308 15.891
	12.406 15.840
	12.503 15.792
	12.597 15.750
	12.689 15.712
	12.779 15.678
	12.866 15.649
	12.952 15.625
	13.035 15.606
	13.116 15.591
	13.195 15.581
	13.271 15.575
	13.346 15.574
	13.418 15.578
	13.488 15.586
	13.556 15.599
	 /
\plot 13.556 15.599 14.091 15.720 /
\linethickness= 0.500pt
\setplotsymbol ({\thinlinefont .})
%
%
%
\plot	11.096 15.720 11.556 15.535
 	11.664 15.489
	11.759 15.444
	11.841 15.399
	11.909 15.354
	12.005 15.266
	12.048 15.181
	12.053 15.104
	12.037 15.043
	11.942 14.969
	11.830 14.966
	11.768 15.043
	11.801 15.154
	11.870 15.206
	11.974 15.255
	12.089 15.292
	12.192 15.306
	12.283 15.299
	12.360 15.271
	12.431 15.242
	12.502 15.235
	12.572 15.250
	12.641 15.287
	12.721 15.391
	12.688 15.514
	12.611 15.567
	12.549 15.583
	12.470 15.592
	12.375 15.595
	12.265 15.591
	12.138 15.580
	12.069 15.571
	11.995 15.562
	11.923 15.549
	11.856 15.533
	11.741 15.490
	11.650 15.433
	11.582 15.361
	11.539 15.274
	11.519 15.173
	11.523 15.057
	11.534 14.994
	11.551 14.927
	11.606 14.801
	11.695 14.699
	11.815 14.620
	11.887 14.589
	11.967 14.564
	12.056 14.545
	12.152 14.532
	12.257 14.525
	12.369 14.524
	12.490 14.528
	12.553 14.533
	12.618 14.539
	12.686 14.546
	12.755 14.555
	12.827 14.566
	12.900 14.577
	 /
\plot 12.900 14.577 14.091 14.779 /
\linethickness= 0.500pt
\setplotsymbol ({\thinlinefont .})
%
%
%
\plot	11.085 14.503 11.905 14.107
 	12.007 14.059
	12.109 14.015
	12.210 13.976
	12.310 13.940
	12.410 13.909
	12.508 13.881
	12.607 13.858
	12.704 13.838
	12.801 13.822
	12.897 13.811
	12.993 13.803
	13.088 13.800
	13.182 13.800
	13.276 13.805
	13.369 13.814
	13.461 13.826
	 /
\plot 13.461 13.826 14.196 13.942 /
%
%
\put{\SetFigFont{12}{14.4}{rm}X} [lB] at  3.829 22.674
%
%
\put{\SetFigFont{12}{14.4}{rm}X} [lB] at 10.751 22.780
%
%
\put{\SetFigFont{12}{14.4}{rm}Y} [lB] at  2.460 16.349
%
%
\put{\SetFigFont{12}{14.4}{rm}Y} [lB] at 11.032 19.395
%
%
\put{\SetFigFont{12}{14.4}{rm}{}where} [lB] at  2.258 17.003
%
%
\put{\SetFigFont{12}{14.4}{rm}=} [lB] at  3.147 16.385
%
%
\put{\SetFigFont{12}{14.4}{rm}Z} [lB] at  9.601 16.300
%
%
\put{\SetFigFont{12}{14.4}{rm}{}and} [lB] at  9.449 16.988
%
%
\put{\SetFigFont{12}{14.4}{rm}{}(see text)} [lB] at  9.102 15.591
%
%
\put{\SetFigFont{12}{14.4}{rm}=} [lB] at 10.355 16.322
\linethickness=0pt
\putrectangle corners at  2.258 26.338 and 14.213 13.750
\endpicture
\endfig

Figure \PImm\ comprises two views of the immersion obtained by
projecting the resulting compressible immersion in $\re^4$ to
an immersion in $\re^3$.  The twist fields have all been replaced
by ripples (as explained in connection with figure \Whit\ above).
This is indicated in the left view by a heavy
black line near the top and the heavy dashed line near the middle.
The heavy line represents a ripple on the {\it outside} of the top
sheet and the dashed line represents a ripple on the {\it inside}
of the top sheet.  The square containing the point $X$ should be 
imagined to be
enlarged and contain a copy of the immersion in $\re^3$ obtained by
projecting figure \mtwist.  Notice that
figure \mtwist\ can be seen to start with a twist to the right
of the line and ends with a twist to the left of the line.  The
projected immersion therefore starts with a ripple on one side
of the sheet and ends with a ripple on the other side.

To the right in figure \PImm\ is a sequence of cross-sections of
the immersions, roughly on the same levels as the left-hand view.
These cross-sections make clear what is happening near the singular
points.  Near the top the ripple turns around to join the `big ripple' 
which is the right half of the cross-cap.  Near the middle it joins 
into the back sheet of the Whitney umbrella.  The small dotted line 
in the fourth section from the bottom at the right is a bridge 
(the 1--handle).  The square marked $X$ should again be imagined
replaced by the sequence in figure \mtwist\ and notice that this
sequence contains a triple point.  The sequence
labelled $Y$ (to be inserted in the square marked $Y$) can be very
simply visualised as an immersion: it is just an immersion in
which a ripple makes a `U' shape.

The final immersion is isotopic to Boy's surface. This is best 
seen by sliding the $X$ sequence down to meet the $Y$ sequence.
The two sequences can then be combined in a single sequence
(with no critical levels) detailed as $Z$ in the figure.  The
immersion now has just three critical levels and can be seen to 
coincide with Philips' projection of Boy's surface [\Ph].

\section{The local proof}    
    
We now prove that the isotopy constructed in the Compression Theorem
can be assumed to be arbitrarily small in the $C^0$ sense.  This
implies that it can be assumed to take place within an arbitrary
neighbourhood of the given embedding.  We can then apply the result to
an immersion by working in an induced regular neighbourhood and this
opens the way to an inductive proof of the `multi-compression theorem'
which states that a number of vector fields can be straightened
simultaneously.
    
What we shall do is to identify a certain submanifold of $M$ which we
call the {\sl downset}.  The closure of the downset is a manifold with
boundary.  We call the boundary the {\sl horizontal set}.  We shall
see that nearly all the straightening process can be made to take
place in an arbitrary neighbourhood of the downset and, by choosing
this neighbourhood to be sufficiently small, the straightening becomes
arbitrarily small.  It will help in understanding the proof to refer
to the pictures in the last section.  In figure \vectp\ the downset is
one point, namely $D$ and note that the initial vector field twists
once around $M$ in an interval centred on $D$.  The straightening
becomes smaller (and tighter) if the initial twist is made to take
place in a smaller interval.  In the surface which determines figure
\mtwist\ the downset is a line and in figure \xtwist\ it is a plane.
In each case it can be seen that the straightening takes place in a
neighbourhood of the downset and can be made small by choosing this
neighborhood sufficiently small.

In figure \Whit\ the downset is a half-open interval (comprising
points like $D$ in the middle of twist fields) with boundary at $H$,
which is the horizontal set in this example.

In the final example (figures \CCap\ to \PImm) the downset again 
comprises a point in the middle of all the twist fields.  It
thus forms an interval in the neighbourhood of which the ripple
in figure \PImm\ is constructed.  It has as boundary the two 
singular points $H_1$ and $H_2$, which form the horizontal set in 
this example.

We start by defining the horizontal set.  This is independent of the
normal vector field and defined for embeddings in $Q\times\re^r$ for
any $r\ge1$.

\sh{The horizontal set}

Suppose given a submanifold $M^m$ in $Q^q\times \re^r$ 
where $m\le q$.
Think of $Q$ as {\it horizontal} and $\re^r$ as {\it vertical}.
Define the {\sl horizontal set} of $M$, denoted $H(M)$,
by $H(M)=\{(x,y)\in M\mid T_{(x,y)}(\{x\}\times\re^r)\subset 
T_{(x,y)}(M)\}$.
Now suppose $Q$ has a metric and give $Q\times\re^r$ the product 
metric, then $H(M)$ is the set of points in $M$ which have 
{\it horizontal\/} normal fibres, ie, (after the obvious identification) 
$H(M)=\{(x,y)\in M\mid \nu_{(x,y)}\subset T_x(Q)\}$ where 
$\nu_{(x,y)}$ denotes the fibre of the normal bundle at $(x,y)$ defined by 
the metric.

\numkey\HSet
\proc{Proposition}After a small isotopy of $M$ in $Q\times\re^r$ the 
horizontal set $H(M)$ can be assumed to be a submanifold of $M$.\rm

The proof (given in the appendix, see proposition A.2) shows more.  
$H(M)$ is given locally by transversality to the appropriate Grassmannian.

\sh{The downset}

After application of the proposition with $r=1$ the corresponding    
horizontal set $H$ is a submanifold of $M$ of codimension $c=q-m+1$.
Let $\psi$ be the unit perpendicular field on $M-H$ 
defined by choosing the $\re$ coordinate to be 
 maximal, so $\psi$ points `upmost' in its normal fibre. 
Suppose $M$ is equipped with a unit perpendicular  
field $\alpha$.  Call $D=\{x\in M\mid \alpha(x)=-\psi(x)\}$ the
{\sl downset} of $M$;  so $D$ is where $\alpha$ points {\it downmost}
in the normal fibre.
   
\numkey\DownMan
\proc{Proposition} After a small isotopy of $\alpha$ we can assume that    
$\alpha$ is transverse to $-\psi$, which implies that $D$ is a 
manifold, and further we can assume that the closure of $D$ is a manifold with 
boundary $H$.\rm   
   
The proof is given in the appendix, see proposition A.3.

\sh{Localisation}

The following considerations explain why the downset is the key
to the straightening process.  Suppose that there is no downset
(and no horizontal set).  Then our perpendicular vector field  $\alpha$ can 
be canonically isotoped to point to the upmost position in each normal
fibre (ie, to the section $\psi$) by isotoping along great circles.  
Once in this upmost position, upwards rotation straightens the field.

More generally, let $\overline D$ denote the closure of the downset 
after proposition \DownMan\ (so that $\overline D$ is a manifold with
boundary $H$) and let $W$ be a tubular neighbourhood of $\overline D$ 
in $M$.  Then $\alpha$ can be canonically isotoped to $\psi$ on 
$\overline{(M-W)}$ and this isotopy extended to $W$ via a collar of
$\d W$ in $W$.  We call this isotopy {\sl localisation} because it
localises the straightening problem in $W$ and we call the resulting
field the {\sl localised field}.  

\rk{Remark}It is possible to give an explicit description of the
localised field in terms of the transversality map which defines
$D$.  We sketch this description (which will not be used).
The field is upmost outside a tubular neighbourhood $V$ of $D$.
A fibre of $V$ comes by transversality from a neighbourhood of
the ``south pole'' in the normal sphere, and by composing with
a standard stretch of this neighbourhood over the sphere, can
be identified with the normal sphere.  The vector field at a point
on the fibre points in the direction thus determined.
\ppar

Now let $\varphi$ be the gradient field on $M$ 
determined by the projection $M\subset Q\times\re\to\re$. Thus for 
$p\in M$ the tangent vector $\varphi_p$, if non-zero, points 
{\it upmost} in the tangent space $T_p(M)$. Note that our 
normal field is grounded if and only if the zeros of $\varphi$ are 
not on $D$.   

The proof of the following lemma is again to be found in the appendix,
see corollary A.6.

\numkey\GPos\proc{Lemma}Suppose that $\alpha$ is perpendicular and
grounded.  Let $U$ be a tubular neighbourhood of $H$.
By a small isotopy of $\alpha$ we may assume that 
$\overline{D-U}$ is in general position with respect to $\varphi$    
in the following sense:  Given $\delta>0$ there    
is a neighbourhood $V$ of $\overline{D-U}$ in $M$ such that each 
component of intersection of an integral curve of $\varphi$ with 
$V$ has length $<\delta$.\rm

We are now ready to prove the main result of this section.
  
\numkey\LCT    
\proc{Local compression theorem}Given $\varepsilon>0$ and the    
hypotheses of the compression theorem, there is an isotopy of     
$M$ to a compressible    
embedding which moves each point a distance  $<\varepsilon$.    
    
\prf  By lemma \Ground\ we can assume that the given vector
field $\alpha$ is perpendicular and grounded.  We shall straighten 
$\alpha$ in two moves.  The first move is a $C^0$--small
move which takes place near the downset $D$, which we shall call the
{\sl local move} and which contains the meat
of the straightening process.  The second move is a 
$C^\infty$--small global move.  Both moves use the modified global 
flow defined in the proof of addendum (i) to \GComp, though in the 
local move, its effect is restricted to a neighbourhood of the
downset.  There is one important observation about the flow near the
horizontal set which we need to make at the outset.   

{\bf Observation}\qua Near the horizontal set the upward 
rotated vector field $\beta$ is 
nearly vertical (since $\alpha$ is nearly horizontal there) and 
the same is true of the globalised field $\gamma$ (see proof of
\GComp).  Hence the (time-dependent) generating field $\gamma^*$ 
for the modified global flow (see remark \TDvf) is initially 
small near $H$.  More
precisely, suppose that we restrict $\gamma$ to a neighbourhood
$W$ of $H$ in $Q\times\re$ and phase it out to be vertical 
outside a slightly larger neighbourhood $W'$, then the isotopy
defined by the corresponding
modified global flow is $C^\infty$--small.  Indeed by choosing $W$, 
$W'$ and $\mu$ (to define the upwards rotation) 
sufficiently small, we can make this isotopy arbitrarily 
small.

The local move will eliminate the downset
except in a small neighbourhood of the horizontal set.  Moreover the 
upwards rotated field $\beta$ changes, if at all, to be more upright.
Thus the global move, which is defined by the vector field left
after the local move, has precisely the form of the flow described 
in the observation and will be arbitrarily $C^\infty$--small.

For the local move we proceed as follows.  Apply \HSet\ and \DownMan\
so that $D$ is a manifold with boundary $H$.  Now choose a small
neighbouhood $U$ of $H$ in $M$ containing a smaller neighbourhood $U'$
say.  Apply \GPos\ to place $\overline{D-U'}$ in general position with
respect to the gradient flow $\varphi$ and choose $\delta>0$ small
enough that $U'$ has distance $>\delta$ from $\overline{M-U}$ and let
$V$ be the neighbourhood of $\overline{D-U'}$ corresponding to
$\delta$ given by \GPos.  We can assume that $\delta$ is also small
compared to the scale of $\alpha$, ie, compared with the distance over
which the direction of $\alpha$ changes.  Localise $\alpha$ using
$U\cup V$.

We now make the constructions made in the proof of \GComp\ and 
addendum (i).  Choose a suitably small $\mu>0$  and apply 
$\mu$--upwards rotation to give the vector field $\beta$.  
Choose a $\nu$--tubular neighbourhood $N$ of $M$ in $Q\times\re$
and globalise $\beta$ to $\gamma$ and let
$\gamma^*$ be the (time-dependent) vector field which gives the
corresponding modified global flow.  We shall make a final 
modification to this vector field.  We need to observe that 
vertical projection near $\overline{D-U'}$ is an immersion in $Q$ 
so we can choose a real number
$\omega$ such that any two points, $x\in\overline{D-U'}$ and 
$y\in M$, which project to the
same point of $Q$ are a distance at least $3\omega$ apart.
(Note that $\omega$ depends on $U'$.)

The final modification to the
flow is to multiply the time-dependent vector field 
$\gamma^*(x,t)$, which generates the flow, by $\rho(t)$ where 
$\rho$ is a bump function with
value $1$ for $t\le\omega$ and $0$ for $t\ge2\omega$.  (The purpose
of this final modification is to terminate
the effect of the flow once the straightening near $D-U'$ has occurred.)
We shall call this final modification {\sl phasing out} the flow.
This modified flow determines the local move.

We now prove that, for any sufficiently small choice of $U'$
(which then defines $\omega$),  
provided the other parameters and neighbourhoods
used to define the flow (namely $\mu$, $\delta$, $V$ and $\nu$) are 
chosen appropriately then the local move will move each point by less than
$\varepsilon\over2$.  We have already observed that the flow is
initially small inside $N\vert U$ (indeed smoothly small) and 
moreover the flow
is initially stationary outside $N\verts V\cup U$ and as we shall see
it is always $C^\infty$--small in these places.  We need to examine the 
flow in $N\vert V$.

In what follows we shall be near $M-H$ where vertical projection 
is an immersion, so we can use local coordinates of three types:
coordinates in $M$, horizontal coordinates perpendicular to $M$, which
we will call {\sl sideways} coordinates and one vertical coordinate.
Without loss we can assume that the fibres of $N$ over $V$ are
compatible with these local coordinates.  This implies that the 
foliation of $N$ given by restricting
$N$ to flow lines of $\varphi$ is locally invariant under sideways and
vertical translations.  We call this foliation $\Upsilon$.
We can also assume without loss that the rotation to vertical which
defines the global field $\gamma$ in terms of $\beta$ is described
in local coordinates of this type.

We now consider the effect of the local flow near $D$ and for
small time.  The key observation is the following.  On $D$ the vector field 
$\alpha$ lies in the same vertical plane as $\varphi$.  This
is just because $\alpha$ is perpendicular to $M$.  Since $\beta$
is obtained from $\alpha$ by rotation in this plane, the same
is true of $\beta$.  Now off $D$, $\alpha$ has a component 
perpendicular to $M$ in the vertical plane of $\varphi$ plus
a sideways component and this is still true after upwards rotation.
By choice of $N$ these facts imply that $\beta$
is tangent to leaves of $\Upsilon$ and by definition of the
global field, the same is true of $\gamma$ and hence initially
of $\gamma^*$. It follows that flowlines 
of the modified global flow lie in leaves of $\Upsilon$ or 
vertical translates.   

Now consider a vector $v$ of $\gamma$ and suppose that it has a
vertically upwards component of less than $\sqrt2\over2$ (ie, it makes
an angle $\pi\over4$ or less above the horizontal).  Then, since it is
obtained by upwards rotation of $\pi/2 -\mu$ or more from the
corresponding vector of $\alpha$, and this vector has a sideways
component and a component under $\varphi$ at most $\pi/2 -\mu$ below
horizontal, the combined sideways component and component under
$\varphi$ of $v$ must have size at least $\sqrt2\over2$.  Rotating
back upwards and subtracting the unit vertical vector to get the
corresponding vector of $\gamma^*$ we can see the following:

{\bf Note}\qua A vector of $\gamma^*$ {\it either} has a vertically
downwards component of at most $1-{\sqrt2\over2}$ {\it or} has a
combined sideways component and component under $\varphi$ of at least
$\sqrt2\over2$.

We are now ready to estimate the distance that a point moves under the
local flow.  It will help to think of the flow generated by $\gamma^*$
on a leaf of $\Upsilon$ in the following terms.  Let $\eta$ be a
flowline of $\varphi$ in $V$ (therefore of length $<\delta$).  Think
of the leaf $L(\eta)$ determined by $\eta$ and its vertical translates
as a river with vertical sides.  The flow is stationary on the banks
(and at many times and places in the middle as well).  The flow is
generated by a disturbance ($\gamma^*$ at time zero on $L(\eta)$)
which moves with unit speed downwards.  The total height $h$ of this
disturbance is controlled by the size of the leaf (determined by
$\delta$ and $N$) and $\nu$, the size of $N$.  Furthermore the
horizontal size of $L(\eta)$ is also controlled by $\delta$ and $\nu$.
Finally notice that the choice of $\delta$ small compared to the scale
of $\alpha$ implies that the direction of any sideways movement is
roughly constant on $L(\eta)$.  It now follows from the note made
above that each point of $\eta$ either reaches the banks of the river
or a point where the disturbance has passed (which then remains
stationary) after a time at most $\sqrt2(h+\delta+\nu)$, which we can make
as small as we please by choice of $\delta$ and $N$.  Moreover since
the speed a point moves is $<\sqrt2$ we may assume that each point is
moved at most $\ep\over2$.

We can also understand the effect of phasing out in terms of the same
river picture.  Suppose $\eta$ is outside $U'$ then the river meets
$M$ again at a distance approximately $3\omega$ from $\eta$; the
approximation is because meetings with $M$ are not necessarily
horizontal---but by choosing the river to be narrow enough (ie
$\delta$ and $\nu$ small enough) this approximation is as good as we
please.  Thus phasing out has the effect of killing the disturbance
after it has done its work on $\eta$ and before it can affect any
other point of $M$.  Finally notice that by choice of $\delta$ any
point of $\overline{V-U}$ lies on such a flowline outside $U'$.

We have proved that, by choosing the parameters suitably, the local
move straightens the vector field on $\overline{V-U}$ 
in a move which moves points at most $\ep\over2$ and moreover the
field outside $U$ does not move any points outside $\overline{V-U}$.
It follows from the observation made near the outset
that $M-\overline{V-U}$ moves only a 
$C^\infty$--small amount. 

This completes the local move.  We observe that the field would at
this point be straight if we had not phased out the flow.  However the
only place where the field might not be straight is now near $H$,
where the phasing out might have stopped the global straightening
before it was finished.  Moreover we can see that the local move has
not adversely affected the upwards rotated field $\beta$ near $H$.
The only substantial effect on $\beta$ (before the flow) comes from
localisation of $\alpha$, which moves vectors in $\alpha$ generally
upwards and has the same effect on $\beta$.  (This point will be
discussed in more detail in the slightly more general setting of the
proof of the full immersion theorem in part III---see [\CompIII;
figure 1] and adjacent text for more detail here.)  The flow also
moves vectors upwards.  Thus we are now in the situation described in
the observation near the start and the field can be finally
straightened by a $C^\infty$--small global move.  By choosing $U$ and
$\mu$ small enough, this final move also moves points less than
$\ep\over2$. \qed

\rk{Remark}In the proof, we could have straightened $\alpha$ in
one $C^0$--small move by not phasing out the local flow.  We could
then have estimated the total effect of this flow: each point is
affected only a finite number of times by the flow near $D$.  However
this would have obscured the local nature of the proof.  Indeed we can
see a local picture generated by the proof which depends on the
behaviour of $\varphi$ near $D$.  In particular it depends on the
intersection multiplicity.  For example the picture given in figure
\xtwist\ corresponds to a point of threefold intersection
multiplicity.  It is possible to describe the whole proof
combinatorially in terms of such local pictures.  The first sequence
of pictures described in section 2 is sufficient for the codimension 2
case.  In this description we would stratify $D$ according to the
intersection multiplicity with $\varphi$ (see [\Boa, \Mat]) and then
construct the straightening isotopy locally, starting at points of
highest multiplicity, and working down the stratification, using the
pictures constructed in section 2 as `templates'.

\sh{Straightening multiple vector fields}

Now let $M$ be embedded in $Q\times\re^n$ and suppose that    
$M$ is equipped with $n$ linearly independent normal vector fields.    
We say that $M$ is {\sl parallel} if the $n$ vector fields are    
parallel to the $n$ coordinate directions in $\re^n$.    
    
\numkey\MCT
\proc{Corollary : Multi-compression Theorem}Suppose that $M^m$ in    
embedded in $Q^q\times\re^n$ with    
$n$ independent normal vector fields and that $q-m\ge1$.    
Then $M$ is isotopic (by a $C^0$--small isotopy) to a parallel embedding.    
    
\prf  Apply the compression theorem to the first normal vector field.    
By the theorem $M$ can be isotoped to make this    
field parallel to the first coordinate axis.  Then $M$ lies    
over an immersion in $Q\times\re^{n-1}$ with $n-1$    
independent normal fields.  (This can be seen by thinking of
the remaining $n-1$ fields as determining a embedding of
$M\times D^{n-1}$ in $Q\times\re^n$ which is compressed into 
$Q\times\re^{n-1}$ by the straightening of the first field.)
Now consider an induced neighbourhood    
of $M$ pulled back by the immersion in $Q\times\re^{n-1}$.  
This neighbourhood is made    
of patches of $Q\times\re^{n-1}$ glued together and we can apply    
the local compression theorem to    
isotope $M$ within this neighbourhood until the second vector field    
is parallel to the second coordinate axis.  This isotopy determines    
a regular homotopy within $Q\times\re^{n-1}$ which lifts to an     
isotopy of $M$ in $Q\times\re^n$ finishing with an embedding     
which has the first two normal fields    
parallel to the first two axes of $\re^n$.  Continue    
in this way until all vector fields are parallel. \qed

\rk{Addenda}    

The local compression theorem admits a number of extensions.
Notice that these addenda improve considerably on the addenda
to theorem \GComp\ with the exception of addendum (iii) to
\GComp\ for which the analogue can be seen to be false:
If a grounded vector field on $M$ in $Q\times\re$ could be compressed 
by a {\it small\/} isotopy when $q-m=0$ then by working in an induced
neighbourhood (as in the proof above) a grounded {\it immersion} could
be compressed.  But it is easy to construct an immersion of
$S^1$ in $\re^2$ with grounded perpendicular vector field,
see figure \figkey\Eight.  This immersion cannot be 
compressed since $S^1$ does not immerse in $\re^1$.
    
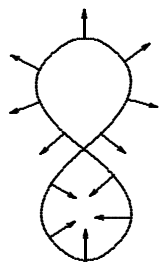
\fig{\Eight: a grounded field on an immersed $S^1$}
\beginpicture
\setcoordinatesystem units <0.350000cm,0.350000cm>
\unitlength=1.00000cm
\linethickness=1pt
\setplotsymbol ({\makebox(0,0)[l]{\tencirc\symbol{'160}}})
\setshadesymbol ({\thinlinefont .})
\setlinear
%
%
\linethickness= 0.500pt
\setplotsymbol ({\thinlinefont .})
\putrule from  7.243 22.703 to  7.243 23.688
%
%
\plot  7.307 23.434  7.243 23.688  7.180 23.434 /
%
%
%
\linethickness= 0.500pt
\setplotsymbol ({\thinlinefont .})
\plot  8.640 21.895  9.500 22.513 /
%
%
\plot  9.330 22.313  9.500 22.513  9.256 22.416 /
%
%
%
\linethickness= 0.500pt
\setplotsymbol ({\thinlinefont .})
\plot  5.704 21.609  4.720 22.117 /
%
%
\plot  4.975 22.057  4.720 22.117  4.917 21.944 /
%
%
%
\linethickness= 0.500pt
\setplotsymbol ({\thinlinefont .})
\plot  8.752 20.593  9.768 20.322 /
%
%
\plot  9.507 20.326  9.768 20.322  9.539 20.449 /
%
%
%
\linethickness= 0.500pt
\setplotsymbol ({\thinlinefont .})
\plot  7.832 19.387  8.657 18.737 /
%
%
\plot  8.418 18.844  8.657 18.737  8.497 18.944 /
%
%
%
\linethickness= 0.500pt
\setplotsymbol ({\thinlinefont .})
\plot  5.704 20.498  4.688 20.117 /
%
%
\plot  4.904 20.265  4.688 20.117  4.949 20.147 /
%
%
%
\linethickness= 0.500pt
\setplotsymbol ({\thinlinefont .})
\plot  6.562 19.418  5.736 18.705 /
%
%
\plot  5.887 18.919  5.736 18.705  5.970 18.823 /
%
%
%
\linethickness= 0.500pt
\setplotsymbol ({\thinlinefont .})
\plot  8.865 16.480  7.641 16.497 /
%
%
\plot  7.896 16.557  7.641 16.497  7.894 16.430 /
%
%
%
\linethickness= 0.500pt
\setplotsymbol ({\thinlinefont .})
\putrule from  7.292 15.037 to  7.292 16.053
%
%
\plot  7.355 15.799  7.292 16.053  7.228 15.799 /
%
%
%
\linethickness= 0.500pt
\setplotsymbol ({\thinlinefont .})
\plot  6.054 15.799  6.911 16.339 /
%
%
\plot  6.730 16.149  6.911 16.339  6.662 16.257 /
%
%
%
\linethickness= 0.500pt
\setplotsymbol ({\thinlinefont .})
\plot  6.085 17.655  6.943 17.164 /
%
%
\plot  6.691 17.235  6.943 17.164  6.754 17.345 /
%
%
%
\linethickness= 0.500pt
\setplotsymbol ({\thinlinefont .})
\plot  8.276 17.958  7.451 17.242 /
%
%
\plot  7.601 17.457  7.451 17.242  7.684 17.361 /
\linethickness= 0.500pt
\setplotsymbol ({\thinlinefont .})
%
%
%
\plot	 5.999 22.038 
 	 5.923 21.955
	 5.856 21.871
	 5.795 21.786
	 5.742 21.700
	 5.696 21.613
	 5.657 21.525
	 5.626 21.436
	 5.601 21.346
	 5.585 21.255
	 5.575 21.163
	 5.573 21.070
	 5.578 20.976
	 5.591 20.881
	 5.610 20.784
	 5.637 20.687
	 5.672 20.589
	 5.713 20.490
	 5.762 20.390
	 5.818 20.288
	 5.882 20.186
	 5.953 20.083
	 6.031 19.978
	 6.116 19.873
	 6.209 19.767
	 6.309 19.659
	 6.416 19.551
	 6.531 19.441
	 6.653 19.331
	 6.717 19.275
	 6.782 19.219
	 6.849 19.163
	 6.919 19.107
	 6.990 19.050
	 7.062 18.993
	 7.137 18.936
	 7.214 18.879
	 7.290 18.821
	 7.365 18.764
	 7.438 18.707
	 7.509 18.650
	 7.578 18.594
	 7.646 18.538
	 7.711 18.482
	 7.775 18.426
	 7.898 18.315
	 8.013 18.205
	 8.122 18.096
	 8.223 17.989
	 8.316 17.882
	 8.403 17.776
	 8.483 17.671
	 8.555 17.567
	 8.620 17.464
	 8.678 17.361
	 8.729 17.260
	 8.773 17.160
	 8.810 17.061
	 8.839 16.963
	 8.861 16.865
	 8.876 16.769
	 8.884 16.674
	 8.885 16.579
	 8.879 16.486
	 8.865 16.393
	 8.844 16.302
	 8.816 16.211
	 8.781 16.121
	 8.739 16.033
	 8.690 15.945
	 8.633 15.858
	 8.569 15.772
	 8.498 15.688
	 8.424 15.606
	 8.350 15.529
	 8.275 15.458
	 8.201 15.392
	 8.127 15.331
	 8.052 15.276
	 7.978 15.226
	 7.904 15.181
	 7.830 15.141
	 7.756 15.106
	 7.682 15.077
	 7.608 15.053
	 7.534 15.034
	 7.461 15.020
	 7.387 15.012
	 7.313 15.009
	 7.239 15.011
	 7.166 15.019
	 7.092 15.031
	 7.019 15.049
	 6.945 15.072
	 6.872 15.100
	 6.798 15.134
	 6.725 15.173
	 6.652 15.217
	 6.579 15.266
	 6.506 15.321
	 6.432 15.380
	 6.359 15.445
	 6.286 15.516
	 6.213 15.591
	 6.140 15.672
	 6.071 15.756
	 6.008 15.841
	 5.953 15.927
	 5.904 16.015
	 5.862 16.103
	 5.826 16.193
	 5.798 16.284
	 5.777 16.376
	 5.762 16.470
	 5.754 16.564
	 5.753 16.660
	 5.759 16.757
	 5.772 16.855
	 5.791 16.954
	 5.818 17.055
	 5.851 17.156
	 5.891 17.259
	 5.938 17.363
	 5.992 17.468
	 6.053 17.574
	 6.121 17.682
	 6.195 17.790
	 6.276 17.900
	 6.364 18.011
	 6.459 18.124
	 6.561 18.237
	 6.670 18.352
	 6.786 18.467
	 6.908 18.584
	 6.972 18.643
	 7.037 18.702
	 7.104 18.762
	 7.173 18.822
	 7.244 18.882
	 7.316 18.942
	 7.389 19.002
	 7.459 19.062
	 7.528 19.122
	 7.595 19.181
	 7.661 19.240
	 7.724 19.299
	 7.846 19.415
	 7.962 19.530
	 8.070 19.644
	 8.171 19.756
	 8.266 19.867
	 8.354 19.976
	 8.434 20.084
	 8.508 20.191
	 8.575 20.296
	 8.635 20.400
	 8.688 20.503
	 8.734 20.604
	 8.773 20.704
	 8.805 20.803
	 8.831 20.900
	 8.849 20.996
	 8.861 21.090
	 8.865 21.183
	 8.863 21.275
	 8.854 21.365
	 8.837 21.454
	 8.814 21.542
	 8.784 21.628
	 8.748 21.713
	 8.704 21.797
	 8.653 21.879
	 8.595 21.960
	 8.531 22.039
	 8.459 22.117
	 8.384 22.192
	 8.309 22.262
	 8.234 22.327
	 8.159 22.387
	 8.083 22.442
	 8.007 22.492
	 7.932 22.536
	 7.856 22.576
	 7.780 22.611
	 7.704 22.641
	 7.628 22.666
	 7.551 22.686
	 7.475 22.701
	 7.398 22.711
	 7.321 22.716
	 7.245 22.716
	 7.168 22.711
	 7.091 22.701
	 7.013 22.686
	 6.936 22.666
	 6.859 22.642
	 6.781 22.612
	 6.703 22.577
	 6.625 22.537
	 6.548 22.492
	 6.470 22.442
	 6.391 22.387
	 6.313 22.327
	 6.235 22.262
	 6.156 22.192
	 6.077 22.118
	 5.999 22.038
	/
\linethickness=0pt
\putrectangle corners at  4.663 23.713 and  9.794 14.988
\endpicture
\endfig

There are similar addenda to the multi-compression theorem which    
we leave the reader to state and prove.\sl     
    
\items
\item{\rm(i)}{\rm(Non-compact relative version)}    
The local compression theorem is true without the    
hypothesis that $M$ is compact, moreover there is a far stronger    
relative version.     
Let $C$ be a closed set in $M$. If $M$ is already compressible in a    
neighbourhood of $C$ then the isotopy can be    
assumed fixed on $C$.  Moreover we can replace the real number 
$\varepsilon$ in the
statement of the theorem by a function $\varepsilon$ from $M$ to 
the positive reals.  The conclusion is that each point $x$ moves   
a distance less than $\ep(x)$.
    
\item{\rm(ii)}If each component of $M$ has relative boundary then the    
dimension condition can be relaxed to $q-m\ge0$.%
{\rm\qua(See note 4.6 below.)}  
    
\item{\rm(iii)}{\rm (Parametrised version)}
Given a parametrised family $M^m_t\subset Q^q_t$    
of embeddings with a normal vector field, $t\in K$,    
where $K$ is a manifold of dimension $k$ and    
$q-m-k\ge 1$, then there is a parametrised family of    
small isotopies to compressible embeddings (and there is a relative version    
similar to (i)).{\rm\qua(Note that the Normal Deformation Theorem 4.7 below
admits a parametrised version with no additional dimensional condition.)}    
    
\item{\rm(iv)}If in (iii) the fields are all perpendicular and
grounded, then the    
dimension condition can be relaxed to $q-m\ge 1$ and there is no
dimension condition on $K$.

\item{\rm(v)}If each component of $M$ has relative boundary then the    
dimension conditions in (iii) and (iv) can each be relaxed by 1.    
    
\item{\rm(vi)}All the above statements apply to immersions of $M$ in
$Q$ (or families of immersions).  The result is a small regular
homotopy (or family) to a compressible immersion (or family).
\enditems

\rm    
    
\rk{Proofs}The relative version in addendum (i), with $M$
compact, follows from the method of proof:  notice that the
first (local) move in the proof of \LCT\ can be assumed fixed
near $C$ and the effect of the second ($C^\infty$--small) move 
near $C$ is to move $M$ a little through compressible embeddings.
This movement can be cancelled by a small local isotopy. 
For the non-compact and variable $\ep$ versions, use the 
relative version and a patch by patch argument.
To prove (ii) use a handle decomposition of $M$ with no     
$m$--handles, and inductively apply the proof to the core of     
each handle, working relative to the union of the previous handles.
Notice that once the core of a handle has been straightened then 
a small neighbourhood is also straight.  But we can shrink the
handle into such a neighbourhood.
For (iii) we notice that the vector field can be assumed to be
grounded by the dimension condition (see corollary A.5)
and the two transversality results used
(\HSet\ and \DownMan) both have parametrised versions (proved
in the appendix, as 
propositions A.2 and A.3) which state that the union of the horizontal
sets and downsets in each fibre can be assumed to be manifolds. We
can then define localisation near the downset in the same way as 
in the unparametrised version and move $D$ in into general position
with respect to the flow on $K\times M$ which is $\varphi$ in each
fibre, see corollary A.6.  The proof then proceeds as before.
Addendum (iv) follows since the dimension condition was used only to
establish groundedness and for (v) we use a handle argument 
as in (ii).   Finally for (vi) we work on an induced regular
neighbourhood as in the proof of \MCT. \qed    
    
\proc{Note}\rm Strictly speaking, the isotopy in part (ii) in the
case $q=m$ is not small.  It is of the form: shrink to the
neighbouhood of a spine and then perform a small isotopy.  Thus it
small in the sense that it takes place inside an arbitrarily small
neighbourhood of the initial embedding.  But notice that this is
precisely what is needed for the proof of the multi-compression
theorem.

\sh{Bundle and deformation versions}

The relative version of the multi-compression theorem leads at once to
a bundle version: suppose that instead of $M^m\subset Q\times\re^n$, we
have $M$ contained in the total space $W$ of a bundle $\eta^n/Q$ such
that there is a subbundle $\xi$ of $TW$ defined at $M$ and isomorphic
to the pull back of $\eta$ to $M$.  Then there is an isotopy of $M$
and $\xi$ realising this isomorphism.  The proof is to work locally
where $\eta$ is trivial and apply the multi-compression theorem to
straighten $\xi$ (ie realise the isomorphism with $\eta$).  This can
be further generalised by replacing $\eta$ by a second subbundle of
$TW$.  It is convenient to state this result as a deformation theorem,
which then admits a parametrised version with no extra conditions on
dimension:

\proc{Normal Deformation Theorem}
Suppose that $M^m\subset W^w$ and that $\xi^n$ is a subbundle of $TW$
defined in a neighbourhood $U$ of $M$ such that $\xi|M$ is normal to
$M$ in $W$ and that $m+n<w$.  Suppose given $\ep >0$ and a homotopy of
$\xi$ through subbundles of $TW$ defined on $U$ finishing with the
subbundle $\xi'$.  Then there is an isotopy of $M$ in $W$ which moves
points at most $\ep$ moving $M$ to $M'$ and covered by a bundle
homotopy of $\xi|M$ to $\xi'|M'$ (and in particular $\xi'|M'$ is
normal to $M'$).
\endproc

\prf  
We work in small compact patches where we can assume that $\xi$ and
$\xi'$ are trivial and we use the proof of the multi-compression
theorem.  Thinking of $\xi$ as comprising $n$ linearly independent
vector fields, consider the first vector field $\alpha$.  By
compactness the total angle that the homotopy of $\xi$ moves $\alpha$
is bounded and we can choose $r$ and a sequence of homotopies
$\xi=\xi_0\simeq\xi_1\simeq \ldots\simeq \xi_r=\xi'$ so that for each
$s=1,\ldots,r-1$ the images of $\alpha$ in $\xi_{s-1}$ move through an
angle less than $\pi/2$ in the homotopy to $\xi_s$.  Now apply the
local proof of the compression theorem.  Think of the image $\alpha_1$
of $\alpha$ in $\xi_1$ as vertically up, then if $\alpha$ is made
perpendicular it is grounded and the proof gives a small isotopy of
$M$ covered by a straightening of $\alpha$ which moves $\alpha$ by a
homotopy which can be seen to be a deformation of the given homotopy
(both homotopies take place in the contractible neighbourhood of
$\alpha_1$ comprising vectors which make an angle $<\pi$ with
$\alpha_1$).  Repeating this $r$ times we find a small isotopy of
$M,\alpha$ to $M',\alpha'$ where $\alpha'$ is the first vector of
$\xi'$.  This moves $\xi$ to $\xi''$ say and then if we consider the
homotopy $h$ of $\xi''$ to $\xi'$ which is the reverse of the isotopy
followed by the given homotopy of $\xi$ to $\xi'$ then $h$ moves
$\alpha'$ through a contractible loop and contracting this loop we
obtain a homotopy of $\xi''$ to $\xi'$ fixing $\alpha'$.  We now
project onto the orthogonal complement of $\alpha'$ (as in the proof
of the mulit-compression theorem) and proceed to straighten the next
vector in the same way, again following the given homotopy. \endprf

\rk{Addenda} The covering bundle homotopy is a deformation of the
given one. There is a relative version which follows directly from the
proof and, since all fields are grounded, there is a parametrised
version with no extra hypotheses, see proof of addenda (iii) and (iv)
to the multi-compression theorem given above.

The proof can readily be modified to construct an isotopy which
``follows'' the given homotopy of $\xi$, in other words if $\xi_i$ is
the position of $\xi$ at time $i$ in the homotopy and $M_i$ the
position of $M$ at time $i$ in the isotopy, then $\xi_i|M_i$ is normal
to $M_i$ for each $i$.  This is done by breaking the homotopy into
very small steps.  However this makes the final isotopy far less
explicit and in that case there is a simpler proof given in Part II
[\CompII; section 3].

Finally there is a codimension 0 ($m+n=w$) version which it is worth
spelling out in detail:

\sl Suppose in the normal deformation theorem that $m+n=w$ and that $M$ 
is open or has boundary and that $X$ is a spine of $M$.  Then there is
an isotopy of $M$ of the form: shrink into a neighbourhood of $X$
followed by a small isotopy in $W$, carrying $M$ to be normal to
$\xi'$.\rm

\rk{Remark}The Normal Deformation Theorem is close to Gromov's theorem 
on directed embeddings [\Gr; 2.4.5 $\rm C'$] from which it can be
readily deduced.  For more detail here see the final remarks in
[\CompII].

\references

\def\section#1{\vskip-\lastskip\penalty-800\vskip 20pt plus10pt minus5pt 
{\large\bf A\quad#1}         
\vskip 8pt plus4pt minus4pt
\nobreak\resultnumber=1}
\def\proc#1{\vskip-\lastskip\ppar\bf%
\noindent#1\ A.\number\resultnumber
\stdspace\sl\global\advance\resultnumber by 1\ignorespaces}

\section{Appendix: transversality and general position}

\small
In this appendix we prove the transversality and general position
results needed in the main proofs (in sections 2 and 4).

Let $M^m$ and $K^k$ be manifolds and let $F\co K\times M\to 
\re^n$ be a {\sl $K$--family of embeddings} of $M$ in $\re^n$. 
Ie, $F$ is a smooth map such that $F\verts \{t\}\times M$ is
a smooth embedding for each $t\in K$.  We denote
$F(\{t\}\times M)$ by $M_t$.  Let $E_n$ 
be the group of isometries of $\re^n$. Define $J\co 
E_n\times K\times M\to\re^n$ by  $J(u,t,x)=uF(t,x)$. Let 
$\G_{n,m}$ be the Grassmannian of $m$-planes in $\re^n$.
For $u\in E_n$ let $F_u\co K\times M\to\G_{n,m}$ be given by 
$F_u(t,x)=T_{J(u,t,x)}J(\{u\}\times\{t\}\times M)$.

\proc{Proposition} Given a neighbourhood $N$ of the identity of 
$E_n$ and a submanifold $W$ of $\G_{n,m}$ there is an 
element $u\in N$ such that  $F_u\co K\times M\to\G_{n,m}$ is 
transverse to $W$.

\prf Let $G\co E_n\times K\times M\to \G_{n,m}$ be given by 
$G(u,t,x)=F_u(t,x)$. It is sufficient to prove the derivative 
$TG_{(u,t,x)}$ is surjective for any $(u,t,x)$, see [\GP, 
page 68]. In fact we prove a stronger result: the derivative on 
tangents to each $E_n\times\{t\}\times\{x\}$ at any $(u,t,x)$ is 
surjective.

Represent a tangent vector to the Grassmannian by 
$\dot\alpha(0)$, where $\alpha\co (-\epsilon,\epsilon)\to 
\G_{n,m}$
and $\alpha(0)=G(u,t,x)$. Let $q\co O_n\to \G_{n,m}$ be the
fibration given by $q(v)=v(\alpha(0))$ and choose $\beta\co(-
\epsilon,\epsilon)\to O_n$ such that  $q\beta=\alpha$.  Finally
define $\gamma\co(-\epsilon,\epsilon)\to E_n$ by 
$\gamma(s)=\tr_{J(u,t,x)}\circ\beta(s)\circ \tr_{-J(u,t,x)}$ where 
$\tr_z$
denotes translation by $z$. Then after identifying $E_n$ with 
$E_n\times\{p\}\times\{x\}$ we have $G\gamma=\alpha$ and 
so 
$TG(\dot\gamma(0))=\dot\alpha(0) $ as required.\qed

Suppose now $Q^q$ is a manifold and we are given a family of 
embeddings $F\co K\times M\to Q\times \re^r$ with $m\leq q$.
Define $H(K\times M)$, the {\sl horizontal set}, to consist of 
those points $(t,x)\in K\times M$ such that the tangents to 
$\re^r$ at $F(t,x)$ are contained in the image of $TF$. 
Similarly define $C(K\times M)$, the {\sl critical set}, to 
consist of points $(t,x)\in K\times M$ such that the image of 
tangents at  $(t,x)$ are tangent to $Q$. In case $Q=\re^q$, by 
considering $\G_{q,m-r}\subset\G_{q+r,m}$ and 
$\G_{q,m}\subset\G_{q+r,m}$, we can assume by, A.1, that 
after a small Euclidean motion $H(K\times M)$ and 
$C(K\times M)$ are submanifolds of codimension $r(q+r-m)$ 
and $rm$ respectively.

The terminology has been chosen so that if $M_t$ 
has a perpendicular normal field then the field must be horizontal 
at points $F(t,x)$ where $(t,x)\in H(K\times M)$, and if 
$(t,x)\in C(K\times M)$ and $r=1$, then $x$ is a critical point 
of the function $M\to\re$ given by $x\mapsto \pi_2F(t,x)$, 
where $\pi_2\co Q\times\re\to \re$ is projection.

\proc{Proposition} Given a family of embeddings 
$F\co K\times M^m\to Q^q\times\re^r$ with $m\leq q$, then 
there is an isotopy of $F$, ie
a homotopy $F_t$ through families of embeddings, such that 
$F=F_0$ and the horizontal set determined by $F_1$ is a 
submanifold of $K\times M$.  Similarly for the critical set.

\prf We prove the case of the horizontal set.  The critical set
is similar.  The case $Q=\re^q$ follows from A.1. 
For the general case we use a standard patch by patch argument.
Choose a locally 
finite cover of $K\times M$ by discs of the form $D=D_1\times D_2$
such that the corresponding half discs 
${1\over2}D={1\over2}D_1\times {1\over2}D_2$ cover and such that 
each disc has image above a Euclidean 
patch in $Q$.  Suppose then $F(D)\subset
U\times\re^r$ where $U$ is a Euclidean patch. By A.1 we can compose 
$F\vert D$ with a small Euclidean motion $e$ of $U\times\re^r$ so that the 
horizontal set in a neighbourhood of ${1\over2}D$
becomes a manifold.  By choosing a path from $e$ to the 
identity in $E_{q+r}$ we can phase out the movement using a 
collar on ${1\over2}D$ in $D$ and then extend by the identity to
$K\times M$. 

We now argue by induction.  Suppose the 
horizontal set in a neighbourhood of the union of the
first $k$ such $1\over2$--discs is a manifold and 
consider the $(k+1)$-st.  By choosing $e$ as above with a sufficiently
short path to the identity, the horizontal set in a 
neighbourhood $(k+1)$-st $1\over2$--disc
becomes a manifold without disturbing the property
that the horizontal set is a manifold in (a possibly smaller)
neighbourhood of the first $k$ $1\over2$--discs.  
Since the cover is locally finite, this inductively defined
isotopy defines the required isotopy of $F$. \qed

We now restrict to the case of application in sections 2 and 4
namely $r=1$ and suppose that $K\times M\to Q\times\re$ is a
$K$--family of embeddings equipped with normal vector fields.
Ie, we assume that the normal vector fields $\alpha$ vary continuously
with $t\in K$.  We identify $K\times M$ with its image in 
$K\times Q\times\re$ by the embedding $(t,x)\mapsto (t,F(t,x))$ and
we think of $\alpha$ as a single vector field on
$K\times M$ in $K\times Q\times\re$.  We denote the vector field 
on $M_t$ in $Q$ corresponding to $t\in K$ by $\alpha_t$.

Choose metrics and assume that $\alpha$ is a unit vector field
such that $\alpha_t$ is perpendicular to $M_t$ in $Q\times\re$ for each
$t\in K$.  We define the {\sl downset} exactly as in section 4.  
Let $\nu$ be the normal bundle on $K\times M$
in $K\times Q\times\re$ given as the union of the normal bundles 
$\nu(M_t\subset Q\times\re)$ for $t\in K$.  Let $\psi$ be the
vector field on $K\times M-H(K\times M)$ in $K\times Q\times\re$ given
by choosing $\psi_t(x)$ to point up the line of steepest ascent 
in $\nu_{(t,x)}$.  Then the downset $D\subset K\times 
M-H(K\times M)$ comprises points $(t,x)$ where 
$\alpha_t(x)=-\psi_t(x)$ ie, $D$ comprises 
all points $(t,x)$ such that $\alpha_t$ points down the line of
steepest descent in $\nu_{(t,x)}$.  Notice that the
downset of $\alpha$ is precisely the union of the downsets of
$\alpha_t$ over $t\in K$.

\proc{Proposition} After a small isotopy of $\alpha$ we can assume that    
$\alpha$ is transverse to $-\psi$, which implies that $D$ is a 
manifold, and further we can assume that the closure of $D$ is a 
manifold with boundary $H$.\rm

\prf By A.2 the horizontal set can be assumed to be a manifold $H$ say.
Since we are in the case $r=1$, the codimension of $H$ in
$K\times M$ is  $c=q+1-m$.  Now assume for simplicity that 
$Q=\re^q$.  Then $H$ is given (by proposition A.1) as the
transverse preimage of $\G_{q,m-1}\subset\G_{q+1,m}$, which
(using the metrics) can also be seen as the transverse preimage of 
$\G_{q,c}\subset\G_{q+1,c}$.

Now we can identify the restriction of the 
canonical disc bundle $\bar\gamma_{q+1,c}$ to $\G_{q,c}$ with a 
closed tubular neighbourhood of $\G_{q,c}$ by 
$(p,P)\mapsto P_p\in G_{q+1,c}$, where $P_p$ is the subspace spanned 
by the vector $(p,-|p|)\in\re^q\times \re^1=\re^{q+1}$ and    
the subspace of $P$ orthogonal to $p$.  
   
Now consider $\alpha$ near $H$. Consider a fibre $D_x$ of a tubular 
neighbourhood  $U$ of $H$ in $K\times M$. Let $\pi\co U\to H$ be the 
projection of $U$.  We can identify the fibres of $\nu$ over $D_x$ 
with $\nu_x$ by choosing an isomorphism of $\pi^*\nu$ with $\nu\vert U$. 
This done we can homotope $\alpha$ to be constant over $D_x$. 
Since $\nu$ is transverse to $\G_{q,c}$ we find that, if $U$ is chosen 
small enough, the map $-\psi$ on $\d D_x$ followed by projection
to the unit sphere is a diffeomorphism.  So $-\psi$ and    
$\alpha$ on $\d D_x$ meet in one point. We now get a collar on $H$ by    
considering variable $x$ and variable radius for $D_x$.
Consequently we have    
$\alpha$ transverse to $-\psi$ on $U-H$, and after a further    
homotopy of $\alpha$we have $\alpha$ transverse to $-\psi$ on    
$K\times M-H$. 
   
The general case  follows by a patch by patch argument as in A.2.    
\qed    

\proc{Lemma}We can realise any small isotopy of $D$ fixed near
$H$ by a small isotopy of $\alpha$.

\prf Consider a small patch in $D$.  Locally $D$ is given as
the tranvserse preimage of $-\psi$ by $\alpha$.  A small local
isotopy of $\alpha$ maintains transversality and has the effect
of moving $D$ inside its normal bundle in $K\times M$ to an
arbitrary section near the zero section.  By a combination of
such moves we can realise any small isotopy of $D$ fixed near
$H$ by a small isotopy of $\alpha$.  \qed

\proc{Corollary}Suppose that $q-m-k\ge1$ then by a small isotopy
of $F$ and $\alpha$ we may assume that $\alpha_t$ is perpendicular and
grounded for each $t\in K^k$.

\prf By A.2 we may assume that the critical set of $F$ is a manifold
of dimension $k$.  But recall that $\alpha_t$ is grounded if and only
if the downset of $\alpha_t$ is disjoint from the critical set of
$M_t$.  The result now follows from the lemma. \qed\ppar

Now let $\varphi_t$ be the gradient vector field on $M_t$ for each
$t\in K$ and let $\varphi$ be the vector field on $K\times M$
such that $\varphi\verts\{t\}\times M=\varphi_t$.

\proc{Corollary}Suppose that $\alpha$ is perpendicular and
grounded.  Let $U$ be a tubular neighbourhood of $H$ in $K\times M$.
By a small isotopy of $\alpha$ we may assume that 
$\overline{D-U}$ is in general position with respect to $\varphi$    
in the following sense:  Given $\delta>0$ there    
is a neighbourhood $V$ of $\overline{D-U}$ in $M$ such that each 
component of intersection of an integral curve of $\varphi$ with 
$V$ has length $<\delta$.
 
\prf  By the lemma and the theory of smooth maps of codimension 0
[\Boa, \Mat], we may assume that each integral curve of $\varphi$ meets 
$\overline{D-U}$ in a discrete set.  The result now follows for compact $K$
and $M$ by a standard accumulation argument.  For the non-compact
case, we use a patch by patch argument as in A.2. \qed 

\bye